\documentclass[12pt]{article}
\usepackage{latexsym,amssymb,amsmath}
\usepackage{stmaryrd}
\begin{document}
\baselineskip=16pt

\newcommand{\la}{\langle}
\newcommand{\ra}{\rangle}
\newcommand{\psp}{\vspace{0.4cm}}
\newcommand{\pse}{\vspace{0.2cm}}
\newcommand{\ptl}{\partial}
\newcommand{\dlt}{\delta}
\newcommand{\sgm}{\sigma}
\newcommand{\al}{\alpha}
\newcommand{\be}{\beta}
\newcommand{\G}{\Gamma}
\newcommand{\gm}{\gamma}
\newcommand{\vs}{\varsigma}
\newcommand{\Lmd}{\Lambda}
\newcommand{\lmd}{\lambda}
\newcommand{\td}{\tilde}
\newcommand{\vf}{\varphi}
\newcommand{\yt}{Y^{\nu}}
\newcommand{\wt}{\mbox{wt}\:}
\newcommand{\rd}{\mbox{Res}}
\newcommand{\ad}{\mbox{ad}}
\newcommand{\stl}{\stackrel}
\newcommand{\ol}{\overline}
\newcommand{\ul}{\underline}
\newcommand{\es}{\epsilon}
\newcommand{\dmd}{\diamond}
\newcommand{\clt}{\clubsuit}
\newcommand{\vt}{\vartheta}
\newcommand{\ves}{\varepsilon}
\newcommand{\dg}{\dagger}
\newcommand{\tr}{\mbox{Tr}}
\newcommand{\ga}{{\cal G}({\cal A})}
\newcommand{\hga}{\hat{\cal G}({\cal A})}
\newcommand{\Edo}{\mbox{End}\:}
\newcommand{\for}{\mbox{for}}
\newcommand{\kn}{\mbox{ker}}
\newcommand{\Dlt}{\Delta}
\newcommand{\rad}{\mbox{Rad}}
\newcommand{\rta}{\rightarrow}
\newcommand{\mbb}{\mathbb}
\newcommand{\lra}{\Longrightarrow}
\newcommand{\X}{{\cal X}}
\newcommand{\Y}{{\cal Y}}
\newcommand{\Z}{{\cal Z}}
\newcommand{\U}{{\cal U}}
\newcommand{\V}{{\cal V}}
\newcommand{\W}{{\cal W}}
\setlength{\unitlength}{3pt}

\begin{center}{\Large \bf  A New Factor from $E_6$-{\bf Mod} to $E_7$-{\bf Mod} } 
\end{center}
\vspace{0.2cm}

\begin{center}{\large Xiaoping Xu}\end{center}
\begin{center}{Hua Loo-Keng Key Mathematical Laboratory}\end{center}
\begin{center}{Institute of Mathematics, Academy of Mathematics \& System Sciences}\end{center}
\begin{center}{Chinese Academy of Sciences, Beijing 100190, P.R. China}
\footnote{Research supported
 by China NSF 10871193}\end{center}

\begin{abstract}{ We find a new
representation of the simple Lie algebra of type $E_7$ on the
polynomial algebra in 27 variables, which gives a fractional
representation of the corresponding Lie group on 27-dimensional
space. Using this representation and Shen's idea of mixed product,
we construct a functor from $E_6$-{\bf Mod} to $E_7$-{\bf Mod}. A
condition for the functor to map a finite-dimensional irreducible
$E_6$-module to an infinite-dimensional irreducible $E_7$-module is
obtained. Our general frame also gives a direct polynomial extension
from irreducible $E_6$-modules to irreducible $E_7$-modules. The
obtained infinite-dimensional irreducible $E_7$-modules are $({\cal
G},K)$-modules in terms of  Lie group representations. The results
could be used in studying the quantum field theory with $E_7$
symmetry and symmetry of partial differential
equations.}\end{abstract}

\section{Introduction}

\quad \quad A quantum field is an operator-valued function on a
certain Hilbert space, which is often a direct sum of
infinite-dimensional irreducible modules of a certain Lie algebra
(group). The Lie algebra of two-dimensional conformal group is
exactly the Virasoro algebra, which contains a central element. The
minimal models of two-dimensional conformal field theory were
constructed from direct sums of certain infinite-dimensional
irreducible modules of the Virasoro algebra, where the value taken
by the central element is an important physical quantity, called the
{\it central charge}.

It is well known that $n$-dimensional projective group gives rise to
a non-homogenous representation of the Lie algebra $sl(n+1,\mbb{C})$
on the polynomial functions of the projective space. Using Shen's
mixed product for Witt algebras, Zhao and the author [ZX]
generalized the above representation of $sl(n+1,\mbb{C})$ to a
non-homogenous representation on the tensor space of any
finite-dimensional irreducible $gl(n,\mbb{C})$-module with the
polynomial space. Moreover, the structure of such a representation
was completely determined. The result can be used to study the
quantum field theory with $sl(n+1,\mbb{C})$ as the symmetry.
 Furthermore, we [XZ] generalize the
conformal representation of $o(n+2,\mbb{C})$ to a non-homogenous
representation of $o(n+2,\mbb{C})$ on the tensor space of any
finite-dimensional irreducible $o(n,\mbb{C})$-module with a
polynomial space by Shen's idea of mixed product for Witt algebras.
It turns out that  a hidden central transformation is involved. More
importantly,  we find a condition on the constant value taken by the
central transformation  such that the generalized conformal
representation is irreducible. The result would be useful in
higher-dimensional conformal field theory.

In [X4], we find a new representation of the simple Lie algebra of
type $E_6$ on the polynomial algebra in 16 variables, which gives a
fractional representation of the corresponding Lie group on
16-dimensional space. Using this representation and Shen's idea of
mixed product, we construct a functor from $D_5$-{\bf Mod} to
$E_6$-{\bf Mod}. A condition for the functor to map a
finite-dimensional irreducible $D_5$-module to an
infinite-dimensional irreducible $E_6$-module is obtained. Our
general frame also gives a direct polynomial extension from
irreducible $D_5$-modules to irreducible $E_6$-modules. The obtained
infinite-dimensional irreducible $E_6$-modules are $({\cal
G},K)$-modules in terms of  Lie group representations. The results
could be used in studying the quantum field theory with $E_6$
symmetry and symmetry of partial differential equations.

This paper is the forth work in the program of studying
quantum-field motivated representations of finite-dimensional simple
Lie algebras. It is well known that the minimal dimension of
irreducible modules over the  simple Lie algebra of type $E_7$ is
56. Based on a grading of the simple Lie algebra of type $E_7$, we
find a first-order differential operator representation of the Lie
algebra on the polynomial algebra in 27 independent variables.  In
fact, the corresponding Lie group representation is given by
fractional transformations on 27-dimensional space. Using this
representation and Shen's idea of mixed product, we construct a new
functor from $E_6$-{\bf Mod} to $E_7$-{\bf Mod}, where a hidden
central transformation is involved. More importantly,  a condition
for the functor to map a finite-dimensional irreducible $E_6$-module
to an infinite-dimensional irreducible $E_7$-module in terms of the
constant value taken by the central transformation is obtained. The
well-known Dickson cubic invariant (cf. [Dl]) plays an important
role in our approach.  Our general frame also gives a direct
polynomial extension from irreducible $E_6$-modules to irreducible
$E_7$-modules, which can be applied to obtain explicit bases of
irreducible $E_7$-modules from those of irreducible $E_6$-modules.
The result could be useful in understanding the quantum field theory
with $E_7$ symmetry. Our fractional representation of the $E_7$ Lie
group could also be  used in symmetry analysis of partial
differential equations just as the conformal representation of
orthogonal Lie groups does. Our infinite-dimensional irreducible
$E_7$-modules are $({\cal G},K)$-modules in terms of the
corresponding Lie group representations.

The $E_7$ Lie algebra and group are important mathematical objects
with broad applications. They are second most complicated
finite-dimensional simple Lie algebra and group. The complexity
indeed imply rich connotation. Here we are only  able to list a
small part of it. Ramond [R] (1977) gave a group theoretical
analysis of a symmetry breaking affected by Higgs fields for the
vector-like unified theory base on $E_7$. Cvitanovi\'{e} [Co] (1981)
studied the $E_7$ symmetry with negative bosonic dimension.  Han,
Kim and Tanii [HKT] (1986) presented a supersymmetrization of the
six-dimensional anomaly-free $E_6\times E_7\times U(1)$ theory with
Lorentz Chern-Simons term. Kato and Kitazawa [KK] (1989) studied the
correlation functions in an $E_7$-type modular-invariant
Wess-Zumino-Witten theory, which is related to the scheme of string
compatification proposed by Gepner. Ma [M] (1990) found the
spectrum-dependent solutions to the Yang-Baxter equation for quantum
$E_7$.

Fern\'{a}dez,  Garcia Fuerres and  Perelomov [FFP] (2005)
re-expressed the quantum Calogero-Sutherland model for the Lie
algebra $E_7$ and the particular value of the coupling constant
$k=1$, using the fundamental irreducible characters of the algebra
as dynamic variables. Moreover, they used the model to obtain
explicitly the characters and Clebsch-Gordan series for the algebra
$E_7$. D'Auria, Ferrara and Trigiante [DFT] (2006) considered
M-theory compactified on seven-torus with fluxes when all the seven
anti-symmetric tensor fields in four dimensions have been dualized
into scalars and thus the $E_{7(7)}$ symmetry was recovered. Duff
and  Ferrara [DF] (2007) proposed that a particular tripartite
entanglement of seven quits, encoded in the Fano plane, is described
by the $E_7$ group and that the entanglement measure is given by
Cartan's quartic $E_7$ invariant. Brink,  Kim and  Ramond [BKR]
(2008) used the Cremmer-Julia $E_{7(7)}$ nonlinear symmetry of $N=8$
supergravity to derive its order $\kappa$ on-shell Hamltonian in
terms of the chiral light-cone superfield. Borsten [Bl] (2008) used
the Freudenthal triple system to sketch  the precise dictionary
relating the 56 charges, parametrizing the general black hole
solution, to the 56-dimensional quantum states of $E_7$.

Part of the mathematical story of $E_7$ is as follows. Brown [Br1]
(1968) constructed a central simple nonassociative algebraic
structure on the Freudenthal triple system of $E_7$ and obtained a
new interpretation of the $E_7$ Lie algebra. He [Br2] (1969) proved
that the group leaving Cartan's quartic form invariant modulo its
center is exactly a simple $E_7$ Lie group. Faulkner [Fau] (1972)
studied the geometry of the planes of the octave symplectic geometry
of the ternary algebra related to the minimal representation of
$E_7$. Ferrar [Fer] (1980) used the Freudenthal triple system  to
classified the simple $E_7$ Lie algebras over algebraic number
fields. Kleidman and Ryba [KR] (1993) proved a Kostant's conjecture
in the case of $E_7$, which says that the group $PSL(2,37)$ embeds
in an $E_7$ Lie group. Moreover, Griess and Ryba [GR] (1994) proved
that the group $U_3(8)$ embeds in an $E_7$ Lie group and the group
$Sz(8)$ does not. Ginzburg [Gd] (1995) showed that the twisted
partial $L$-function on the 56-dimensional representation of
$GE_7(\mbb{C})$ is entire except the points 0 and 1.

Cooperstein [Cb] (1995) proved that  Cartan's  quartic
$E_7$-invariant is the unique fundamental invariant over the basic
module. Moreover, Shult [Se] (1997) showed that the basic
representation of $E_7$ provides the absolutely universal embedding
of the point-line geometry $E_{7,1}$. Plotkin [P] (1998)
investigated the stability of the $K_1$-functor the $E_7$ group.
Dokovi\'{c} [Dd1] (1999) classified nilpotent adjoint orbits of real
simple noncompact groups of type $E_7$ by means of Caley triples. He
[Dd2] (2001) also worked out the partial order on nilpotent orbits
in the split real Lie algebra of type $E_7$. The $E_7$ root system
was used by Sekiguchi [Sj] (1999) to study the configurations of
seven lines on the real projective plane. Kleidman, Meierfrankenfeld
and Ryba [KMR] (1999) proved that the group $HS$ embeds in $E_7(5)$.
Choi and Yoon [CY] (1999) calculated the homology of the double and
triple loop spaces of $E_7$.

Ukai [U] (2001) computed the $b$-functions of the prehomogeneous
vector space arising from a cuspidal character sheaf of $E_7$.
Garibaidi [Gr] (2001) gave  explicit descriptions of the homogeneous
projective varieties associated with the $E_7$ group with trivial
Tits algebras. A characterization of quadratic forms of type $E_7$
was given by De Medits [Dt] (2002). Cherenousov [Cv] (2003) used
Hasse principle to prove that the Rost invariant has trivial kernel
for a quasi-split $E_7$ group. Kono, Lin and Nishimura [KLN] (2003)
characterized the mod 3 cohomology of $E_7$. Weiss [W] (2006) solved
a fundamental question about the structure of the automorphism group
of the Moufang quadrangles of type $E_7$. Moreover, we found in [X3]
that the weight matrices of $E_7$ on its minimal irreducible module
and adjoint module  generate ternary orthogonal codes with large
minimal distances.

This work further reveals new beauties of the simple Lie algebra of
type $E_7$. In Section 2, we construct the 27-dimensional
representation of the simple Lie algebra of type $E_6$ in terms of
first-order differential operators on the polynomial algebra in 27
independent variables from the lattice-construction of the simple
Lie algebra of type $E_7$.  In Section 3, we realize the simple Lie
algebra of type $E_7$ in terms of first-order differential operators
on the polynomial algebra in 27 independent variables. Section 4 is
devoted to the explicit presentation of the functor from $E_6$-{\bf
Mod} to $E_7$-{\bf Mod}. Finally in Section 6, we determine a
condition  for the functor to map a finite-dimensional irreducible
$E_6$-module to an infinite-dimensional irreducible $E_7$-module.

 \section{Polynomial Representation of $E_6$ via $E_7$}

\quad \quad In this section, we will explicitly construct the
27-dimensional basic irreducible representation of $E_6$.

 For convenience, we will use the notion
$$\ol{i,i+j}=\{i,i+1,i+2,...,i+j\}\eqno(2.1)$$
for integer $i$ and positive integer $j$ throughout this paper. We
start with the root lattice construction of the simple Lie algebra
of type $E_7$. As we all known, the Dynkin diagram of $E_7$ is as
follows:

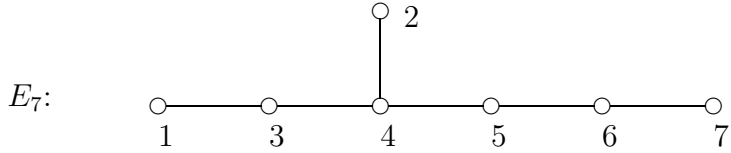
\begin{picture}(93,20)
\put(2,0){$E_7$:}\put(21,0){\circle{2}}\put(21,
-5){1}\put(22,0){\line(1,0){12}}\put(35,0){\circle{2}}\put(35,
-5){3}\put(36,0){\line(1,0){12}}\put(49,0){\circle{2}}\put(49,
-5){4}\put(49,1){\line(0,1){10}}\put(49,12){\circle{2}}\put(52,10){2}\put(50,0){\line(1,0){12}}
\put(63,0){\circle{2}}\put(63,-5){5}\put(64,0){\line(1,0){12}}\put(77,0){\circle{2}}\put(77,
-5){6}\put(78,0){\line(1,0){12}}\put(91,0){\circle{2}}\put(91,
-5){7}
\end{picture}
\vspace{0.7cm}

 \noindent Let $\{\al_i\mid i\in\ol{1,7}\}$ be the
simple positive roots corresponding to the vertices in the diagram,
and let $\Phi_{E_7}$ be the root system of $E_7$. Set
$$Q_{E_7}=\sum_{i=1}^7\mbb{Z}\al_i,\eqno(2.2)$$ the root lattice of type
$E_7$. Denote by $(\cdot,\cdot)$ the symmetric $\mbb{Z}$-bilinear
form on $Q_{E_7}$ such that
$$\Phi_{E_7}=\{\al\in Q_{E_7}\mid (\al,\al)=2\}.\eqno(2.3)$$
Define $F(\cdot,\cdot):\; Q_{E_7}\times  Q_{E_7}\rta \{\pm 1\}$ by
$$F(\sum_{i=1}^7k_i\al_i,\sum_{j=1}^7l_j\al_j)=(-1)^{\sum_{i=1}^7k_il_i+\sum_{1\leq i<j\leq 7}k_il_j
(\al_i,\al_j)},\qquad k_i,l_j\in\mbb{Z}.\eqno(2.4)$$ Then for
$\al,\be,\gm\in  Q_{E_7}$,
$$F(\al+\be,\gm)=F(\al,\gm)F(\be,\gm),\;\;F(\al,\be+\gm)=F(\al,\be)F(\al,\gm),\eqno(2.5)
$$
$$F(\al,\be)F(\be,\al)^{-1}=(-1)^{(\al,\be)},\;\;F(\al,\al)=(-1)^{(\al,\al)/2}.
\eqno(2.6)$$ In particular,
$$F(\al,\be)=-F(\be,\al)\qquad
\mbox{if}\;\;\al,\be,\al+\be\in \Phi_{E_7}.\eqno(2.7)$$

Denote
$$H_{E_7}=\sum_{i=1}^7\mbb{R}\al_i.\eqno(2.8)$$
The simple Lie algebra of type $E_7$ is
$${\cal
G}^{E_7}=H_{E_7}\oplus\bigoplus_{\al\in
\Phi_{E_7}}\mbb{R}E_{\al}\eqno(2.9)$$ with the Lie bracket
$[\cdot,\cdot]$ determined by:
 $$[H_{E_7},H_{E_7}]=0,\;\;[h,E_{\al}]=(h,\al)E_{\al},\;\;[E_{\al},E_{-\al}]=-\al,
 \eqno(2.10)$$
 $$[E_{\al},E_{\be}]=\left\{\begin{array}{ll}0&\mbox{if}\;\al+\be\not\in \Phi_{E_7},\\
 F(\al,\be)E_{\al+\be}&\mbox{if}\;\al+\be\in\Phi_{E_7}\end{array}\right.\eqno(2.11)$$
for $\al,\be\in\Phi_{E_7}$ and $h\in H_{E_7}$ (e.g., cf. [K, X1]).
Moreover, we define a bilinear form $(\cdot|\cdot)$ on ${\cal
G}^{E_6}$ by
$$(h_1|h_2)=(h_1,h_2),\;\; (h|E_{\al})=0,\;\;
 (E_{\al}|E_{\be})=-\dlt_{\al+\be,0}\eqno(2.12)$$
for $h_1,h_2\in H$ and $\al,\be\in \Phi_{E_7}$. It can be verified
that $(\cdot|\cdot)$ is a ${\cal G}^{E_7}$-invariant form, that is,
$$([u,v]|w)=-(v|[u,w])\qquad\for\;\;u,v\in{\cal
G}^{E_7}.\eqno(2.13)$$

Note that
$$Q_{E_6}=\sum_{i=1}^6\mbb{Z}\al_i\subset Q_{E_7}\eqno(2.14)$$
is the root lattice of $E_6$ and
$$\Phi_{E_6}=Q_{E_6}\bigcap \Phi_{E_7}\eqno(2.15)$$
is the root system of $E_6$. Set
$$H_{E_6}=\sum_{i=1}^6\mbb{R}\al_i.\eqno(2.16)$$
Then the subalgebra
$${\cal G}^{E_6}=H_{E_6}\oplus\bigoplus_{\al\in
\Phi_{E_6}}\mbb{R}E_{\al}\eqno(2.17)$$ of ${\cal G}^{E_7}$ is
exactly the simple Lie algebra of type $E_6$. Denote by
$\Phi_{E_6}^+$ the set of positive roots of $E_6$ and by
$\Phi_{E_7}^+$ the set of positive roots of $E_7$. The elements of
$\Phi_{E_6}^+$ are:
$$\al_1+2\al_2+2\al_3+3\al_4+2\al_5+\al_6,\eqno(2.18)$$
$$\{\al_1+\sum_{r=3}^j\al_r\mid j\in\ol{2,6}\}\bigcup
\{\sum_{r=i+1}^j\al_r\mid 2\leq i<j\leq 6\},\eqno(2.19)$$
$$\{\sum_{s=2}^j\al_s+\sum_{t=4}^k\al_t\mid 2\leq j<k\leq
6\}\eqno(2.20)$$ and
$$\{\sum_{\iota=1}^i\al_\iota+
\sum_{s=3}^j\al_s+\sum_{t=4}^k\al_t\mid 2\leq i< j<k\leq
6\}.\eqno(2.21)$$

Denote by $\bar\Phi_{E_7}^+$ the set of the following positive
roots:
$$\al_1+\sum_{r=3}^7\al_r,\qquad\al_3+2\al_4+\al_5+\sum_{i=1}^6\al_i+\sum_{r=1}^7\al_r,\eqno(2.22)$$
$$\{2\sum_{s=1}^6\al_s-\al_1+\al_4-\al_6+\sum_{r=i+1}^7\al_r\mid
i\in\ol{1,6}\} ,\qquad \{\sum_{r=i+1}^7\al_r\mid
i\in\ol{2,6}\},\eqno(2.23)$$
$$\{\sum_{s=2}^j\al_s+\sum_{t=4}^7\al_t\mid j\in\ol{2,6}\},\qquad \{\sum_{\iota=1}^i\al_\iota+
\sum_{s=3}^j\al_s+\sum_{t=4}^7\al_t\mid 2\leq i< j\leq
6\}.\eqno(2.24)$$ Then
$$\Phi_{E_7}^+=\Phi_{E_6}^+\bigcup \bar\Phi_{E_7}^+.\eqno(2.25)$$

For convenience, we also denote
$$E_{(k_1,...,k_r)}=E_\al,\;\;E'_{(k_1,...,k_r)}=E_{-\al}\qquad\for\;\;\al=\sum_{s=1}^rk_s\al_s\in\Phi_{E_7}^+,\;k_r\neq
0.\eqno(2.26)$$ Write
$$\xi_1=E_{(0,0,0,0,0,0,1)},\;\;\xi_2=E_{(0,0,0,0,0,1,1)},\;\;\xi_3=E_{(0,0,0,0,1,1,1)},\;\;\xi_4=E_{(0,0,0,1,1,1,1)},\eqno(2.27)$$
$$\xi_5=E_{(0,0,1,1,1,1,1)},\;\;\xi_6=E_{(0,1,0,1,1,1,1)},\;\;\xi_7=E_{(0,1,1,1,1,1,1)},\;\;\xi_8=E_{(1,0,1,1,1,1,1)},\eqno(2.28)$$
$$\xi_9=E_{(0,1,1,2,1,1,1)},\;\;\xi_{10}=E_{(1,1,1,1,1,1,1)},\;\;\xi_{11}=E_{(0,1,1,2,2,1,1)},\;\;\xi_{12}=E_{(1,1,1,2,1,1,1)},\eqno(2.29)$$
$$\xi_{13}=E_{(1,1,1,2,2,1,1)},\;\;\xi_{14}=E_{(0,1,1,2,2,2,1)},\;\;\xi_{15}=E_{(1,1,2,2,1,1,1)},\;\;\xi_{16}=E_{(1,1,2,2,2,1,1)},\eqno(2.30)$$
$$\xi_{17}=E_{(1,1,1,2,2,2,1)},\;\;\xi_{18}=E_{(1,1,2,3,2,1,1)},\;\;\xi_{19}=E_{(1,1,2,2,2,2,1)},\;\;\xi_{20}=E_{(1,2,2,3,2,1,1)},\eqno(2.31)$$
$$\xi_{21}=E_{(1,1,2,3,2,2,1)},\;\;\xi_{22}=E_{(1,1,2,3,3,2,1)},\;\;\xi_{23}=E_{(1,2,2,3,2,2,1)},\;\;\xi_{24}=E_{(1,2,2,3,3,2,1)},\eqno(2.32)$$
$$\xi_{25}=E_{(1,2,2,4,3,2,1)},\;\;\xi_{26}=E_{(1,2,3,4,3,2,1)},\;\;\xi_{27}=E_{(2,2,3,4,3,2,1)},\eqno(2.33)$$
$$\eta_1=E'_{(0,0,0,0,0,0,1)},\;\;\eta_2=E'_{(0,0,0,0,0,1,1)},\;\;\eta_3=E'_{(0,0,0,0,1,1,1)},\;\;\eta_4=E'_{(0,0,0,1,1,1,1)},\eqno(2.34)$$
$$\eta_5=E'_{(0,0,1,1,1,1,1)},\;\;\eta_6=E'_{(0,1,0,1,1,1,1)},\;\;\eta_7=E'_{(0,1,1,1,1,1,1)},\;\;\eta_8=E'_{(1,0,1,1,1,1,1)},\eqno(2.35)$$
$$\eta_9=E'_{(0,1,1,2,1,1,1)},\;\;\eta_{10}=E'_{(1,1,1,1,1,1,1)},\;\;\eta_{11}=E'_{(0,1,1,2,2,1,1)},\;\;\eta_{12}=E'_{(1,1,1,2,1,1,1)},\eqno(2.36)$$
$$\eta_{13}=E'_{(1,1,1,2,2,1,1)},\;\;\eta_{14}=E'_{(0,1,1,2,2,2,1)},\;\;\eta_{15}=E'_{(1,1,2,2,1,1,1)},\;\;\eta_{16}=E'_{(1,1,2,2,2,1,1)},\eqno(2.37)$$
$$\eta_{17}=E'_{(1,1,1,2,2,2,1)},\;\;\eta_{18}=E'_{(1,1,2,3,2,1,1)},\;\;\eta_{19}=E'_{(1,1,2,2,2,2,1)},\;\;\eta_{20}=E'_{(1,2,2,3,2,1,1)},\eqno(2.38)$$
$$\eta_{21}=E'_{(1,1,2,3,2,2,1)},\;\;\eta_{22}=E'_{(1,1,2,3,3,2,1)},\;\;\eta_{23}=E'_{(1,2,2,3,2,2,1)},\;\;\eta_{24}=E'_{(1,2,2,3,3,2,1)},\eqno(2.39)$$
$$\eta_{25}=E'_{(1,2,2,4,3,2,1)},\;\;\eta_{26}=E'_{(1,2,3,4,3,2,1)},\;\;\eta_{27}=E'_{(2,2,3,4,3,2,1)}.\eqno(2.40)$$
Set
$${\cal G}_-=\sum_{i=1}^{27}\mbb{C}\eta_i,\qquad {\cal G}_0={\cal
G}^{E_6}+\mbb{C}\al_7,\qquad {\cal
G}_+=\sum_{i=1}^{27}\mbb{C}\eta_i.\eqno(2.41)$$ Then ${\cal G}_\pm$
are abelian subalgebras of ${\cal G}^{E_7}$ and ${\cal G}_0$ is a
maximal reductive Lie subalgebra of ${\cal G}^{E_7}$. Moreover,
$$[{\cal G}_+,{\cal G}_-]\subset {\cal G}_0,\qquad[{\cal G}_0,{\cal
G}_\pm]\subset{\cal G}_\pm,\qquad{\cal G}^{E_7}={\cal
G}_-\oplus{\cal G}_0\oplus {\cal G}_+.\eqno(2.42)$$ Denote by
$\lmd_i$ the $i$th fundamental weight of ${\cal G}^{E_6}$. With
respect to the adjoint representation of ${\cal G}^{E_7}$, ${\cal
G}_+$ forms an irreducible ${\cal G}^{E_6}$-module with highest
weight $\lmd_1$ and ${\cal G}_-$ forms an irreducible ${\cal
G}^{E_6}$-module with highest weight $\lmd_6$.

Set
$${\cal A}=\mbb{C}[x_1,x_2,...,x_{27}],\eqno(2.43)$$
the polynomial algebra in $x_1,x_2,...,x_{27}$. Write
$$[u,\eta_i]=\sum_{j=1}^{27}\vf_{i,j}(u)\eta_j\qquad\for\;\;i\in\ol{1,27},\;u\in{\cal
G}_0,\eqno(2.44)$$ where $\vf_{i,j}(u)\in\mbb{C}$. Define an action
of ${\cal G}_0$ on ${\cal A}$ by
$$u(f)=\sum_{i,j=1}^{27}\vf_{i,j}(u)x_j\ptl_{x_i}(f)\qquad\for\;\;u\in{\cal
G}_0,\;f\in{\cal A}.\eqno(2.45)$$ Then ${\cal A}$ forms a ${\cal
G}_0$-module and the subspace
$$V=\sum_{i=1}^{27}\mbb{C}x_i\eqno(2.46)$$
forms a ${\cal G}_0$-submodule isomorphic to ${\cal G}_-$, where the
isomorphism is determined by $x_i\mapsto \eta_i$  for
$i\in\ol{1,27}$.

Denote by $\mbb{N}$ the set of nonnegative integers. Write
$$x^\al=\prod_{i=1}^{27}x_i^{\al_i},\;\;\ptl^\al=\prod_{i=1}^{27}\ptl_{x_i}^{\al_i}\qquad\for\;\;
\al=(\al_1,\al_2,...,\al_{27})\in\mbb{N}^{27}.\eqno(2.47)$$ Let
$$\mbb{A}=\sum_{\al\in\mbb{N}^{27}}{\cal A}\ptl^\al\eqno(2.48)$$
be the algebra of differential operators on ${\cal A}$. Then the
linear transformation $\tau$ determined by
$$\tau(x^\be\ptl^\gm)=x^\gm\ptl^\be\qquad\for\;\;\be,\gm\in\mbb{N}^{27}\eqno(2.49)$$
is an involutive anti-automorphism of $\mbb{A}$.

Thanks to (2.10), (2.11), (2.44) and (2.45), we find the following
representation formulas of ${\cal G}_0$:
$$E_{\al_1}|_{\cal A}=x_5\ptl_{x_8}+x_7\ptl_{x_{10}}+x_9\ptl_{x_{12}}+x_{11}
\ptl_{x_{13}}+x_{14}\ptl_{x_{17}}+x_{26}\ptl_{x_{27}},\eqno(2.50)$$
$$E_{\al_2}|_{\cal A}=x_4\ptl_{x_6}+x_5\ptl_{x_7}+x_8\ptl_{x_{10}}
-x_{18}\ptl_{x_{20}}-x_{21}\ptl_{x_{23}}-x_{22}\ptl_{x_{24}},\eqno(2.51)$$
$$E_{\al_3}|_{\cal A}=x_4\ptl_{x_5}+x_6\ptl_{x_7}+x_{12}\ptl_{x_{15}}+x_{13}\ptl_{x_{16}}
+x_{17}\ptl_{x_{19}}-x_{25}\ptl_{x_{26}},\eqno(2.52)$$
$$E_{\al_4}|_{\cal A}=x_3\ptl_{x_4}-x_7\ptl_{x_9}-x_{10}\ptl_{x_{12}}
-x_{16}\ptl_{x_{18}}-x_{19}\ptl_{x_{21}}-x_{24}\ptl_{x_{25}},\eqno(2.53)$$
$$E_{\al_5}|_{\cal A}=x_2\ptl_{x_3}-x_9\ptl_{x_{11}}
-x_{12}\ptl_{x_{13}}-x_{15}\ptl_{x_{16}}-x_{21}\ptl_{x_{22}}-x_{23}\ptl_{x_{24}},\eqno(2.54)$$
$$E_{\al_6}|_{\cal A}=x_1\ptl_{x_2}-x_{11}\ptl_{x_{14}}-x_{13}\ptl_{x_{17}}
-x_{16}\ptl_{x_{19}}-x_{18}\ptl_{x_{21}}-x_{20}\ptl_{x_{23}},\eqno(2.55)$$
$$E_{(1,0,1)}|_{\cal A}=x_4\ptl_{x_8}+x_6\ptl_{x_{10}}-x_9\ptl_{x_{15}}+x_{11}\ptl_{x_{16}}-x_{14}\ptl_{x_{19}}
+x_{25}\ptl_{x_{27}},\eqno(2.56)$$
$$E_{(0,1,0,1)}|_{\cal A}=x_3\ptl_{x_6}+x_5\ptl_{x_9}+x_8\ptl_{x_{12}}+x_{16}\ptl_{x_{20}}
+x_{19}\ptl_{x_{23}}-x_{22}\ptl_{x_{25}},\eqno(2.57)$$
$$E_{(0,0,1,1)}|_{\cal A}=x_3\ptl_{x_5}+x_6\ptl_{x_9}-x_{10}\ptl_{x_{15}}+x_{13}\ptl_{x_{18}}
+x_{17}\ptl_{x_{21}}+x_{24}\ptl_{x_{26}},\eqno(2.58)$$
$$E_{(0,0,0,1,1)}|_{\cal A}=x_2\ptl_{x_4}-x_7\ptl_{x_{11}}-x_{10}\ptl_{x_{13}}+x_{15}\ptl_{x_{18}}
-x_{19}\ptl_{x_{22}}+x_{23}\ptl_{x_{25}},\eqno(2.59)$$
$$E_{(0,0,0,0,1,1)}|_{\cal A}=x_1\ptl_{x_3}-x_9\ptl_{x_{14}}-x_{12}\ptl_{x_{17}}
-x_{15}\ptl_{x_{19}}+x_{18}\ptl_{x_{22}}+x_{20}\ptl_{x_{24}},\eqno(2.60)$$
$$E_{(1,0,1,1)}|_{\cal A}=x_3\ptl_{x_8}+x_6\ptl_{x_{12}}+x_7\ptl_{x_{15}}-x_{11}\ptl_{x_{18}}
-x_{14}\ptl_{x_{21}} +x_{24}\ptl_{x_{27}},\eqno(2.61)$$
$$E_{(0,1,1,1)}|_{\cal A}=x_3\ptl_{x_7}-
x_4\ptl_{x_9}+x_8\ptl_{x_{15}}-x_{13}\ptl_{x_{20}}-x_{17}\ptl_{x_{23}}+x_{22}\ptl_{x_{26}},\eqno(2.62)$$
$$E_{(0,1,0,1,1)}|_{\cal A}=x_2\ptl_{x_6}+x_5\ptl_{x_{11}}+x_8\ptl_{x_{13}}-x_{15}\ptl_{x_{20}}
+x_{19}\ptl_{x_{x_{24}}}+x_{21}\ptl_{x_{25}},\eqno(2.63)$$
$$E_{(0,0,1,1,1)}|_{\cal A}=x_2\ptl_{x_5}+x_6\ptl_{x_{11}}-x_{10}\ptl_{x_{16}}-x_{12}\ptl_{x_{18}}
+x_{17}\ptl_{x_{22}}-x_{23}\ptl_{x_{26}},\eqno(2.64)$$
$$E_{(0,0,0,1,1,1)}|_{\cal A}=x_1\ptl_{x_4}-x_7\ptl_{x_{14}}-x_{10}\ptl_{x_{17}}+x_{15}\ptl_{x_{21}}+x_{16}\ptl_{x_{22}}-x_{20}\ptl_{x_{25}},
\eqno(2.65)$$
$$E_{(1,1,1,1)}|_{\cal A}=x_3\ptl_{x_{10}}-x_4\ptl_{x_{12}}-x_5\ptl_{x_{15}}+x_{11}\ptl_{x_{20}}+x_{14}\ptl_{x_{23}}
+x_{22}\ptl_{x_{27}},\eqno(2.66)$$
$$E_{(1,0,1,1,1)}|_{\cal A}=x_2\ptl_{x_8}+x_6\ptl_{x_{13}}+x_7\ptl_{x_{16}}+x_9\ptl_{x_{18}}-x_{14}\ptl_{x_{22}}-x_{23}\ptl_{x_{27}},
\eqno(2.67)$$
$$E_{(0,1,1,1,1)}|_{\cal A}=x_2\ptl_{x_7}-x_4\ptl_{x_{11}}+x_8\ptl_{x_{16}}+x_{12}\ptl_{x_{20}}
-x_{17}\ptl_{x_{24}}-x_{21}\ptl_{x_{26}},\eqno(2.68)$$
$$E_{(0,1,0,1,1,1)}|_{\cal A}=x_1\ptl_{x_6}+x_5\ptl_{x_{14}}
+x_8\ptl_{x_{17}}-x_{15}\ptl_{x_{23}}-x_{16}\ptl_{x_{24}}-x_{18}\ptl_{x_{25}},\eqno(2.69)$$
$$E_{(0,0,1,1,1,1)}|_{\cal A}=x_1\ptl_{x_5}+x_6\ptl_{x_{14}}-x_{10}\ptl_{x_{19}}-x_{12}\ptl_{x_{21}}-x_{13}\ptl_{x_{22}}
+x_{20}\ptl_{x_{26}},\eqno(2.70)$$
$$E_{(1,1,1,1,1}|_{\cal A}=x_2\ptl_{x_{10}}-x_4\ptl_{x_{13}}-x_5\ptl_{x_{16}}-x_9\ptl_{x_{20}}+x_{14}\ptl_{x_{24}}-x_{21}\ptl_{x_{27}},
\eqno(2.71)$$
$$E_{(1,0,1,1,1,1)}|_{\cal A}=x_1\ptl_{x_8}+x_6\ptl_{x_{17}}+x_7\ptl_{x_{19}}+x_9\ptl_{x_{21}}+x_{11}\ptl_{x_{22}}
+x_{20}\ptl_{x_{27}},\eqno(2.72)$$
$$E_{(0,1,1,2,1)}|_{\cal A}=x_2\ptl_{x_9}-x_3\ptl_{x_{11}}+x_8\ptl_{x_{18}}-x_{10}\ptl_{x_{20}}
-x_{17}\ptl_{x_{25}}+x_{19}\ptl_{x_{26}},\eqno(2.73)$$
$$E_{(0,1,1,1,1,1)}|_{\cal
A}=x_1\ptl_{x_7}-x_4\ptl_{x_{14}}+x_8\ptl_{x_{19}}+x_{12}\ptl_{x_{23}}+x_{13}\ptl_{x_{24}}
+x_{18}\ptl_{x_{26}},\eqno(2.74)$$
$$E_{(1,1,1,2,1)}|_{\cal A}=x_2\ptl_{x_{12}}-x_3\ptl_{x_{13}}
-x_5\ptl_{x_{18}}+x_7\ptl_{x_{20}}+x_{14}\ptl_{x_{25}}+x_{19}\ptl_{x_{27}},
\eqno(2.75)$$
$$E_{(1,1,1,1,1,1)}|_{\cal A}=x_1\ptl_{x_{10}}
-x_4\ptl_{x_{17}}-x_5\ptl_{x_{19}}-x_9\ptl_{x_{23}}-x_{11}\ptl_{x_{24}}+x_{18}\ptl_{x_{27}},
\eqno(2.76)$$
$$E_{(0,1,1,2,1,1)}|_{\cal A}=x_1\ptl_{x_9}-x_3\ptl_{x_{14}}+x_8\ptl_{x_{21}}-x_{10}\ptl_{x_{23}}
+x_{13}\ptl_{x_{25}}-x_{16}\ptl_{x_{26}},\eqno(2.77)$$
$$E_{(1,1,2,2,1)}|_{\cal A}=x_2\ptl_{x_{15}}-x_3\ptl_{x_{16}}
+x_4\ptl_{x_{18}}-x_6\ptl_{x_{20}}-x_{14}\ptl_{x_{26}}-x_{17}\ptl_{x_{27}},
\eqno(2.78)$$
$$E_{(1,1,1,2,1,1)}|_{\cal A}=x_1\ptl_{x_{12}}
-x_3\ptl_{x_{17}}-x_5\ptl_{x_{21}}+x_7\ptl_{x_{23}}-x_{11}\ptl_{x_{25}}-x_{16}\ptl_{x_{27}},
\eqno(2.79)$$
$$E_{(0,1,1,2,2,1)}|_{\cal A}=x_1\ptl_{x_{11}}
-x_2\ptl_{x_{14}}+x_8\ptl_{x_{22}}-x_{10}\ptl_{x_{24}}
-x_{12}\ptl_{x_{25}}+x_{15}\ptl_{x_{26}},\eqno(2.80)$$
$$E_{(1,1,2,2,1,1)}|_{\cal A}=x_1\ptl_{x_{15}}
-x_3\ptl_{x_{19}}+x_4\ptl_{x_{21}}-x_6\ptl_{x_{23}}+x_{11}\ptl_{x_{26}}+x_{13}\ptl_{x_{27}},
\eqno(2.81)$$
$$E_{(1,1,1,2,2,1)}|_{\cal A}=x_1\ptl_{x_{13}}
-x_2\ptl_{x_{17}}-x_5\ptl_{x_{22}}+x_7\ptl_{x_{24}}
+x_9\ptl_{x_{25}}+x_{15}\ptl_{x_{27}},\eqno(2.82)$$
$$E_{(1,1,2,2,2,1)}|_{\cal A}=x_1\ptl_{x_{16}}
-x_2\ptl_{x_{19}}+x_4\ptl_{x_{22}}-x_6\ptl_{x_{24}}
-x_9\ptl_{x_{26}}-x_{12}\ptl_{x_{27}},\eqno(2.83)$$
$$E_{(1,1,2,3,2,1)}|_{\cal A}=x_1\ptl_{x_{18}}
-x_2\ptl_{x_{21}}+x_3\ptl_{x_{22}}-x_6\ptl_{x_{25}}+x_7\ptl_{x_{26}}+x_{10}\ptl_{x_{27}},
\eqno(2.84)$$
$$E_{(1,2,2,3,2,1)}|_{\cal A}=x_1\ptl_{x_{20}}
-x_2\ptl_{x_{23}}+x_3\ptl_{x_{24}}-x_4\ptl_{x_{25}}+x_5\ptl_{x_{26}}+x_8\ptl_{x_{27}},
\eqno(2.85)$$
$$E_{-\al}|_{\cal A}=-\tau(E_{\al}|_{\cal
A})\qquad\for\;\;\al\in\Phi_{E_6}^+,\eqno(2.86)$$
$$\al_r|_{\cal A}=\sum_{i=1}^{27}a_{i,r}x_i\ptl_{x_i}\qquad\for\;\;r\in\ol{1,7}\eqno(2.87)$$
with $a_{i,r}$ given in the following table:
 \begin{center}{\bf \large Table 1}\end{center}
{\small\begin{center}\begin{tabular}{|r||r|r|r|r|r|r|r||r||r|r|r|r|r|r|r|}\hline
$i$&$a_{i,1}$&$a_{i,2}$&$a_{i,3}$&$a_{i,4}$&$a_{i,5}$&$a_{i,6}$&$a_{i,7}$&
$i$&$a_{i,1}$&$a_{i,2}$&$a_{i,3}$&$a_{i,4}$&$a_{i,5}$&
$a_{i,6}$&$a_{i,7}$
\\\hline\hline 1&0&0&0&0&0&1&$-2$&2&0&0&0&0&1&$-1$&$-1$\\\hline
3&0&0&0&1&$-1$&0&$-1$&4&$0$&1&1&$-1$&0&0&$-1$
\\\hline 5&1&1&$-1$&0&0&0&$-1$&6&0&$-1$&1&0&0&0&$-1$
\\\hline 7&1&$-1$&$-1$&1&0&0&$-1$&8&$-1$&1&0&$0$&0&0&$-1$
\\\hline 9&1&0&0&$-1$&1&0&$-1$&10&$-1$&$-1$&0&1&0&0&$-1$
\\\hline 11&1&0&0&0&$-1$&1&$-1$&12&$-1$&0&1&$-1$&1&0&$-1$
\\\hline 13&$-1$&0&1&0&$-1$&1&$-1$&14&1&0&0&0&0&$-1$&0
\\\hline 15&0&0&$-1$&0&1&0&$-1$&16&0&0&$-1$&1&$-1$&1&$-1$
\\\hline 17&$-1$&0&1&0&0&$-1$&0&18&0&1&0&$-1$&0&1&$-1$
\\\hline 19&0&0&$-1$&1&0&$-1$&0&20&0&$-1$&0&0&0&1&$-1$
\\\hline 21&0&1&0&$-1$&1&$-1$&0&22&0&1&0&0&$-1$&0&0
\\\hline 23&0&$-1$&0&0&1&$-1$&0&24&0&$-1$&0&1&$-1$&0&0
\\\hline 25&0&0&1&$-1$&0&0&0&26&1&0&$-1$&0&0&0&0
\\\hline 27&$-1$&0&0&0&0&0&0
&&&&&&&&\\\hline\end{tabular}\end{center}}

Let
$$\widehat\al=2\al_1+3\al_2+4\al_3+6\al_4+5\al_5+4\al_6+3\al_7.\eqno(2.88)$$
Then
$$(\widehat\al,\al_r)=0\qquad\for\;\;r\in\ol{1,6}\eqno(2.89)$$
by the Dynkin diagram of $E_7$. Thanks to (2.10),
$$[\widehat\al,{\cal G}^{D_6}]=0.\eqno(2.90)$$
By Schur's Lemma, $\widehat\al|_V=c\sum_{i=1}^{27}x_i\ptl_{x_i}$.
According to the coefficients of $x_1\ptl_{x_1}$ in (2.87) and Table
1, we have that
$$\widehat\al|_{\cal A}=\sum_{i=1}^{16}x_i\ptl_{x_i}=D\eqno(2.91)$$
is the degree operator on ${\cal A}$.

 \section{Extended Representation of $E_7$}

\quad \quad In this section, we realize the simple Lie algebra of
type $E_7$ in terms of first-order differential operators on the
polynomial algebra in 27 independent variables.

Recall that  a singular vector of ${\cal G}^{E_6}$ is a nonzero
weight vector annihilated by positive root vectors. According to
Table 1, we find a ${\cal G}^{E_6}$-singular vector
$$\zeta_1=x_1x_{14}+x_2x_{11}+x_3x_9+x_4x_7-x_5x_6
\eqno(3.1)$$ of weight $\lmd_1$ as (3.1)-(3.6) in [X2]. Note
$$E_{-\al_1}|_{\cal A}=-x_8\ptl_{x_5}-x_{10}\ptl_{x_7}-x_{12}\ptl_{x_9}-x_{13}
\ptl_{x_{11}}-x_{17}\ptl_{x_{14}}-x_{27}\ptl_{x_{26}},\eqno(3.2)$$
$$E_{-\al_2}|_{\cal A}=-x_6\ptl_{x_4}-x_7\ptl_{x_5}-x_{10}\ptl_{x_8}
+x_{20}\ptl_{x_{18}}+x_{23}\ptl_{x_{21}}+x_{24}\ptl_{x_{22}},\eqno(3.3)$$
$$E_{-\al_3}|_{\cal A}=-x_5\ptl_{x_4}-x_7\ptl_{x_6}-x_{15}\ptl_{x_{12}}-x_{16}\ptl_{x_{13}}
-x_{19}\ptl_{x_{17}}+x_{26}\ptl_{x_{25}},\eqno(3.4)$$
$$E_{-\al_4}|_{\cal A}=-x_4\ptl_{x_3}+x_9\ptl_{x_7}+x_{12}\ptl_{x_{10}}
+x_{18}\ptl_{x_{16}}+x_{21}\ptl_{x_{19}}+x_{25}\ptl_{x_{24}},\eqno(3.5)$$
$$E_{-\al_5}|_{\cal A}=-x_3\ptl_{x_2}+x_{11}\ptl_{x_9}
+x_{13}\ptl_{x_{12}}+x_{16}\ptl_{x_{15}}+x_{22}\ptl_{x_{21}}+x_{24}\ptl_{x_{23}},\eqno(3.6)$$
$$E_{-\al_6}|_{\cal A}=-x_2\ptl_{x_1}+x_{14}\ptl_{x_{11}}+x_{17}\ptl_{x_{13}}
+x_{19}\ptl_{x_{16}}+x_{21}\ptl_{x_{18}}+x_{23}\ptl_{x_{20}}\eqno(3.7)$$
by (2.50)-(2.55) and (2.86).

As (3.7)-(3.24) in [X2], we set by (3,2)-(3.7) that
$$\zeta_2=x_1x_{17}+x_2x_{13}+x_3x_{12}+x_4x_{10}-x_6x_8,
\eqno(3.8)$$
$$\zeta_3=x_1x_{19}+x_2x_{16}+x_3x_{15}+x_5x_{10}-x_7x_8,
\eqno(3.9)$$
$$\zeta_4=-x_1x_{21}-x_2x_{18}+x_4x_{15}-x_5x_{12}
+x_8x_9, \eqno(3.10)$$
$$\zeta_5=x_1x_{22}-x_3x_{18}-x_4x_{16}+x_5x_{13}
-x_8x_{11}, \eqno(3.11)$$
$$\zeta_6=x_1x_{23}+x_2x_{20}+x_6x_{15}-x_7x_{12}
+x_9x_{10}, \eqno(3.12)$$
$$\zeta_7=-x_1x_{24}+x_3x_{20}-x_6x_{16}+x_7x_{13}
-x_{10}x_{11}, \eqno(3.13)$$
$$\zeta_8=x_2x_{22}+x_3x_{21}+x_4x_{19}-x_5x_{17}
+x_8x_{14}, \eqno(3.14)$$
$$\zeta_9=x_1x_{25}+x_4x_{20}+x_6x_{18}-x_9x_{13}
+x_{11}x_{12}, \eqno(3.15)$$
$$\zeta_{10}=-x_2x_{24}-x_3x_{23}+x_6x_{19}-x_7x_{17}
+x_{10}x_{14}, \eqno(3.16)$$
$$\zeta_{11}=x_1x_{26}-x_5x_{20}-x_7x_{18}+x_9x_{16}
-x_{11}x_{15}, \eqno(3.17)$$
$$\zeta_{12}=x_2x_{25}-x_4x_{23}-x_6x_{21}+x_9x_{17}
-x_{12}x_{14}, \eqno(3.18)$$
$$\zeta_{13}=x_3x_{25}+x_4x_{24}+x_6x_{22}-x_{11}x_{17}
+x_{13}x_{14}, \eqno(3.19)$$
$$\zeta_{14}=x_1x_{27}-x_8x_{20}-x_{10}x_{18}+x_{12}x_{16}
-x_{13}x_{15}, \eqno(3.20)$$
$$\zeta_{15}=
x_2x_{26}+x_5x_{23}+x_7x_{21}-x_9x_{19}+x_{14}x_{15}, \eqno(3.21)$$
$$\zeta_{16}=x_3x_{26}-x_5x_{24}-x_7x_{22}+x_{11}x_{19}-x_{14}x_{16},
\eqno(3.22)$$
$$\zeta_{17}=
-x_2x_{27}-x_8x_{23}-x_{10}x_{21}+x_{12}x_{19}-x_{15}x_{17},
\eqno(3.23)$$
$$\zeta_{18}=x_4x_{26}+x_5x_{25}+x_9x_{22}-x_{11}x_{21}+x_{14}x_{18},
\eqno(3.24)$$
$$\zeta_{19}=-x_3x_{27}+x_8x_{24}+x_{10}x_{22}-x_{13}x_{19}+x_{16}x_{17},
\eqno(3.25)$$
$$\zeta_{20}=-x_6x_{26}-x_7x_{25}+x_9x_{24}-x_{11}x_{23}+x_{14}x_{20},
\eqno(3.26)$$
$$\zeta_{21}=-x_4x_{27}-x_8x_{25}-x_{12}x_{22}+x_{13}x_{21}-x_{17}x_{18},
\eqno(3.27)$$
$$\zeta_{22}=x_5x_{27}-x_8x_{26}+x_{15}x_{22}-x_{16}x_{21}+x_{18}x_{19},
\eqno(3.28)$$
$$\zeta_{23}=x_6x_{27}+x_{10}x_{25}-x_{12}x_{24}+x_{13}x_{23}-x_{17}x_{20},
\eqno(3.29)$$
$$\zeta_{24}=-x_7x_{27}+x_{10}x_{26}+x_{15}x_{24}-x_{16}x_{23}+x_{19}x_{20},
\eqno(3.30)$$
$$\zeta_{25}=-x_9x_{27}+x_{12}x_{26}+x_{15}x_{25}-x_{18}x_{23}+x_{20}x_{21},
\eqno(3.31)$$
$$\zeta_{26}=-x_{11}x_{27}+x_{13}x_{26}+x_{16}x_{25}-x_{18}x_{24}+x_{20}x_{22},
\eqno(3.32)$$
$$\zeta_{27}=-x_{14}x_{27}+x_{17}x_{26}+x_{19}x_{25}-x_{21}x_{24}+x_{22}x_{23}.
\eqno(3.33)$$ Then
$$U=\sum_{i=1}^{27}\mbb{C}\zeta_i\eqno(3.34)$$
forms an irreducible ${\cal G}^{E_6}$-module with highest weight
$\lmd_1$.

By (2.50)-(2.55) and (3.1)-(3.33), we calculate
$$E_{\al_1}|_U=\zeta_1\ptl_{\zeta_2}+\zeta_{11}\ptl_{\zeta_{14}}-\zeta_{15}\ptl_{\zeta_{17}}
-\zeta_{16}\ptl_{\zeta_{19}}-\zeta_{18}\ptl_{\zeta_{21}}-\zeta_{20}\ptl_{\zeta_{23}},\eqno(3.35)$$
$$E_{\al_2}|_U=\zeta_4\ptl_{\zeta_6}+\zeta_5\ptl_{\zeta_7}+\zeta_8\ptl_{\zeta_{10}}
-\zeta_{18}\ptl_{\zeta_{20}}-\zeta_{21}\ptl_{\zeta_{23}}-\zeta_{22}\ptl_{\zeta_{24}},\eqno(3.36)$$
$$E_{\al_3}|_U=\zeta_2\ptl_{\zeta_3}-\zeta_9\ptl_{\zeta_{11}}
-\zeta_{12}\ptl_{\zeta_{15}}-\zeta_{13}\ptl_{\zeta_{16}}-\zeta_{21}\ptl_{\zeta_{22}}-\zeta_{23}\ptl_{\zeta_{24}},\eqno(3.37)$$
$$E_{\al_4}|_U=\zeta_3\ptl_{\zeta_4}+\zeta_7\ptl_{\zeta_9}+\zeta_{10}\ptl_{\zeta_{12}}
+\zeta_{16}\ptl_{\zeta_{18}}+\zeta_{19}\ptl_{\zeta_{21}}-\zeta_{24}\ptl_{\zeta_{25}},\eqno(3.38)$$
$$E_{\al_5}|_U=\zeta_4\ptl_{\zeta_5}+\zeta_6\ptl_{\zeta_7}+\zeta_{12}\ptl_{\zeta_{13}}+\zeta_{15}\ptl_{\zeta_{16}}
+\zeta_{17}\ptl_{\zeta_{19}}-\zeta_{25}\ptl_{\zeta_{26}},\eqno(3.39)$$
$$E_{\al_6}|_U=\zeta_5\ptl_{\zeta_8}+\zeta_7\ptl_{\zeta_{10}}+\zeta_9\ptl_{\zeta_{12}}+\zeta_{11}
\ptl_{\zeta_{15}}
-\zeta_{14}\ptl_{\zeta_{17}}-\zeta_{26}\ptl_{\zeta_{27}},\eqno(3.40)$$
$$E_{-\al_1}|_U=-\zeta_2\ptl_{\zeta_1}-\zeta_{14}\ptl_{\zeta_{11}}+\zeta_{17}\ptl_{\zeta_{15}}
+\zeta_{19}\ptl_{\zeta_{16}}+\zeta_{21}\ptl_{\zeta_{18}}+\zeta_{23}\ptl_{\zeta_{20}},\eqno(3.41)$$
$$E_{-\al_2}|_U=-\zeta_6\ptl_{\zeta_4}-\zeta_7\ptl_{\zeta_5}-\zeta_{10}\ptl_{\zeta_8}
+\zeta_{20}\ptl_{\zeta_{18}}+\zeta_{23}\ptl_{\zeta_{21}}+\zeta_{24}\ptl_{\zeta_{22}},\eqno(3.42)$$
$$E_{-\al_3}|_U=-\zeta_3\ptl_{\zeta_2}+\zeta_{11}\ptl_{\zeta_9}
+\zeta_{15}\ptl_{\zeta_{12}}+\zeta_{16}\ptl_{\zeta_{13}}+\zeta_{22}\ptl_{\zeta_{21}}+\zeta_{24}\ptl_{\zeta_{23}},\eqno(3.43)$$
$$E_{-\al_4}|_U=-\zeta_4\ptl_{\zeta_3}-\zeta_9\ptl_{\zeta_7}-\zeta_{12}\ptl_{\zeta_{10}}
-\zeta_{18}\ptl_{\zeta_{16}}-\zeta_{21}\ptl_{\zeta_{19}}+\zeta_{25}\ptl_{\zeta_{24}},\eqno(3.44)$$
$$E_{-\al_5}|_U=-\zeta_5\ptl_{\zeta_4}-\zeta_7\ptl_{\zeta_6}-\zeta_{13}\ptl_{\zeta_{12}}-\zeta_{16}\ptl_{\zeta_{15}}
-\zeta_{19}\ptl_{\zeta_{17}}+\zeta_{26}\ptl_{\zeta_{25}},\eqno(3.45)$$
$$E_{-\al_6}|_U=-\zeta_8\ptl_{\zeta_5}-\zeta_{10}\ptl_{\zeta_7}-\zeta_{12}\ptl_{\zeta_9}-\zeta_{15}
\ptl_{\zeta_{11}}
+\zeta_{17}\ptl_{\zeta_{14}}+\zeta_{27}\ptl_{\zeta_{26}}.\eqno(3.46)$$
Moreover, (2.84) and Table 1 give
$$\al_r|_U=\sum_{i=1}^{27}b_{i,r}\zeta_i\ptl_{\zeta_i}\qquad\for\;\;r\in\ol{1,6}\eqno(3.47)$$
with $b_{i,r}$ given in the following table:
\begin{center}{\bf \large Table 2}\end{center}
\begin{center}\begin{tabular}{|r||r|r|r|r|r|r||r||r|r|r|r|r|r|}\hline
$i$&$b_{i,1}$&$b_{i,2}$&$b_{i,3}$&$b_{i,4}$&$b_{i,5}$&$b_{i,6}$&$i$&$b_{i,1}$&$b_{i,2}$&$b_{i,3}$&$b_{i,4}$&$b_{i,5}$&
$b_{i,6}$
\\\hline\hline 1&1&0&0&0&0&0&2&$-1$&0&1&0&0&0\\\hline 3&0&0&$-1$&1&0&0&4&$0$&1&0&$-1$&1&0
\\\hline 5&0&1&$0$&0&$-1$&1&6&$0$&$-1$&0&$0$&1&0 \\\hline
7&0&$-1$&$0$&1&$-1$&1&8&$0$&1&0&$0$&0&$-1$\\\hline
9&0&$0$&$1$&$-1$&0&1&10&$0$&$-1$&0&$1$&0&$-1$\\\hline
11&1&$0$&$-1$&$0$&0&1&12&$0$&0&1&$-1$&1&$-1$\\\hline
13&0&0&1&0&$-1$&0&14&$-1$&0&0&0&0&1\\\hline
15&1&0&$-1$&0&1&$-1$&16&1&0&$-1$&1&$-1$&0\\\hline
17&$-1$&0&0&0&1&$-1$&18&1&1&$0$&$-1$&0&0\\\hline
19&$-1$&0&0&1&$-1$&0&20&1&$-1$&$0$&0&0&0\\\hline
21&$-1$&1&1&$-1$&0&0&22&0&1&$-1$&0&0&0\\\hline
23&$-1$&$-1$&1&0&0&0&24&0&$-1$&$-1$&1&0&0\\\hline
25&0&0&0&$-1$&1&0&26&0&0&0&$0$&$-1$&1\\\hline 27&0&0&0&0&0&$-1$
&&&&&&&\\\hline\end{tabular}\end{center}\psp

Recall that the Lie bracket in the algebra $\mbb{A}$ of differential
operators (cf. (2.48)) is given by the commutator
$$[d_1,d_2]=d_1d_2-d_2d_1.\eqno(3.48)$$
Write
$${\cal D}=\sum_{i=1}^{27}\mbb{C}\ptl_{x_i}.\eqno(3.49)$$
Then ${\cal D}$ forms a ${\cal G}^{E_6}$-module with respect to the
action
$$u(\ptl)=[u|_{\cal A},\ptl]\qquad\for\;\;u\in{\cal G}^{E_6},\;\ptl\in{\cal D}.\eqno(3.50)$$ On the other hand,
${\cal G}_\pm$ (cf. (2.25)-(2.39)) form ${\cal G}^{E_6}$-modules
with respect to the adjoint representation.
 According to (2.44)
and (2.45), the linear map determined by $\eta_i\mapsto x_i$ for
$i\in\ol{1,27}$ gives a ${\cal G}^{E_6}$-module isomorphism from
${\cal G}_-$ to $V$. Moreover,
$$(\xi_i|\eta_j)=-\dlt_{i,j}\qquad\for\;\;i,j\in\ol{1,27}\eqno(3.51)$$
by (2.12). The fact $\ptl_{x_i}(x_j)=\dlt_{i,j}$ for
$i,j\in\ol{1,27}$, (2.13) and (3.51) show that the linear map
determined by $\xi_i\mapsto \ptl_{x_i}$ for $i\in\ol{1,27}$ gives a
${\cal G}^{E_6}$-module isomorphism from ${\cal G}_+$ to ${\cal D}$.
Hence we define the action of ${\cal G}_+$ on ${\cal A}$ by
$$\xi_i|_{\cal
A}=\ptl_{x_i}\qquad\for\;\;i\in\ol{1,27}.\eqno(3.52)$$

Recall the Witt Lie subalgebra of $\mbb{A}$:
$${\cal W}_{27}=\sum_{i=1}^{27}{\cal A}\ptl_{x_i}.\eqno(3.53)$$
Now we want to find the differential operators
$P_1,P_2,...,P_{27}\in {\cal W}_{24}$ such that the following action
matches the structure of ${\cal G}^{E_7}$:
$$\eta_i|_{\cal A}=P_i\qquad\for\;\;i\in\ol{1,27}.\eqno(3.54)$$
Comparing the weights in Table 1 and Table 2, we use (3.1) and
(3.8)-(3.33) to assume
\begin{eqnarray*}\qquad P_1&=&x_1D+c_1\zeta_1\ptl_{x_{14}}+c_2\zeta_2\ptl_{x_{17}}+c_3\zeta_3\ptl_{x_{19}}
+c_4\zeta_4\ptl_{x_{21}}+c_5\zeta_5\ptl_{x_{22}}
\\&&+c_6\zeta_6\ptl_{x_{23}}+c_7\zeta_7\ptl_{x_{24}}+c_9\zeta_9\ptl_{x_{25}}+c_{11}\zeta_{11}\ptl_{x_{26}}
+c_{14}\zeta_{14}\ptl_{x_{27}},\hspace{2.7cm}(3.55)\end{eqnarray*}
where $c_i\in\mbb{C}$. Imposing
\begin{eqnarray*}\qquad [\ptl_{x_1},P_1]&=&[\xi_1,\eta_1]|_{\cal A}=-\al_7|_{\cal
A}=2x_1\ptl_{x_1}+\sum_{i=2}^{13}x_i\ptl_{x_i}\\&
&+x_{15}\ptl_{x_{15}} +x_{16}\ptl_{x_{16}}+x_{18}\ptl_{x_{18}}
+x_{20}\ptl_{x_{20}}\hspace{4.9cm}(3.56)\end{eqnarray*} by (2.87)
and Table 1, we get
$$c_1=c_2=c_3=-c_4=c_5=c_6=-c_7=c_9=c_{11}=c_{14}=-1\eqno(3.57)$$
by (3.1), (3.8)-(3.13), (3.15), (3.17) and (3.20). Thus
\begin{eqnarray*}\qquad\qquad
P_1&=&x_1D-\zeta_1\ptl_{x_{14}}-\zeta_2\ptl_{x_{17}}-\zeta_3\ptl_{x_{19}}
+\zeta_4\ptl_{x_{21}}-\zeta_5\ptl_{x_{22}}
\\&&-\zeta_6\ptl_{x_{23}}+\zeta_7\ptl_{x_{24}}-\zeta_9\ptl_{x_{25}}-\zeta_{11}\ptl_{x_{26}}
-\zeta_{14}\ptl_{x_{27}}.\hspace{3.9cm}(3.58)\end{eqnarray*}

According to (3.1), (3.8)-(3.13), (3.15), (3.17) and (3.20),
 we find
$$[\ptl_{x_2},P_1]=x_1\ptl_{x_2}-x_{11}\ptl_{x_{14}}-x_{13}\ptl_{x_{17}}
-x_{16}\ptl_{x_{19}}-x_{18}\ptl_{x_{21}}-x_{20}\ptl_{x_{23}},\eqno(3.59)$$
$$[\ptl_{x_3},P_1]=x_1\ptl_{x_3}-x_9\ptl_{x_{14}}-x_{12}\ptl_{x_{17}}
-x_{15}\ptl_{x_{19}}+x_{18}\ptl_{x_{22}}+x_{20}\ptl_{x_{24}},\eqno(3.60)$$
$$[\ptl_{x_4},P_1]=x_1\ptl_{x_4}-x_7\ptl_{x_{14}}-x_{10}\ptl_{x_{17}}
+x_{15}\ptl_{x_{21}}+x_{16}\ptl_{x_{22}}-x_{20}\ptl_{x_{25}},\eqno(3.61)$$
$$[\ptl_{x_5},P_1]=x_1\ptl_{x_5}+x_6\ptl_{x_{14}}-x_{10}\ptl_{x_{19}}
-x_{12}\ptl_{x_{21}}-x_{13}\ptl_{x_{22}}+x_{20}\ptl_{x_{25}},\eqno(3.62)$$
$$[\ptl_{x_6},P_1]=x_1\ptl_{x_6}+x_5\ptl_{x_{14}}+x_8\ptl_{x_{17}}
-x_{15}\ptl_{x_{23}}-x_{16}\ptl_{x_{24}}-x_{18}\ptl_{x_{25}},\eqno(3.63)$$
$$[\ptl_{x_7},P_1]=x_1\ptl_{x_7}-x_4\ptl_{x_{14}}+x_8\ptl_{x_{19}}
+x_{12}\ptl_{x_{23}}+x_{13}\ptl_{x_{24}}+x_{18}\ptl_{x_{26}},\eqno(3.64)$$
$$[\ptl_{x_8},P_1]=x_1\ptl_{x_8}+x_6\ptl_{x_{17}}+x_7\ptl_{x_{19}}
+x_9\ptl_{x_{21}}+x_{11}\ptl_{x_{22}}+x_{20}\ptl_{x_{27}},\eqno(3.65)$$
$$[\ptl_{x_9},P_1]=x_1\ptl_{x_9}-x_3\ptl_{x_{14}}+x_8\ptl_{x_{21}}
-x_{10}\ptl_{x_{23}}+x_{13}\ptl_{x_{25}}-x_{16}\ptl_{x_{26}},\eqno(3.66)$$
$$[\ptl_{x_{10}},P_1]=x_1\ptl_{x_{10}}-x_4\ptl_{x_{17}}-x_5\ptl_{x_{19}}
-x_9\ptl_{x_{23}}-x_{11}\ptl_{x_{24}}+x_{18}\ptl_{x_{27}},\eqno(3.67)$$
$$[\ptl_{x_{11}},P_1]=x_1\ptl_{x_{11}}-x_2\ptl_{x_{14}}+x_8\ptl_{x_{22}}
-x_{10}\ptl_{x_{24}}-x_{12}\ptl_{x_{25}}+x_{15}\ptl_{x_{26}},\eqno(3.68)$$
$$[\ptl_{x_{12}},P_1]=x_1\ptl_{x_{12}}-x_3\ptl_{x_{17}}-x_5\ptl_{x_{21}}
+x_7\ptl_{x_{23}}-x_{11}\ptl_{x_{25}}-x_{16}\ptl_{x_{27}},\eqno(3.69)$$
$$[\ptl_{x_{13}},P_1]=x_1\ptl_{x_{13}}-x_2\ptl_{x_{17}}-x_5\ptl_{x_{22}}
+x_7\ptl_{x_{24}}+x_9\ptl_{x_{25}}+x_{15}\ptl_{x_{27}},\eqno(3.70)$$
$$[\ptl_{x_{15}},P_1]=x_1\ptl_{x_{15}}-x_3\ptl_{x_{19}}+x_4\ptl_{x_{21}}
-x_6\ptl_{x_{23}}+x_{11}\ptl_{x_{26}}+x_{13}\ptl_{x_{27}},\eqno(3.71)$$
$$[\ptl_{x_{16}},P_1]=x_1\ptl_{x_{16}}-x_2\ptl_{x_{19}}+x_4\ptl_{x_{22}}
-x_6\ptl_{x_{24}}-x_9\ptl_{x_{26}}-x_{12}\ptl_{x_{27}},\eqno(3.72)$$
$$[\ptl_{x_{18}},P_1]=x_1\ptl_{x_{18}}-x_2\ptl_{x_{21}}+x_3\ptl_{x_{22}}
-x_6\ptl_{x_{25}}+x_7\ptl_{x_{26}}+x_{10}\ptl_{x_{27}},\eqno(3.73)$$
$$[\ptl_{x_{20}},P_1]=x_1\ptl_{x_{20}}-x_2\ptl_{x_{23}}+x_3\ptl_{x_{24}}
-x_4\ptl_{x_{25}}+x_5\ptl_{x_{26}}+x_8\ptl_{x_{27}},\eqno(3.74)$$
$$[x_r,P_1]=0\qquad\for\;\;r=14,17,19,21,22,23,24,25,26,27.\eqno(3.75)$$
On the other hand, if $[\xi_i,\eta_1]\neq 0$ for $i\geq 2$, then it
is equal to the vector obtained by deleting the last coordinate $1$
by (2.27)-(2.33). For instance, $\xi_{18}=E_{(1,1,2,3,2,1,1)}$ and
$[\xi_{18},\eta_1]=E_{(1,1,2,3,2,1)}$. By  (3.58)-(3.73) and
correspondingly (2.55), (2.60), (2.65), (2.70), (2.69), (2.74),
(2.72), (2.77), (2.76), (2.80), (2.79), (2.82), (2.81),
(2.83)-(2.85), we have:
$$[\ptl_{x_i},P_1]=[\xi_i,\eta_1]|_{\cal
A}\qquad\for\;\;i\in\ol{1,27}.\eqno(3.76)$$

Expressions (2.50)-(2.55), (3.35)-(3.40) and (3.48) imply that $P_1$
is a ${\cal G}^{E_6}$-singular vector in ${\cal W}_{27}$ with weight
$\lmd_6$. Since $-[E_{-\al_6},x_1]=x_2$, we set
\begin{eqnarray*}
\qquad P_2&=&-[E_{-\al_6},P_1]=
x_2D-\zeta_1\ptl_{x_{11}}-\zeta_2\ptl_{x_{13}}-\zeta_3\ptl_{x_{16}}
+\zeta_4\ptl_{x_{18}}\\ &
&-\zeta_6\ptl_{x_{20}}-\zeta_8\ptl_{x_{22}}
+\zeta_{10}\ptl_{x_{24}}-\zeta_{12}\ptl_{x_{25}}-\zeta_{15}\ptl_{x_{26}}
+\zeta_{17}\ptl_{x_{27}}\hspace{3cm}(3.77)\end{eqnarray*}
 by (3.7), (3.46) and (3.58). As $-[E_{-\al_5},x_2]=x_3$, we take
\begin{eqnarray*}
\qquad P_3&=&-[E_{-\al_5},P_2]=
x_3D-\zeta_1\ptl_{x_9}-\zeta_2\ptl_{x_{12}}-\zeta_3\ptl_{x_{15}}
+\zeta_5\ptl_{x_{18}}\\ &
&-\zeta_7\ptl_{x_{20}}-\zeta_8\ptl_{x_{21}}
+\zeta_{10}\ptl_{x_{23}}-\zeta_{13}\ptl_{x_{25}}-\zeta_{16}\ptl_{x_{26}}
+\zeta_{19}\ptl_{x_{27}}\hspace{3cm}(3.78)\end{eqnarray*}
 by (3.6), (3.45) and (3.77). Thanks to $-[E_{-\al_4},x_3]=x_4$, we
 have
\begin{eqnarray*}
\qquad P_4&=&-[E_{-\al_4},P_3]=
x_4D-\zeta_1\ptl_{x_7}-\zeta_2\ptl_{x_{10}}-\zeta_4\ptl_{x_{15}}
+\zeta_5\ptl_{x_{16}}\\ &
&-\zeta_8\ptl_{x_{19}}-\zeta_9\ptl_{x_{20}}
+\zeta_{12}\ptl_{x_{23}}-\zeta_{13}\ptl_{x_{24}}-\zeta_{18}\ptl_{x_{26}}
+\zeta_{21}\ptl_{x_{27}}\hspace{3cm}(3.79)\end{eqnarray*}
 by (3.5), (3.44) and (3.78). Due to $-[E_{-\al_3},x_4]=x_5$, we
 write
\begin{eqnarray*}
\qquad P_5&=&-[E_{-\al_3},P_4]=
x_5D+\zeta_1\ptl_{x_6}-\zeta_3\ptl_{x_{10}}+\zeta_4\ptl_{x_{12}}
-\zeta_5\ptl_{x_{13}}\\ &
&+\zeta_8\ptl_{x_{17}}+\zeta_{11}\ptl_{x_{20}}
-\zeta_{15}\ptl_{x_{23}}+\zeta_{16}\ptl_{x_{24}}-\zeta_{18}\ptl_{x_{25}}
-\zeta_{22}\ptl_{x_{27}}\hspace{2.8cm}(3.80)\end{eqnarray*}
 by (3.4), (3.43) and (3.79). Since
$-[E_{-\al_2},x_4]=x_6$, we
 denote
\begin{eqnarray*}
\qquad P_6&=&-[E_{-\al_2},P_4]=
x_6D+\zeta_1\ptl_{x_5}+\zeta_2\ptl_{x_8}-\zeta_6\ptl_{x_{15}}
+\zeta_7\ptl_{x_{16}}\\ &
&-\zeta_9\ptl_{x_{18}}-\zeta_{10}\ptl_{x_{19}}
+\zeta_{12}\ptl_{x_{21}}-\zeta_{13}\ptl_{x_{22}}+\zeta_{20}\ptl_{x_{26}}
-\zeta_{23}\ptl_{x_{27}}\hspace{2.8cm}(3.81)\end{eqnarray*}
 by (3.3), (3.42) and (3.79). As
$-[E_{-\al_3},x_6]=x_7$, we
 denote
\begin{eqnarray*}
\qquad P_7&=&-[E_{-\al_3},P_6]=
x_7D-\zeta_1\ptl_{x_4}+\zeta_3\ptl_{x_8}+\zeta_6\ptl_{x_{12}}
-\zeta_7\ptl_{x_{13}}\\ &
&+\zeta_{10}\ptl_{x_{17}}+\zeta_{11}\ptl_{x_{18}}
-\zeta_{15}\ptl_{x_{21}}+\zeta_{16}\ptl_{x_{22}}+\zeta_{20}\ptl_{x_{25}}
+\zeta_{24}\ptl_{x_{27}}\hspace{2.7cm}(3.82)\end{eqnarray*}
 by (3.4), (3.43) and (3.81). Thanks to $-[E_{-\al_1},x_5]=x_8$, we
 find
\begin{eqnarray*}
\qquad P_8&=&-[E_{-\al_2},P_5]=
x_8D+\zeta_2\ptl_{x_6}+\zeta_3\ptl_{x_7}-\zeta_4\ptl_{x_9}
+\zeta_5\ptl_{x_{11}}\\ &
&-\zeta_8\ptl_{x_{14}}+\zeta_{14}\ptl_{x_{20}}
+\zeta_{17}\ptl_{x_{23}}-\zeta_{19}\ptl_{x_{24}}+\zeta_{21}\ptl_{x_{25}}
+\zeta_{22}\ptl_{x_{26}}\hspace{2.8cm}(3.83)\end{eqnarray*}
 by (3.2), (3.41) and (3.80).

Since $[E_{-\al_4},x_7]=x_9$, we take
\begin{eqnarray*}
\qquad P_9&=&[E_{-\al_4},P_7]=
x_9D-\zeta_1\ptl_{x_3}-\zeta_4\ptl_{x_8}-\zeta_6\ptl_{x_{10}}
+\zeta_9\ptl_{x_{13}}\\ &
&-\zeta_{12}\ptl_{x_{17}}-\zeta_{11}\ptl_{x_{16}}
+\zeta_{15}\ptl_{x_{19}}-\zeta_{18}\ptl_{x_{22}}-\zeta_{20}\ptl_{x_{24}}
+\zeta_{25}\ptl_{x_{27}}\hspace{2.7cm}(3.84)\end{eqnarray*}
 by (3.5), (3.44) and (3.82). As $[E_{-\al_1},x_7]=-x_{10}$, we let
\begin{eqnarray*}
\qquad P_{10}&=&-[E_{-\al_1},P_7]=
x_{10}D-\zeta_2\ptl_{x_4}-\zeta_3\ptl_{x_5}-\zeta_6\ptl_{x_9}
+\zeta_7\ptl_{x_{11}}\\ &
&-\zeta_{10}\ptl_{x_{14}}+\zeta_{14}\ptl_{x_{18}}
+\zeta_{17}\ptl_{x_{21}}-\zeta_{19}\ptl_{x_{22}}-\zeta_{23}\ptl_{x_{25}}
-\zeta_{24}\ptl_{x_{26}}\hspace{2.6cm}(3.85)\end{eqnarray*}
 by (3.2), (3.41) and (3.82). Due to $[E_{-\al_5},x_9]=x_{11}$, we take
\begin{eqnarray*}
\qquad P_{11}&=&[E_{-\al_4},P_7]=
x_{11}D-\zeta_1\ptl_{x_2}+\zeta_5\ptl_{x_8}+\zeta_7\ptl_{x_{10}}
-\zeta_9\ptl_{x_{12}}\\ &
&+\zeta_{11}\ptl_{x_{15}}+\zeta_{13}\ptl_{x_{17}}
-\zeta_{16}\ptl_{x_{19}}+\zeta_{18}\ptl_{x_{21}}+\zeta_{20}\ptl_{x_{23}}
+\zeta_{26}\ptl_{x_{27}}\hspace{2.6cm}(3.86)\end{eqnarray*}
 by (3.6), (3.45) and (3.84). Thanks to
$-[E_{-\al_1},x_9]=x_{12}$, we find
\begin{eqnarray*}
\qquad
P_{12}&=&-[E_{-\al_1},P_9]=x_{12}D-\zeta_2\ptl_{x_3}+\zeta_4\ptl_{x_5}+\zeta_6\ptl_{x_7}
-\zeta_9\ptl_{x_{11}}\\ &
&+\zeta_{12}\ptl_{x_{14}}-\zeta_{14}\ptl_{x_{16}}
-\zeta_{17}\ptl_{x_{19}}+\zeta_{21}\ptl_{x_{22}}+\zeta_{23}\ptl_{x_{24}}
-\zeta_{25}\ptl_{x_{26}}\hspace{2.6cm}(3.87)\end{eqnarray*}
 by (3.2), (3.41) and (3.84). Due to $[E_{-\al_5},x_{12}]=x_{13}$, we
 write
\begin{eqnarray*}
\qquad
P_{13}&=&[E_{-\al_5},P_{12}]=x_{13}D-\zeta_2\ptl_{x_2}-\zeta_5\ptl_{x_5}-\zeta_7\ptl_{x_7}
+\zeta_9\ptl_{x_9}\\ &
&-\zeta_{13}\ptl_{x_{14}}+\zeta_{14}\ptl_{x_{15}}
+\zeta_{19}\ptl_{x_{19}}-\zeta_{21}\ptl_{x_{21}}-\zeta_{23}\ptl_{x_{23}}
-\zeta_{26}\ptl_{x_{26}}\hspace{2.5cm}(3.88)\end{eqnarray*}
 by (3.6), (3.45) and (3.87). Equation $[E_{-\al_6},x_{11}]=x_{14}$, we
 let
\begin{eqnarray*}
\qquad P_{14}&=&[E_{-\al_6},P_{11}]=
x_{14}D-\zeta_1\ptl_{x_1}-\zeta_8\ptl_{x_8}-\zeta_{10}\ptl_{x_{10}}
+\zeta_{12}\ptl_{x_{12}}\\ &
&-\zeta_{13}\ptl_{x_{13}}-\zeta_{15}\ptl_{x_{15}}
+\zeta_{16}\ptl_{x_{16}}-\zeta_{18}\ptl_{x_{18}}-\zeta_{20}\ptl_{x_{20}}
+\zeta_{27}\ptl_{x_{27}}\hspace{2.5cm}(3.89)\end{eqnarray*}
 by (3.7), (3.46) and (3.86).

As $-[E_{-\al_3},x_{12}]=x_{15}$, we get
\begin{eqnarray*}
\qquad
P_{15}&=&-[E_{-\al_3},P_{12}]=x_{15}D-\zeta_3\ptl_{x_3}-\zeta_4\ptl_{x_4}-\zeta_6\ptl_{x_6}
+\zeta_{11}\ptl_{x_{11}}\\ &
&+\zeta_{14}\ptl_{x_{13}}-\zeta_{15}\ptl_{x_{14}}
+\zeta_{17}\ptl_{x_{17}}-\zeta_{22}\ptl_{x_{22}}-\zeta_{24}\ptl_{x_{24}}
-\zeta_{25}\ptl_{x_{25}}\hspace{2.6cm}(3.90)\end{eqnarray*}
 by (3.4), (3.43) and (3.87). Thanks to $[E_{-\al_5},x_{15}]=x_{16}$, we let
\begin{eqnarray*}
\qquad
P_{16}&=&[E_{-\al_5},P_{15}]=x_{16}D-\zeta_3\ptl_{x_2}+\zeta_5\ptl_{x_4}+\zeta_7\ptl_{x_6}
-\zeta_{11}\ptl_{x_9}\\ &
&-\zeta_{14}\ptl_{x_{12}}+\zeta_{16}\ptl_{x_{14}}
-\zeta_{19}\ptl_{x_{17}}+\zeta_{22}\ptl_{x_{21}}+\zeta_{24}\ptl_{x_{23}}
-\zeta_{26}\ptl_{x_{25}}\hspace{2.6cm}(3.91)\end{eqnarray*}
 by (3.6), (3.45) and (3.90). Due to $-[E_{-\al_1},x_{14}]=x_{17}$, we get
\begin{eqnarray*}
\qquad P_{17}&=&=-[E_{-\al_1},P_{14}]=
x_{17}D-\zeta_2\ptl_{x_1}+\zeta_8\ptl_{x_5}+\zeta_{10}\ptl_{x_7}
-\zeta_{12}\ptl_{x_9}\\ &
&+\zeta_{13}\ptl_{x_{11}}+\zeta_{17}\ptl_{x_{15}}
-\zeta_{19}\ptl_{x_{16}}+\zeta_{21}\ptl_{x_{18}}+\zeta_{23}\ptl_{x_{20}}
-\zeta_{27}\ptl_{x_{26}}\hspace{2.5cm}(3.92)\end{eqnarray*}
 by (3.2), (3.41) and (3.89). Since  $[E_{-\al_4},x_{16}]=x_{18}$, we
 find
 \begin{eqnarray*}
\qquad
P_{18}&=&[E_{-\al_4},P_{16}]=x_{18}D+\zeta_4\ptl_{x_2}+\zeta_5\ptl_{x_3}-\zeta_9\ptl_{x_6}
+\zeta_{11}\ptl_{x_7}\\ &
&+\zeta_{14}\ptl_{x_{10}}-\zeta_{18}\ptl_{x_{14}}
+\zeta_{21}\ptl_{x_{17}}-\zeta_{22}\ptl_{x_{19}}+\zeta_{25}\ptl_{x_{23}}
+\zeta_{26}\ptl_{x_{24}}\hspace{2.6cm}(3.93)\end{eqnarray*}
 by (3.5), (3.44) and (3.91).

Thanks to $[E_{-\al_6},x_{16}]=x_{19}$, we
 calculate
 \begin{eqnarray*}
\qquad
P_{19}&=&[E_{-\al_6},P_{16}]=x_{19}D-\zeta_3\ptl_{x_1}-\zeta_8\ptl_{x_4}-\zeta_{10}\ptl_{x_6}
+\zeta_{15}\ptl_{x_9}\\ &
&-\zeta_{17}\ptl_{x_{12}}-\zeta_{16}\ptl_{x_{11}}
+\zeta_{19}\ptl_{x_{13}}-\zeta_{22}\ptl_{x_{18}}-\zeta_{24}\ptl_{x_{20}}
-\zeta_{27}\ptl_{x_{25}}\hspace{2.6cm}(3.94)\end{eqnarray*}
 by (3.7), (3.46) and (3.91). As $[E_{-\al_2},x_{18}]=x_{20}$, we
 have
 \begin{eqnarray*}
\qquad
P_{20}&=&[E_{-\al_2},P_{18}]=x_{20}D-\zeta_6\ptl_{x_2}-\zeta_7\ptl_{x_3}-\zeta_9\ptl_{x_4}
+\zeta_{11}\ptl_{x_5}\\ &
&+\zeta_{14}\ptl_{x_8}-\zeta_{20}\ptl_{x_{14}}
+\zeta_{23}\ptl_{x_{17}}-\zeta_{24}\ptl_{x_{19}}-\zeta_{25}\ptl_{x_{21}}
-\zeta_{26}\ptl_{x_{22}}\hspace{2.6cm}(3.95)\end{eqnarray*}
 by (3.3), (3.42) and (3.93). Since  $[E_{-\al_4},x_{19}]=x_{21}$, we
 let
 \begin{eqnarray*}
\qquad
P_{21}&=&[E_{-\al_4},P_{19}]=x_{21}D+\zeta_4\ptl_{x_1}-\zeta_8\ptl_{x_3}+\zeta_{12}\ptl_{x_6}
-\zeta_{15}\ptl_{x_7}\\ &
&+\zeta_{17}\ptl_{x_{10}}+\zeta_{18}\ptl_{x_{11}}
-\zeta_{21}\ptl_{x_{13}}+\zeta_{22}\ptl_{x_{16}}-\zeta_{25}\ptl_{x_{20}}
+\zeta_{27}\ptl_{x_{24}}\hspace{2.6cm}(3.96)\end{eqnarray*}
 by (3.5), (3.44) and (3.94). Equation $[E_{-\al_5},x_{21}]=x_{22}$, we
 write
 \begin{eqnarray*}
\qquad
P_{22}&=&[E_{-\al_5},P_{21}]=x_{22}D-\zeta_5\ptl_{x_1}-\zeta_8\ptl_{x_2}-\zeta_{13}\ptl_{x_6}
+\zeta_{16}\ptl_{x_7}\\ &
&-\zeta_{19}\ptl_{x_{10}}-\zeta_{18}\ptl_{x_9}
+\zeta_{21}\ptl_{x_{12}}-\zeta_{22}\ptl_{x_{15}}-\zeta_{26}\ptl_{x_{20}}
-\zeta_{27}\ptl_{x_{23}}\hspace{2.6cm}(3.97)\end{eqnarray*}
 by (3.6), (3.45) and (3.96).

Note that the equation $[E_{-\al_2},x_{21}]=x_{23}$ gives rise to
 \begin{eqnarray*}
\qquad
P_{23}&=&[E_{-\al_2},P_{21}]=x_{23}D-\zeta_6\ptl_{x_1}+\zeta_{10}\ptl_{x_3}+\zeta_{12}\ptl_{x_4}
-\zeta_{15}\ptl_{x_5}\\ &
&+\zeta_{17}\ptl_{x_8}+\zeta_{20}\ptl_{x_{11}}
-\zeta_{23}\ptl_{x_{13}}+\zeta_{24}\ptl_{x_{16}}+\zeta_{25}\ptl_{x_{18}}
-\zeta_{27}\ptl_{x_{22}}\hspace{2.6cm}(3.98)\end{eqnarray*}
 by (3.3), (3.42) and (3.96). Since $[E_{-\al_2},x_{22}]=x_{24}$, we
 write
 \begin{eqnarray*}
\qquad
P_{24}&=&[E_{-\al_2},P_{22}]=x_{24}D+\zeta_7\ptl_{x_1}+\zeta_{10}\ptl_{x_2}-\zeta_{13}\ptl_{x_4}
+\zeta_{16}\ptl_{x_5}\\ &
&-\zeta_{19}\ptl_{x_8}-\zeta_{20}\ptl_{x_9}
+\zeta_{23}\ptl_{x_{12}}-\zeta_{24}\ptl_{x_{15}}+\zeta_{26}\ptl_{x_{18}}
+\zeta_{27}\ptl_{x_{21}}\hspace{2.6cm}(3.99)\end{eqnarray*}
 by (3.3), (3.42) and (3.97). As $[E_{-\al_4},x_{24}]=x_{25}$, we
 write
 \begin{eqnarray*}
\qquad
P_{25}&=&[E_{-\al_4},P_{24}]=x_{25}D-\zeta_9\ptl_{x_1}-\zeta_{12}\ptl_{x_2}-\zeta_{13}\ptl_{x_3}
-\zeta_{18}\ptl_{x_5}\\ &
&+\zeta_{21}\ptl_{x_8}+\zeta_{20}\ptl_{x_7}
-\zeta_{23}\ptl_{x_{10}}-\zeta_{25}\ptl_{x_{15}}-\zeta_{26}\ptl_{x_{16}}
-\zeta_{27}\ptl_{x_{19}}\hspace{2.4cm}(3.100)\end{eqnarray*}
 by (3.5), (3.44) and (3.99). Thanks to $[E_{-\al_3},x_{25}]=x_{26}$, we
 write
 \begin{eqnarray*}
\qquad
P_{26}&=&[E_{-\al_3},P_{25}]=x_{26}D-\zeta_{11}\ptl_{x_1}-\zeta_{15}\ptl_{x_2}-\zeta_{16}\ptl_{x_3}
-\zeta_{18}\ptl_{x_4}\\ &
&+\zeta_{22}\ptl_{x_8}+\zeta_{20}\ptl_{x_6}
-\zeta_{24}\ptl_{x_{10}}-\zeta_{25}\ptl_{x_{12}}-\zeta_{26}\ptl_{x_{13}}
-\zeta_{27}\ptl_{x_{17}}\hspace{2.5cm}(3.101)\end{eqnarray*}
 by (3.4), (3.43) and (3.100). Due to $-[E_{-\al_1},x_{26}]=x_{27}$, we
 write
 \begin{eqnarray*}
\qquad
P_{27}&=&-[E_{-\al_1},P_{26}]=x_{27}D-\zeta_{14}\ptl_{x_1}+\zeta_{17}\ptl_{x_2}+\zeta_{19}\ptl_{x_3}
+\zeta_{21}\ptl_{x_4}\\ &
&-\zeta_{22}\ptl_{x_5}-\zeta_{23}\ptl_{x_6}
+\zeta_{24}\ptl_{x_7}+\zeta_{25}\ptl_{x_9}+\zeta_{26}\ptl_{x_{11}}
+\zeta_{27}\ptl_{x_{14}}\hspace{2.8cm}(3.102)\end{eqnarray*}
 by (3.2), (3.41) and (3.101).

Set
$${\cal P}=\sum_{i=1}^{27}\mbb{C}P_i,\qquad{\cal C}_0={\cal G}^{E_6}|_{\cal A}+\mbb{C}D\eqno(3.103)$$ (cf. (2.50)-(2.87) and (2.91)) and
$${\cal
C}={\cal P}+{\cal C}_0+{\cal D}\eqno(3.104)$$ (cf. (3.49)). Then we
have:\psp

{\bf Theorem 3.1}. {\it The space ${\cal C}$ forms a Lie subalgebra
of the Witt algebra ${\cal W}_{27}$ (cf. (3.53)). Moreover, the
linear map $\vt$ determined by
$$\vt(\xi_i)=\ptl_{x_i},\;\;\vt(\eta_i)=P_i,\;\;\vt(u)=u|_{\cal
A}\qquad\for\;\;i\in\ol{1,27},\;u\in{\cal G}^{E_6}\eqno(3.105)$$
 and
$$\vt(\al_7)=-2x_1\ptl_{x_1}-\sum_{i=2}^{13}x_i\ptl_{x_i}-x_{15}\ptl_{x_{15}}-x_{16}\ptl_{x_{16}}-x_{18}\ptl_{x_{18}}-x_{20}\ptl_{x_{20}}\eqno(3.106)$$
(cf. (2.87) and Table 1) gives a Lie algebra isomorphism from ${\cal
G}^{E_7}$ to ${\cal C}$.}

{\it Proof}. Since ${\cal D}\cong {\cal G}_+$ as ${\cal
G}^{E_6}$-modules, we have
$${\cal G}_0+{\cal G}_+\stl{\vt}{\cong}{\cal C}_0+{\cal D}\eqno(3.107)$$ as Lie algebras. Denote by
$U({\cal G})$ the universal enveloping algebra of a Lie algebra
${\cal G}$. Note that
$${\cal B}_-={\cal G}_0+{\cal G}_-,\qquad{\cal B}_+={\cal G}_0+{\cal G}_+\eqno(3.108)$$
are also  Lie subalgebras of ${\cal G}^{E_7}$ and
$${\cal G}^{E_7}={\cal B}_-\oplus {\cal G}_+={\cal G}_-\oplus {\cal B}_+.\eqno(3.109)$$
We define a one-dimensional ${\cal B}_-$-module $\mbb{C}u_0$ by
$$w(u_0)=0\qquad\for\;\;w\in{\cal B}_-\bigcap{\cal G}^{E_6},\;\;\widehat\al(u_0)=-27u_0\eqno(3.110)$$
(cf. (2.88)). Let
$$\Psi=U({\cal G}^{E_7})\otimes_{{\cal B}_-}\mbb{C}u_0\cong U({\cal G}_+)\otimes_\mbb{C} \mbb{C}u_0\eqno(3.111)$$
be the induced ${\cal G}^{E_7}$-module.

Recall that $\mbb{N}$ is the set of nonnegative integers. Let
$${\cal
A}'=\mbb{C}[\ptl_{x_1},\ptl_{x_2},...,\ptl_{x_{27}}].\eqno(3.112)$$
We define an action of the associative algebra $\mbb{A}$ (cf.
(2.48)) on ${\cal A}'$ by
$$\ptl_{x_i}(\prod_{j=1}^{27}\ptl_{x_j}^{\be_j})=\ptl_{x_i}^{\be_i+1}\prod_{i\neq
j\in\ol{1,27}}\prod_{j=1}^{27}\ptl_{x_j}^{\be_j}\eqno(3.113)$$ and
$$x_i(\prod_{j=1}^{27}\ptl_{x_j}^{\be_j})=-\be_i\ptl_{x_i}^{\be_i-1}\prod_{i\neq
j\in\ol{1,27}}\prod_{j=1}^{27}\ptl_{x_j}^{\be_j}\eqno(3.114)$$ for
$i\in\ol{1,27}$. Since
$$[-x_i,\ptl_{x_j}]=[\ptl_{x_i},x_j]=\dlt_{i,j}\qquad\for\;\;i,j\in\ol{1,27},\eqno(3.115)$$
the above action gives an associative algebra representation of
$\mbb{A}$. Thus it also gives a Lie algebra representation of
$\mbb{A}$ (cf. (3.48)).
 It is straightforward to verify that
$$[d|_{{\cal A}'},\ptl|_{{\cal A}'}]=[d,\ptl]|_{{\cal
A}'}\qquad\for\;\;d\in{\cal C}_0,\;\ptl\in{\cal D}.\eqno(3.116)$$

Define linear map $\vs: \Psi\rta {\cal A}'$ by
$$\vs(\prod_{i=1}^{27}\xi_i^{\ell_i}\otimes
u_0)=\prod_{i=1}^{27}\ptl_{x_i}^{\ell_i}\qquad(\ell_1,...,\ell_{27})\in\mbb{N}^{27}.\eqno(3.117)$$
According to (2.50)-(2.87), (3.113) and (3.114),
$$D(1)=-27,\;\;d(1)=0\qquad\for\;\;d\in o(10,\mbb{C})|_{\cal A}.\eqno(3.118)$$
Moreover, (3.110), (3.111),  (3.113), (3.114) and (3.118) imply
$$\vs(\xi(v))=\vt(\xi)\vs(v)\qquad\for\;\;\xi\in{\cal
G}_0,\;v\in\Psi.\eqno(3.119)$$

Now (3.111) and (3.113) imply
$$\vs(w(u))=\vt(w)(\vs(u))\qquad\for\;\;w\in{\cal
B}_+,\;\;u\in\Psi.\eqno(3.120)$$ Thus  we have
$$\vs w|_{\Psi}\vs^{-1}=\vt(w)|_{{\cal A}'}\qquad\for\;\;w\in{\cal
B}_+.\eqno(3.121)$$ On the other hand, the linear map
$$\psi(v)=\vs v|_{\Psi}\vs^{-1}\qquad\for\;\; v\in{\cal
G}^{E_7}\eqno(3.122)$$ is a Lie algebra monomorphism from ${\cal
G}^{E_7}$ to $\mbb{A}|_{{\cal A}'}$. According to (3.76) and
(3.115),
$$\psi(\eta_1)=P_1|_{{\cal A}'}.\eqno(3.123)$$
By the constructions of $P_2,...,P_{27}$ in (2.77)-(2.102), we have
$$\psi(\eta_i)=P_i|_{{\cal
A}'}\qquad\for\;\;i\in\ol{2,27}.\eqno(3.124)$$ Therefore, we have
$$\psi(v)=\vt(v)|_{{\cal A}'}\qquad\for\;\;v\in{\cal
G}^{E_7}.\eqno(3.125)$$ In particular, ${\cal C}|_{{\cal
A}'}=\vt({\cal G}^{E_7})|_{{\cal A}'}=\psi({\cal G}^{E_7})$ forms a
Lie algebra. Since the linear map $d\mapsto d|_{{\cal A}'}$ for
$d\in {\cal C}$ is injective, we have that ${\cal C}$ forms a Lie
subalgebra of $\mbb{A}$ and $\vt$ is a Lie algebra
isomorphism.$\qquad\Box$ \psp

By the above theorem, a Lie group of type $E_7$ is generated by the
linear transformations $\{e^{bu|_{\cal A}}\mid b\in\mbb{R},\; u\in
{\cal G}^{E_6}\}$ associated with (2.50)-(2.87), the real
translations and dilations in $\sum_{i=1}^{27}\mbb{R}x_i$, and the
fractional transformations $\{e^{bP_i}\mid  b\in\mbb{R},\;
i\in\ol{1,27}\}$ such as
$$e^{bP_1}(x_i)=\frac{x_i}{1-bx_1},\;\;i\in\{\ol{1,13},15,16,18,20\},\eqno(3.126)$$
$$e^{bP_1}(x_{14})
=x_{14}-\frac{b(x_2x_{11}+x_3x_9+x_4x_7-x_5x_6)}{1-bx_1},\eqno(3.127)$$
$$e^{bP_1}(x_{17})
=x_{17}-\frac{b(x_2x_{13}+x_3x_{12}+x_4x_{10}-x_6x_8)}{1-bx_1},\eqno(3.128)$$
$$e^{bP_1}(x_{19})
=x_{19}-\frac{b(x_2x_{16}+x_3x_{15}+x_5x_{10}-x_7x_8)}{1-bx_1},\eqno(3.129)$$
$$e^{bP_1}(x_{21})
=x_{21}-\frac{b(x_2x_{18}-x_4x_{15}+x_5x_{12}
-x_8x_9)}{1-bx_1},\eqno(3.130)$$
$$e^{bP_1}(x_{22})
=x_{22}+\frac{b(x_3x_{18}+x_4x_{16}-x_5x_{13}
+x_8x_{11})}{1-bx_1},\eqno(3.131)$$
$$e^{bP_1}(x_{23})
=x_{23}-\frac{b(x_2x_{20}+x_6x_{15}-x_7x_{12}
+x_9x_{10})}{1-bx_1},\eqno(3.132)$$
$$e^{bP_1}(x_{24})
=x_{24}+\frac{b(x_3x_{20}-x_6x_{16}+x_7x_{13}
-x_{10}x_{11})}{1-bx_1},\eqno(3.133)$$
$$e^{bP_1}(x_{25})
=x_{25}-\frac{b(x_4x_{20}+x_6x_{18}-x_9x_{13}
+x_{11}x_{12})}{1-bx_1},\eqno(3.134)$$
$$e^{bP_1}(x_{26})
=x_{26}+\frac{b(x_5x_{20}+x_7x_{18}-x_9x_{16}+x_{11}x_{15})}{1-bx_1},\eqno(3.135)$$
$$e^{bP_1}(x_{27})
=x_{27}+\frac{b(x_8x_{20}+x_{10}x_{18}-x_{12}x_{16}+x_{13}x_{15}
)}{1-bx_1}\eqno(3.136)$$ by (3.1)-(3.7), (3.9), (3.11), (3.14) and
(3.58).\psp

Later, we need the following Dickson's ${\cal G}^{E_6}$-invariant:
\begin{eqnarray*}\chi&=&x_1(x_{17}x_{26}+x_{19}x_{25}-x_{21}x_{24}+x_{22}x_{23})
-x_{14}(x_1x_{27}-x_8x_{20}-x_{10}x_{18}+x_{12}x_{16} \\ &
&-x_{13}x_{15})
+x_2(x_{13}x_{26}+x_{16}x_{25}-x_{18}x_{24}+x_{20}x_{22})
-x_{11}(x_2x_{27}+x_8x_{23}+x_{10}x_{21}\\ &
&-x_{12}x_{19}+x_{15}x_{17})
+x_3(x_{12}x_{26}+x_{15}x_{25}-x_{18}x_{23}+x_{20}x_{21})
-x_9(x_3x_{27}-x_8x_{24}\\ &
&-x_{10}x_{22}+x_{13}x_{19}-x_{16}x_{17})
+x_4(x_{10}x_{26}+x_{15}x_{24}-x_{16}x_{23}+x_{19}x_{20})
-x_7(x_4x_{27}\\ &
&+x_8x_{25}+x_{12}x_{22}-x_{13}x_{21}+x_{17}x_{18})+x_5(x_{10}x_{25}-x_{12}x_{24}+x_{13}x_{23}-x_{17}x_{20})
\\ &
&+x_6(x_5x_{27}-x_8x_{26}+x_{15}x_{22}-x_{16}x_{21}+x_{18}x_{19}).
\hspace{5cm}(3.137)\end{eqnarray*} In [X2], we used partial
differential equations to prove that any ${\cal G}^{E_6}$-singular
vector in ${\cal A}$ is a polynomial in $x_1,\;\zeta_1,\;\chi$.
Moreover, $$\sum_{14\neq
i\in\ol{1,26}}x_i\zeta_{28-i}-x_{14}\zeta_{14}-x_{27}\zeta_1=-3\chi\eqno(3.138)$$
by Table 1, Table 2, (2.50)-(2.55), (3.1), (3.8)-(3.33) and
(3.35)-(3.40).

\section{Functor from $E_6$-Mod to $E_7$-Mod}

\quad \quad In this section, we construct a new functor from
$E_6$-{\bf Mod} to $E_7$-{\bf Mod}.

Note that
$${\cal G}^{E_6}_{\cal A}=(\bigoplus_{i=1}^6{\cal A}\al_i)\oplus\bigoplus_{\al\in
\Phi_{E_6}}{\cal A}E_{\al}\eqno(4.1)$$ forms a Lie algebra with the
Lie bracket:
$$[fu_1,gu_2]=fg[u_1,u_2]\qquad\for\;\;f,g\in{\cal A},\;u_1,u_2\in
{\cal G}^{E_6}.\eqno(4.2)$$ Moreover, we define the Lie algebra
$${\cal K}={\cal G}^{E_6}_{\cal A}\oplus {\cal A}\kappa\eqno(4.3)$$
with the Lie bracket:
$$[\xi_1+f\kappa,\xi_2+g\kappa]=[\xi_1,\xi_2]\qquad\for\;\;\xi_1,\xi_2\in
{\cal G}^{E_6}_{\cal A},\;f,g\in{\cal A}.\eqno(4.4)$$ Similarly,
$gl(27,{\cal A})$ becomes a Lie algebra with the Lie bracket as that
in (4.2). Recall the Witt algebra ${\cal
W}_{27}=\sum_{i=1}^{27}{\cal A}\ptl_{x_i}$, and Shen [Sg1-3] found a
monomorphism $\Im$ from the Lie algebra ${\cal W}_{27}$ to the Lie
algebra of semi-product ${\cal W}_{27}+gl(27,{\cal A})$ defined by
$$\Im(\sum_{i=1}^{27}f_i\ptl_{x_i})=\sum_{i=1}^{27}f_i\ptl_{x_i}+\Im_1(\sum_{i=1}^{27}f_i\ptl_{x_i}),\;\;
\Im_1(\sum_{i=1}^{27}f_i\ptl_{x_i})=\sum_{i,j=1}^{27}\ptl_{x_i}(f_j)E_{i,j}.
\eqno(4.5)$$ According to our construction of $P_1$-$P_{27}$ in
(3.58) and (3.77)-(3.102),
$$\Im_1(P_i)=\sum_{r=1}^{27}x_r\Im_1([\xi_r,\eta_i]|_{\cal A})\qquad\for\;\;i\in\ol{1,27}.\eqno(4.6)$$

On the other hand,
$$\widehat{\cal K}={\cal W}_{27}\oplus {\cal K}\eqno(4.7)$$
becomes a Lie algebra with the Lie bracket
\begin{eqnarray*}&&[d_1+f_1u_1+f_2\kappa,d_2+g_1u_2+g_2\kappa]\\
&=&[d_1,d_2]+f_1g_1[u_1,u_2]+d_1(g_2)u_2
-d_2(g_1)u_1+(d_1(g_2)-d_2(g_1))\kappa\hspace{2.8cm}(4.8)\end{eqnarray*}
for $f_1,f_2,g_1,g_2\in{\cal A},\;u_1,u_2\in {\cal G}^{E_6}$ and
$d_1,d_2\in{\cal W}_{27}$. Note
$${\cal G}_0={\cal
G}^{E_6}\oplus\mbb{C}\widehat\al\eqno(4.9)$$ by (2.41) and (2.88).
So there exists a Lie algebra monomorphism $\varrho:{\cal G}_0\rta
{\cal K}$ determined by
$$\varrho(\widehat\al)=3\kappa,\;\;\varrho(u)=u\qquad\for\;\;u\in{\cal
G}^{E_6}.\eqno(4.10)$$ Since $\Im$ is a Lie algebra monomorphism,
our construction of $P_1$-$P_{27}$ in (3.48) and (3.58)-(3.102) show
that we have a Lie algebra monomorphism $\iota: {\cal G}^{E_7}\rta
\widehat{\cal K}$ given by
$$\iota(u)=u|_{\cal A}+\varrho(u)\qquad\for\;\;u\in{\cal
G}_0,\eqno(4.11)$$
$$\iota(\xi_i)=\ptl_{x_i},\;\;\iota(\eta_i)=P_i+\sum_{r=1}^{27}x_r\varrho([\xi_r,\eta_i])\qquad\for
\;\;i\in\ol{1,27}.\eqno(4.12)$$

According to (2.4), (2.11) and (2.34)-(2.40),
\begin{eqnarray*} \iota(\eta_1)&=&P_1+x_2E_{\al_6}
+x_3E_{(0,0,0,0,1,1)}+x_4E_{(0,0,0,1,1,1)}
+x_5E_{(0,0,1,1,1,1)}+x_6E_{(0,1,0,1,1,1)}\\
&
&-x_1\varrho(\al_7)+x_7E_{(0,1,1,1,1,1)}+x_8E_{(1,0,1,1,1,1)}+x_9E_{(0,1,1,2,1,1)}+x_{10}E_{(1,1,1,1,1,1)}
\\& & +x_{11}E_{(0,1,1,2,2,1)}+x_{12}E_{(1,1,1,2,1,1)}+x_{13}E_{(1,1,1,2,2,1)}+x_{15}E_{(1,1,2,2,1,1)}
\\& &
+x_{16}E_{(1,1,2,2,2,1)}+x_{18}E_{(1,1,2,3,2,1)}+x_{20}E_{(1,2,2,3,2,1)}.
\hspace{4.3cm}(4.13)\end{eqnarray*} Moreover, (3.7) yields
\begin{eqnarray*}
\iota(\eta_2)&=&-\iota([E_{-\al_6},\eta_1])=-[\iota(E_{-\al_6}),\iota(\eta_1)]=
-[E_{-\al_6}|_{\cal A}+E_{-\al_6},\iota(\eta_1)]\\
&=&P_2-x_1E_{-\al_6} +x_3E_{\al_5}+x_4E_{(0,0,0,1,1)}
+x_5E_{(0,0,1,1,1)}+x_6E_{(0,1,0,1,1)}\\
&&-x_2\varrho(\al_6+\al_7)+x_7E_{(0,1,1,1,1)}+x_8E_{(1,0,1,1,1)}+x_9E_{(0,1,1,2,1)}+x_{10}E_{(1,1,1,1,1)}
\\& & +x_{12}E_{(1,1,1,2,1)}-x_{14}E_{(0,1,1,2,2,1)}+x_{15}E_{(1,1,2,2,1)}-x_{17}E_{(1,1,1,2,2,1)}
\\& &
-x_{19}E_{(1,1,2,2,2,1)}-x_{21}E_{(1,1,2,3,2,1)}-x_{23}E_{(1,2,2,3,2,1)}.
\hspace{4.3cm}(4.14)\end{eqnarray*}  Expression (3.6) gives
\begin{eqnarray*}
\iota(\eta_3)&=&-\iota([E_{-\al_5},\eta_2])=-[\iota(E_{-\al_5}),\iota(\eta_2)]=
-[E_{-\al_5}|_{\cal A}+E_{-\al_5},\iota(\eta_2)]\\
&=&P_3-x_1E'_{(0,0,0,0,1,1)} -x_2E_{-\al_5}+x_4E_{\al_4}
+x_5E_{(0,0,1,1)}+x_6E_{(0,1,0,1)}\\
&&-x_3\varrho(\al_5+\al_6+\al_7)+x_7E_{(0,1,1,1)}+x_8E_{(1,0,1,1)}+x_{10}E_{(1,1,1,1)}-x_{11}E_{(0,1,1,2,1)}
\\& & -x_{13}E_{(1,1,1,2,1)}-x_{14}E_{(0,1,1,2,1,1)}-x_{16}E_{(1,1,2,2,1)}-x_{17}E_{(1,1,1,2,1,1)}
\\& &
-x_{19}E_{(1,1,2,2,1,1)}+x_{22}E_{(1,1,2,3,2,1)}+x_{24}E_{(1,2,2,3,2,1)}.
\hspace{4.3cm}(4.15)\end{eqnarray*} Furthermore,
\begin{eqnarray*}
\iota(\eta_4)&=&-\iota([E_{-\al_4},\eta_3])=-[\iota(E_{-\al_4}),\iota(\eta_3)]=
-[E_{-\al_4}|_{\cal A}+E_{-\al_4},\iota(\eta_3)]\\
&=&P_4-x_1E'_{(0,0,0,1,1,1)} -x_2E'_{(0,0,0,1,1)}-x_3E_{-\al_4}
+x_5E_{\al_3}+x_6E_{\al_2}\\
&&-x_4\varrho(\sum_{i=4}^7\al_i)-x_9E_{(0,1,1,1)}+x_8E_{(1,0,1)}-x_{12}E_{(1,1,1,1)}-x_{11}E_{(0,1,1,1,1)}
\\& & -x_{13}E_{(1,1,1,1,1)}-x_{14}E_{(0,1,1,1,1,1)}-x_{17}E_{(1,1,1,1,1,1)}+x_{18}E_{(1,1,2,2,1)}
\\& &
+x_{21}E_{(1,1,2,2,1,1)}+x_{22}E_{(1,1,2,2,2,1)}-x_{25}E_{(1,2,2,3,2,1)}
\hspace{4.3cm}(4.16)\end{eqnarray*} by (3.5),
\begin{eqnarray*}
\iota(\eta_5)&=&-\iota([E_{-\al_3},\eta_4])=-[\iota(E_{-\al_3}),\iota(\eta_4)]=
-[E_{-\al_3}|_{\cal A}+E_{-\al_3},\iota(\eta_4)]\\
&=&P_5-x_1E'_{(0,0,1,1,1,1)} -x_2E'_{(0,0,1,1,1)}-x_3E'_{(0,0,1,1)}
-x_4E_{-\al_3}+x_7E_{\al_2}\\
&&-x_5\varrho(\sum_{i=3}^7\al_i)+x_8E_{\al_1}+x_9E_{(0,1,0,1)}-x_{15}E_{(1,1,1,1)}+x_{11}E_{(0,1,0,1,1)}
\\& & -x_{16}E_{(1,1,1,1,1)}+x_{14}E_{(0,1,0,1,1,1)}-x_{19}E_{(1,1,1,1,1,1)}-x_{18}E_{(1,1,1,2,1)}
\\& &
-x_{21}E_{(1,1,1,2,1,1)}-x_{22}E_{(1,1,1,2,2,1)}+x_{26}E_{(1,2,2,3,2,1)}
\hspace{4.3cm}(4.17)\end{eqnarray*} by (3.4),
\begin{eqnarray*}
\iota(\eta_6)&=&-\iota([E_{-\al_2},\eta_4])=-[\iota(E_{-\al_2}),\iota(\eta_4)]=
-[E_{-\al_2}|_{\cal A}+E_{-\al_2},\iota(\eta_4)]\\
&=&P_6-x_1E'_{(0,1,0,1,1,1)} -x_2E'_{(0,1,0,1,1)}-x_3E'_{(0,1,0,1)}
-x_4E_{-\al_2}+x_7E_{\al_3}\\
&&-x_6\varrho(\sum_{i=2,4,5,6,7}\al_i)+x_9E_{(0,0,1,1)}+x_{10}E_{(1,0,1)}+x_{12}E_{(1,0,1,1)}+x_{11}E_{(0,0,1,1,1)}
\\& & +x_{13}E_{(1,0,1,1,1)}+x_{14}E_{(0,0,1,1,1,1)}+x_{17}E_{(1,0,1,1,1,1)}-x_{20}E_{(1,1,2,2,1)}
\\& &
-x_{23}E_{(1,1,2,2,1,1)}-x_{24}E_{(1,1,2,2,2,1)}-x_{25}E_{(1,1,2,3,2,1)}
\hspace{4.4cm}(4.18)\end{eqnarray*}
 by (3.3),
\begin{eqnarray*}
\iota(\eta_7)&=&-\iota([E_{-\al_2},\eta_5])=-[\iota(E_{-\al_2}),\iota(\eta_5)]=
-[E_{-\al_2}|_{\cal A}+E_{-\al_2},\iota(\eta_5)]\\
&=&P_7-x_1E'_{(0,1,1,1,1,1)} -x_2E'_{(0,1,1,1,1)}-x_3E'_{(0,1,1,1)}
-x_5E_{-\al_2}-x_6E_{-\al_3}\\
&&-x_7\varrho(\sum_{i=2}^7\al_i)-x_9E_{\al_4}+x_{10}E_{\al_1}+x_{15}E_{(1,0,1,1)}-x_{11}E_{(0,0,0,1,1)}
\\& & +x_{16}E_{(1,0,1,1,1)}-x_{14}E_{(0,0,0,1,1,1)}+x_{19}E_{(1,0,1,1,1,1)}+x_{20}E_{(1,1,1,2,1)}
\\& &
+x_{23}E_{(1,1,1,2,1,1)}+x_{24}E_{(1,1,1,2,2,1)}+x_{26}E_{(1,1,2,3,2,1)}
\hspace{4.3cm}(4.19)\end{eqnarray*} by (3.3),
\begin{eqnarray*}
\iota(\eta_8)&=&-\iota([E_{-\al_1},\eta_5])=-[\iota(E_{-\al_1}),\iota(\eta_5)]=
-[E_{-\al_1}|_{\cal A}+E_{-\al_1},\iota(\eta_5)]\\
&=&P_8-x_1E'_{(1,0,1,1,1,1)} -x_2E'_{(1,0,1,1,1)}-x_3E'_{(1,0,1,1)}
-x_4E'_{(1,0,1)}-x_5E_{-\al_1}\\
&&-x_8\varrho(\sum_{2\neq
i\in\ol{1,7}}\al_i)+x_{10}E_{\al_2}+x_{12}E_{(0,1,0,1)}+x_{15}E_{(0,1,1,1)}+x_{13}E_{(0,1,0,1,1)}
\\& & +x_{16}E_{(0,1,1,1,1)}+x_{17}E_{(0,1,0,1,1,1)}+x_{19}E_{(0,1,1,1,1,1)}+x_{18}E_{(0,1,1,2,1)}
\\& &
+x_{21}E_{(0,1,1,2,1,1)}+x_{22}E_{(0,1,1,2,2,1)}+x_{27}E_{(1,2,2,3,2,1)}
\hspace{4.3cm}(4.20)\end{eqnarray*} by (3.2),
\begin{eqnarray*}
\iota(\eta_9)&=&\iota([E_{-\al_4},\eta_7])=[\iota(E_{-\al_4}),\iota(\eta_7)]=
[E_{-\al_4}|_{\cal A}+E_{-\al_4},\iota(\eta_7)]\\
&=&P_9-x_1E'_{(0,1,1,2,1,1)} -x_2E'_{(0,1,1,2,1)}+x_4E'_{(0,1,1,1)}
-x_5E'_{(0,1,0,1)}-x_6E'_{(0,0,1,1)}\\
&&-x_9\varrho(\al_4+\sum_{i=2}^7\al_i)+x_7E_{-\al_4}+x_{12}E_{\al_1}-x_{15}E_{(1,0,1)}-x_{11}E_{\al_5}
\\& & +x_{18}E_{(1,0,1,1,1)}-x_{14}E_{(0,0,0,0,1,1)}+x_{21}E_{(1,0,1,1,1,1)}-x_{20}E_{(1,1,1,1,1)}
\\& &
-x_{23}E_{(1,1,1,1,1,1)}+x_{25}E_{(1,1,1,2,2,1)}-x_{26}E_{(1,1,2,2,2,1)}
\hspace{4.3cm}(4.21)\end{eqnarray*} by (3.5),
\begin{eqnarray*}
\iota(\eta_{10})&=&-\iota([E_{-\al_1},\eta_7])=-[\iota(E_{-\al_1}),\iota(\eta_7)]=
-[E_{-\al_1}|_{\cal A}+E_{-\al_1},\iota(\eta_7)]\\
&=&P_{10}-x_1E'_{(1,1,1,1,1,1)}
-x_2E'_{(1,1,1,1,1)}-x_3E'_{(1,1,1,1)}
-x_8E_{-\al_2}-x_6E'_{(1,0,1)}\\
&&-x_{10}\varrho(\sum_{i=1}^7\al_i)-x_7E_{-\al_1}-x_{12}E_{\al_4}-x_{15}E_{(0,0,1,1)}-x_{13}E_{(0,0,0,1,1)}
\\& & -x_{16}E_{(0,0,1,1,1)}-x_{17}E_{(0,0,0,1,1,1)}-x_{19}E_{(0,0,1,1,1,1)}-x_{20}E_{(0,1,1,2,1)}
\\& &
-x_{23}E_{(0,1,1,2,1,1)}-x_{24}E_{(0,1,1,2,2,1)}+x_{27}E_{(1,1,2,3,2,1)}
\hspace{4.3cm}(4.22)\end{eqnarray*} by (3.2),
\begin{eqnarray*}
\iota(\eta_{11})&=&\iota([E_{-\al_5},\eta_9])=[\iota(E_{-\al_5}),\iota(\eta_9)]=
[E_{-\al_5}|_{\cal A}+E_{-\al_5},\iota(\eta_9)]\\
&=&P_{11}-x_1E'_{(0,1,1,2,2,1)}
+x_3E'_{(0,1,1,2,1)}+x_4E'_{(0,1,1,1,1)}
-x_5E'_{(0,1,0,1,1)}\\
&&-x_{11}\varrho(\al_4+\al_5+\sum_{i=2}^7\al_i)-x_6E'_{(0,0,1,1,1)}+x_7E'_{(0,0,0,1,1)}+x_9E_{-\al_5}\\&
&+x_{13}E_{\al_1}-x_{16}E_{(1,0,1)}
 -x_{18}E_{(1,0,1,1)}-x_{14}E_{\al_6}+x_{20}E_{(1,1,1,1)}\\& &
+x_{22}E_{(1,0,1,1,1,1)}-x_{24}E_{(1,1,1,1,1,1)}-x_{25}E_{(1,1,1,2,1,1)}+x_{26}E_{(1,1,2,2,1,1)}
\hspace{1.3cm}(4.23)\end{eqnarray*} by (3.6),
\begin{eqnarray*}
\iota(\eta_{12})&=&\iota([E_{-\al_4},\eta_{10}])=[\iota(E_{-\al_4}),\iota(\eta_{10})]=
[E_{-\al_4}|_{\cal A}+E_{-\al_4},\iota(\eta_{10})]\\
&=&P_{12}-x_1E'_{(1,1,1,2,1,1)}
-x_2E'_{(1,1,1,2,1)}+x_4E'_{(1,1,1,1)}-x_6E'_{(1,0,1,1)}
\\
&&-x_{12}\varrho(\al_4+\sum_{i=1}^7\al_i)-x_8E'_{(0,1,0,1)}-x_9E_{-\al_1}+x_{10}E_{-\al_4}\\&
&-x_{13}E_{\al_5}+x_{15}E_{\al_3}
 -x_{17}E_{(0,0,0,0,1,1)}-x_{18}E_{(0,0,1,1,1)}+x_{20}E_{(0,1,1,1,1)}\\& &-x_{21}E_{(0,0,1,1,1,1)}
+x_{23}E_{(0,1,1,1,1,1)}-x_{25}E_{(0,1,1,2,2,1)}-x_{27}E_{(1,1,2,2,2,1)}
\hspace{1.3cm}(4.24)\end{eqnarray*}
 by (3.5),
\begin{eqnarray*}
\iota(\eta_{13})&=&\iota([E_{-\al_5},\eta_{12}])=[\iota(E_{-\al_5}),\iota(\eta_{12})]=
[E_{-\al_5}|_{\cal A}+E_{-\al_5},\iota(\eta_{12})]\\
&=&P_{13}-x_1E'_{(1,1,1,2,2,1)}
+x_3E'_{(1,1,1,2,1)}+x_4E'_{(1,1,1,1,1)}
-x_8E'_{(0,1,0,1,1)}\\
&&-x_{13}\varrho(\al_4+\al_5+\sum_{i=1}^7\al_i)-x_6E'_{(1,0,1,1,1)}+x_{10}E'_{(0,0,0,1,1)}-x_{11}E_{-\al_1}\\&
&+x_{12}E_{-\al_5}+x_{16}E_{\al_3}
 -x_{17}E_{\al_6}+x_{18}E_{(0,0,1,1)}-x_{20}E_{(0,1,1,1)}\\& &-x_{22}E_{(0,0,1,1,1,1)}
+x_{24}E_{(0,1,1,1,1,1)}+x_{25}E_{(0,1,1,2,1,1)}+x_{27}E_{(1,1,2,2,1,1)}
\hspace{1.3cm}(4.25)\end{eqnarray*} by (3.6),
\begin{eqnarray*}
\iota(\eta_{14})&=&\iota([E_{-\al_6},\eta_{11}])=[\iota(E_{-\al_6}),\iota(\eta_{11})]=
[E_{-\al_6}|_{\cal A}+E_{-\al_6},\iota(\eta_{11})]\\
&=&P_{14}+x_2E'_{(0,1,1,2,2,1)}
+x_3E'_{(0,1,1,2,1,1)}+x_4E'_{(0,1,1,1,1,1)}
-x_5E'_{(0,1,0,1,1,1)}\\
&&-x_{14}\varrho(\sum_{r=4}^6\al_r+\sum_{i=2}^7\al_i)-x_6E'_{(0,0,1,1,1,1)}+x_7E'_{(0,0,0,1,1,1)}+x_9E'_{(0,0,0,0,1,1)}\\&
&+x_{11}E_{-\al_6}+x_{17}E_{\al_1}
 -x_{19}E_{(1,0,1)}-x_{21}E_{(1,0,1,1)}-x_{22}E_{(1,0,1,1,1)}\\& &+x_{23}E_{(1,1,1,1)}
+x_{24}E_{(1,1,1,1,1)}+x_{25}E_{(1,1,1,2,1)}-x_{26}E_{(1,1,2,2,1)}
\hspace{2.5cm}(4.26)\end{eqnarray*} by (3.7),
\begin{eqnarray*}
\iota(\eta_{15})&=&-\iota([E_{-\al_3},\eta_{12}])=-[\iota(E_{-\al_2}),\iota(\eta_{12})]=
-[E_{-\al_3}|_{\cal A}+E_{-\al_3},\iota(\eta_{12})]\\
&=&P_{15}-x_1E'_{(1,1,2,2,1,1)}
-x_2E'_{(1,1,2,2,1)}+x_5E'_{(1,1,1,1)}
-x_7E'_{(1,0,1,1)}\\
&&-x_{15}\varrho(\al_3+\al_4+\sum_{i=1}^7\al_i)-x_8E'_{(0,1,1,1)}+x_9E'_{(1,0,1)}+x_{10}E'_{(0,0,1,1)}\\&
&-x_{12}E_{-\al_3}-x_{16}E_{\al_5} +x_{18}E_{(0,0,0,1,1)}
-x_{19}E_{(0,0,0,0,1,1)}-x_{20}E_{(0,1,0,1,1)}\\&
&+x_{21}E_{(0,0,0,1,1,1)}
-x_{23}E_{(0,1,0,1,1,1)}+x_{26}E_{(0,1,1,2,2,1)}+x_{27}E_{(1,1,1,2,2,1)}
\hspace{1.3cm}(4.27)\end{eqnarray*}
 by (3.4),
\begin{eqnarray*}\iota(\eta_{16})&=&\iota([E_{-\al_5},\eta_{15}])=[\iota(E_{-\al_5}),\iota(\eta_{15})]=
[E_{-\al_5}|_{\cal A}+E_{-\al_5},\iota(\eta_{15})]\\
&=&P_{16}-x_1E'_{(1,1,2,2,2,1)}
+x_3E'_{(1,1,2,2,1)}+x_5E'_{(1,1,1,1,1)}
-x_7E'_{(1,0,1,1,1)}\\
&&-x_{16}\varrho(\sum_{r=3}^5\al_r+\sum_{i=1}^7\al_i)-x_8E'_{(0,1,1,1,1)}+x_{10}E'_{(0,0,1,1,1)}+x_{11}E'_{(1,0,1)}\\&
&-x_{13}E_{-\al_3}+x_{15}E_{-\al_5}
 -x_{18}E_{\al_4}-x_{19}E_{\al_6}+x_{20}E_{(0,1,0,1)}+x_{22}E_{(0,0,0,1,1,1)}\\& &
-x_{24}E_{(0,1,0,1,1,1)}-x_{26}E_{(0,1,1,2,1,1)}-x_{27}E_{(1,1,1,2,1,1)}
\hspace{4.3cm}(4.28)\end{eqnarray*}
 by (3.6),
\begin{eqnarray*}
\iota(\eta_{17})&=&\iota([E_{-\al_6},\eta_{13}])=[\iota(E_{-\al_6}),\iota(\eta_{13})]=
[E_{-\al_6}|_{\cal A}+E_{-\al_6},\iota(\eta_{13})]\\
&=&P_{17}+x_2E'_{(1,1,1,2,2,1)}
+x_3E'_{(1,1,1,2,1,1)}+x_4E'_{(1,1,1,1,1,1)}
-x_6E'_{(1,0,1,1,1,1)}\\
&&-x_{17}\varrho(\sum_{r=4}^6\al_r+\sum_{i=1}^7\al_i)-x_8E'_{(0,1,0,1,1,1)}+x_{10}E'_{(0,0,0,1,1,1)}-x_{14}E_{-\al_1}\\&
&+x_{12}E'_{(0,0,0,0,1,1)}+x_{13}E_{-\al_6}+x_{19}E_{\al_3}
 +x_{21}E_{(0,0,1,1)}-x_{23}E_{(0,1,1,1)}\\& &+x_{22}E_{(0,0,1,1,1)}
-x_{24}E_{(0,1,1,1,1)}-x_{25}E_{(0,1,1,2,1)}-x_{27}E_{(1,1,2,2,1)}
\hspace{2.2cm}(4.29)\end{eqnarray*} by (3.7)
\begin{eqnarray*}\iota(\eta_{18})&=&\iota([E_{-\al_4},\eta_{16}])=[\iota(E_{-\al_4}),\iota(\eta_{16})]=
[E_{-\al_4}|_{\cal A}+E_{-\al_4},\iota(\eta_{16})]\\
&=&P_{18}-x_1E'_{(1,1,2,3,2,1)}
-x_4E'_{(1,1,2,2,1)}+x_5E'_{(1,1,1,2,1)}
-x_8E'_{(0,1,1,2,1)}\\
&&-x_{18}\varrho(\al_4+\sum_{r=3}^5\al_r+\sum_{i=1}^7\al_i)-x_9E'_{(1,0,1,1,1)}+x_{11}E'_{(1,0,1,1)}+x_{12}E'_{(0,0,1,1,1)}\\&
&-x_{13}E'_{(0,0,1,1)}-x_{15}E'_{(0,0,0,1,1)}
 +x_{16}E_{-\al_4}-x_{21}E_{\al_6}-x_{20}E_{\al_2}+x_{22}E_{(0,0,0,0,1,1)}\\& &
-x_{25}E_{(0,1,0,1,1,1)}+x_{26}E_{(0,1,1,1,1,1)}+x_{27}E_{(1,1,1,1,1,1)}
\hspace{4.3cm}(4.30)\end{eqnarray*}
 by (3.5),
\begin{eqnarray*}\iota(\eta_{19})&=&\iota([E_{-\al_6},\eta_{16}])=[\iota(E_{-\al_6}),\iota(\eta_{16})]=
[E_{-\al_6}|_{\cal A}+E_{-\al_6},\iota(\eta_{16})]\\
&=&P_{19}+x_2E'_{(1,1,2,2,2,1)}
+x_3E'_{(1,1,2,2,1,1)}+x_5E'_{(1,1,1,1,1,1)}
-x_7E'_{(1,0,1,1,1,1)}\\
&&-x_{19}\varrho(\sum_{r=3}^6\al_r+\sum_{i=1}^7\al_i)-x_8E'_{(0,1,1,1,1,1)}+x_{10}E'_{(0,0,1,1,1,1)}+x_{14}E'_{(1,0,1)}\\&
&+x_{15}E'_{(0,0,0,0,1,1)}+x_{16}E_{-\al_6}-x_{17}E_{-\al_3}
 -x_{21}E_{\al_4}+x_{23}E_{(0,1,0,1)}\\& &-x_{22}E_{(0,0,0,1,1)}
+x_{24}E_{(0,1,0,1,1)}+x_{26}E_{(0,1,1,2,1)}+x_{27}E_{(1,1,1,2,1)}
\hspace{2.3cm}(4.31)\end{eqnarray*}
 by (3.7),
\begin{eqnarray*}\iota(\eta_{20})&=&\iota([E_{-\al_2},\eta_{18}])=[\iota(E_{-\al_2}),\iota(\eta_{18})]=
[E_{-\al_2}|_{\cal A}+E_{-\al_2},\iota(\eta_{18})]\\
&=&P_{20}-x_1E'_{(1,2,2,3,2,1)}
+x_6E'_{(1,1,2,2,1)}-x_7E'_{(1,1,1,2,1)}
+x_{10}E'_{(0,1,1,2,1)}\\
&&-x_{20}\varrho(\al_4+\sum_{r=2}^5\al_r+\sum_{i=1}^7\al_i)+x_9E'_{(1,1,1,1,1)}-x_{11}E'_{(1,1,1,1)}-x_{12}E'_{(0,1,1,1,1)}\\&
&+x_{13}E'_{(0,1,1,1)}+x_{15}E'_{(0,1,0,1,1)}
 -x_{16}E_{(0,1,0,1)}+x_{18}E_{-\al_2}-x_{23}E_{\al_6}\\& &+x_{24}E_{(0,0,0,0,1,1)}
-x_{25}E_{(0,0,0,1,1,1)}+x_{26}E_{(0,0,1,1,1,1)}+x_{27}E_{(1,0,1,1,1,1)}
\hspace{1.4cm}(4.32)\end{eqnarray*}
 by (3.3),
\begin{eqnarray*}\iota(\eta_{21})&=&\iota([E_{-\al_4},\eta_{19}])=[\iota(E_{-\al_4}),\iota(\eta_{19})]=
[E_{-\al_4}|_{\cal A}+E_{-\al_4},\iota(\eta_{19})]\\
&=&P_{21}+x_2E'_{(1,1,2,3,2,1)}
-x_4E'_{(1,1,2,2,1,1)}+x_5E'_{(1,1,1,2,1,1)}
-x_8E'_{(0,1,1,2,1,1)}\\
&&-x_{21}\varrho(\al_4+\sum_{r=3}^6\al_r+\sum_{i=1}^7\al_i)-x_9E'_{(1,0,1,1,1,1)}+x_{12}E'_{(0,0,1,1,1,1)}+x_{14}E'_{(1,0,1,1)}\\&
&-x_{15}E'_{(0,0,0,1,1,1)}-x_{17}E'_{(0,0,1,1)}+x_{18}E_{-\al_6}
 +x_{19}E_{-\al_4}-x_{23}E_{\al_2}\\& &-x_{22}E_{\al_5}
+x_{25}E_{(0,1,0,1,1)}-x_{26}E_{(0,1,1,1,1)}-x_{27}E_{(1,1,1,1,1)}
\hspace{3.3cm}(4.33)\end{eqnarray*}
 by (3.5),
\begin{eqnarray*}\iota(\eta_{22})&=&\iota([E_{-\al_5},\eta_{21}])=[\iota(E_{-\al_5}),\iota(\eta_{21})]=
[E_{-\al_5}|_{\cal A}+E_{-\al_5},\iota(\eta_{21})]\\
&=&P_{22}-x_3E'_{(1,1,2,3,2,1)} -x_4E'_{(1,1,2,2,2,1)}
-x_{22}\varrho(\al_4+\al_5+\sum_{r=3}^6\al_r+\sum_{i=1}^7\al_i)\\
&&+x_5E'_{(1,1,1,2,2,1)}-x_8E'_{(0,1,1,2,2,1)}-x_{11}E'_{(1,0,1,1,1,1)}+x_{13}E'_{(0,0,1,1,1,1)}+x_{14}E'_{(1,0,1,1,1)}\\&
&-x_{16}E'_{(0,0,0,1,1,1)}-x_{17}E'_{(0,0,1,1.1)}-x_{18}E'_{(0,0,0,0,1,1)}
 +x_{19}E'_{(0,0,0,1,1)}\\& &+x_{21}E_{-\al_5}-x_{24}E_{\al_2}
-x_{25}E_{(0,1,0,1)}+x_{26}E_{(0,1,1,1)}+x_{27}E_{(1,1,1,1)}
\hspace{2cm}(4.34)\end{eqnarray*}
 by (3.6),
\begin{eqnarray*}\iota(\eta_{23})&=&\iota([E_{-\al_2},\eta_{21}])=[\iota(E_{-\al_2}),\iota(\eta_{21})]=
[E_{-\al_2}|_{\cal A}+E_{-\al_2},\iota(\eta_{21})]\\
&=&P_{23}+x_2E'_{(1,2,2,3,2,1)}
+x_6E'_{(1,1,2,2,1,1)}-x_7E'_{(1,1,1,2,1,1)}
+x_9E'_{(1,1,1,1,1,1)}\\
&&-x_{23}\varrho(\al_4+\sum_{r=2}^6\al_r+\sum_{i=1}^7\al_i)+x_{10}E'_{(0,1,1,2,1,1)}-x_{12}E'_{(0,1,1,1,1,1)}-x_{14}E'_{(1,1,1,1)}\\&
&+x_{15}E'_{(0,1,0,1,1,1)}+x_{17}E'_{(0,1,1,1)}+x_{20}E_{-\al_6}
 -x_{19}E_{(0,1,0,1)}+x_{21}E_{-\al_2}\\& &-x_{24}E_{\al_5}
+x_{25}E_{(0,0,0,1,1)}-x_{26}E_{(0,0,1,1,1)}-x_{27}E_{(1,0,1,1,1)}
\hspace{3.3cm}(4.35)\end{eqnarray*}
 by (3.3),
\begin{eqnarray*}\iota(\eta_{24})&=&\iota([E_{-\al_2},\eta_{22}])=[\iota(E_{-\al_2}),\iota(\eta_{22})]=
[E_{-\al_2}|_{\cal A}+E_{-\al_2},\iota(\eta_{22})]\\
&=&P_{24}-x_3E'_{(1,2,2,3,2,1)} +x_6E'_{(1,1,2,2,2,1)}
-x_{24}\varrho(\al_4+\al_5+\sum_{r=2}^6\al_r+\sum_{i=1}^7\al_i)\\
&&-x_7E'_{(1,1,1,2,2,1)}+x_{10}E'_{(0,1,1,2,2,1)}+x_{11}E'_{(1,1,1,1,1,1)}-x_{13}E'_{(0,1,1,1,1,1)}\\&
&-x_{14}E'_{(1,1,1,1,1)}+x_{16}E'_{(0,1,0,1,1,1)}+x_{17}E'_{(0,1,1,1,1)}-x_{19}E'_{(0,1,0,1,1)}-x_{20}E'_{(0,0,0,0,1,1)}
 \\& &+x_{22}E_{-\al_2}+x_{23}E_{-\al_5}
-x_{25}E_{\al_4}+x_{26}E_{(0,0,1,1)}+x_{27}E_{(1,0,1,1)}
\hspace{2.5cm}(4.36)\end{eqnarray*}
 by (3.3),
\begin{eqnarray*}\iota(\eta_{25})&=&\iota([E_{-\al_4},\eta_{24}])=[\iota(E_{-\al_4}),\iota(\eta_{24})]=
[E_{-\al_4}|_{\cal A}+E_{-\al_4},\iota(\eta_{24})]\\
&=&P_{25}+x_4E'_{(1,2,2,3,2,1)} +x_6E'_{(1,1,2,3,2,1)}
-x_{25}\varrho(2\al_4+\al_5+\sum_{r=2}^6\al_r+\sum_{i=1}^7\al_i)\\
&&-x_9E'_{(1,1,1,2,2,1)}+x_{12}E'_{(0,1,1,2,2,1)}+x_{11}E'_{(1,1,1,2,1,1)}-x_{13}E'_{(0,1,1,2,1,1)}\\&
&-x_{14}E'_{(1,1,1,2,1)}+x_{17}E'_{(0,1,1,2,1)}+x_{18}E'_{(0,1,0,1,1,1)}+x_{20}E'_{(0,0,0,1,1,1)}-x_{21}E'_{(0,1,0,1,1)}
 \\& &+x_{22}E'_{(0,1,0,1)}-x_{23}E'_{(0,0,0,1,1)}
+x_{24}E_{-\al_4}-x_{26}E_{\al_3}-x_{27}E_{(1,0,1)}
\hspace{2cm}(4.37)\end{eqnarray*}
 by (3.5),
\begin{eqnarray*}\iota(\eta_{26})&=&\iota([E_{-\al_3},\eta_{25}])=[\iota(E_{-\al_3}),\iota(\eta_{25})]=
[E_{-\al_3}|_{\cal A}+E_{-\al_3},\iota(\eta_{25})]\\
&=&P_{26}-x_5E'_{(1,2,2,3,2,1)} -x_7E'_{(1,1,2,3,2,1)}
-x_{26}\varrho(\al_4+\sum_{s=3}^5+\sum_{r=2}^6\al_r+\sum_{i=1}^7\al_i)\\
&&+x_9E'_{(1,1,2,2,2,1)}-x_{15}E'_{(0,1,1,2,2,1)}-x_{11}E'_{(1,1,2,2,1,1)}+x_{16}E'_{(0,1,1,2,1,1)}\\&
&+x_{14}E'_{(1,1,2,2,1)}-x_{19}E'_{(0,1,1,2,1)}-x_{18}E'_{(0,1,1,1,1,1)}-x_{20}E'_{(0,0,1,1,1,1)}+x_{21}E'_{(0,1,1,1,1)}
 \\& &-x_{22}E'_{(0,1,1,1)}+x_{23}E'_{(0,0,1,1,1)}
-x_{24}E'_{(0,0,1,1)}+x_{25}E_{-\al_3}+x_{27}E_{\al_1}
\hspace{1.8cm}(4.38)\end{eqnarray*}
 by (3.4),
\begin{eqnarray*}\iota(\eta_{27})&=&-\iota([E_{-\al_1},\eta_{26}])=-[\iota(E_{-\al_1}),\iota(\eta_{26})]=
-[E_{-\al_1}|_{\cal A}+E_{-\al_1},\iota(\eta_{26})]\\
&=&P_{27}-x_8E'_{(1,2,2,3,2,1)} -x_{10}E'_{(1,1,2,3,2,1)}
-x_{27}\varrho(\al_4+\sum_{s=3}^5+\sum_{r=1}^6\al_r+\sum_{i=1}^7\al_i)\\
&&+x_{12}E'_{(1,1,2,2,2,1)}-x_{15}E'_{(1,1,1,2,2,1)}-x_{13}E'_{(1,1,2,2,1,1)}+x_{16}E'_{(1,1,1,2,1,1)}\\&
&+x_{17}E'_{(1,1,2,2,1)}-x_{19}E'_{(1,1,1,2,1)}-x_{18}E'_{(1,1,1,1,1,1)}-x_{20}E'_{(1,0,1,1,1,1)}+x_{21}E'_{(1,1,1,1,1)}
 \\& &-x_{22}E'_{(1,1,1,1)}+x_{23}E'_{(1,0,1,1,1)}
-x_{24}E'_{(1,0,1,1)}+x_{25}E'_{(1,0,1)}-x_{26}E_{-\al_1}
\hspace{1.2cm}(4.39)\end{eqnarray*}
 by (3.2). Note that the coefficient of $x_i$ in the above
 $\iota(\eta_i)$ is exactly $\varrho(\be)$ if $\eta_i=E_\be$.

Recall ${\cal A}=\mbb{C}[x_1,...,x_{27}]$. Let $M$ be a ${\cal
G}^{E_6}$-module and set
$$\widehat{M}={\cal A}\otimes_{\mbb{C}}M.\eqno(4.40)$$
We identify
$$f\otimes v=fv\qquad\for\;\;f\in{\cal A},\;v\in M.\eqno(4.41)$$
Recall the Lie algebra $\widehat{\cal K}$ defined via (4.1)-(4.8).
Fix $c\in\mbb{C}$. Then $\widehat M$ becomes a $\widehat{\cal
K}$-module with the action defined by
$$d(fw)=d(f)w,\;\;\kappa(fw)=cfw,\;\;(gu)(fw)=fg
u(w)\eqno(4.42)$$ for $d\in {\cal W}_{27},\;f,g\in{\cal A},\;w\in M$
and $u\in{\cal G}^{E_6}$.

Since the linear map $\iota: {\cal G}^{E_7}\rta \widehat{\cal K}$
defined in (4.10)-(4.12) is a  Lie algebra monomorphism,
$\widehat{M}$ becomes a ${\cal G}^{E_7}$-module with the action
defined by
$$\xi(w)=\iota(\xi)(w)\qquad\for\;\;\xi\in{\cal
G}^{E_7},\;w\in\widehat{M}.\eqno(4.43)$$ In fact, we have:\psp

{\bf Theorem 4.1}. {\it The map $M\mapsto \widehat M$ gives a
functor from ${\cal G}^{E_6}$-{\bf Mod} to ${\cal G}^{E_7}$-{\bf
Mod}.}\psp

We remark that  the module $\widehat M$ is not a generalized module
in general because it may not be equal to $U({\cal G})(M)=U({\cal
G}_-)(M)$. \psp

{\bf Proposition 4.2}. {\it If $M$ is an irreducible ${\cal
G}^{E_6}$-module, then $U({\cal G}_-)(M)$ is an irreducible ${\cal
G}^{E_7}$-module.}

{\it Proof.} Note that for any $i\in\ol{1,27}$, $f\in{\cal A}$ and
$v\in M$, (4.12), (4.42) and (4.43) imply
$$\xi_i(fv)=\ptl_{x_i}(f)v.\eqno(4.44)$$
Let $W$ be any nonzero ${\cal G}^{E_7}$-submodule. The above
expression shows that $W\bigcap M\neq\{0\}$. According to (4.42),
$W\bigcap M$ is a ${\cal G}^{E_6}$-submodule. By the irreducibility
of $M$, $M\subset W$. Thus $U({\cal G}_-)(M)\subset W$. So $U({\cal
G}_-)(M)=W$ is irreducible. $\qquad\Box$\psp

By the above proposition, the map $M\mapsto U({\cal G}_-)(M)$ is a
polynomial extension from  irreducible ${\cal G}^{E_6}$-modules to
irreducible ${\cal G}^{E_7}$-modules.

\section{Irreducibility}

\quad \quad In this section, we want to determine the irreducibility
of ${\cal G}^{E_6}$-modules $\widehat M$.

Note that $\widehat{M}$ can be viewed as a ${\cal G}^{E_6}$-module.
Indeed, (4.11) and (4.43) show
$$u(fv)=u(f)v+fu(v)\qquad\for\;\;u\in{\cal
G}^{E_6},\;f\in{\cal A},\;v\in M\eqno(5.1)$$ (cf. (2.50)-(2.87)). So
$\widehat{M}={\cal A}\otimes_{\mbb{C}}M$ is a tensor module of
${\cal G}^{E_6}$. Write
$$\eta^\al=\prod_{i=1}^{27}\eta_{i}^{\alpha_i},\;\;
 |\al|=\sum_{i=1}^{27}\al_i\;\; \for\;\;
\al=(\al_1,\al_2,...,\al_{27})\in\mbb{N}^{27}\eqno(5.2)$$ (cf.
(2.34)-(2.39)).
 Recall the Lie subalgebras ${\cal G}_\pm$ and ${\cal G}_0$ of ${\cal G}^{E_6}$ defined
 in (2.41).  For $k\in\mbb{N}$, we set
$${\cal A}_k=\mbox{Span}_{\mathbb{C}}\{x^\al\mid
\al \in\mbb{N}^{27};|\al|=k\},\;\; \widehat
M_{{\langle}k\rangle}={\cal A}_kM\eqno(5.3)$$ (cf. (2.47), (4.41))
and
 $$(U({\cal G}_-)(
M))_{{\langle}k\rangle}=\mbox{Span}_{\mathbb{C}}\{\eta^\al( M)\mid
\al \in\mbb{N}^{27}, \ |\al|=k\}. \eqno(5.4)$$
 Moreover,
$$(U({\cal G}_-)(
M))_{{\langle}0\rangle}=\widehat M_{{\langle}0\rangle}=
M.\eqno(5.5)$$ Furthermore,
 $$\widehat M=\bigoplus\limits_{k=0}^\infty\widehat M_{\langle
 k\rangle},\qquad
 U({\cal G}_-)(M)=\bigoplus\limits_{k=0}^\infty(U({\cal G}_-)(M))_{\langle k\rangle}.\eqno(5.6)$$

Next we define a linear transformation $\vf$ on  $\widehat M$
determined by
$$\vf(x^\al v)=\eta^\al(
v)\qquad\for\;\;\al\in\mbb{N}^{27},\;v\in M.\eqno(5.7)$$ Note that
${\cal A}_1=\sum_{i=1}^{27}\mbb{C}x_i$ forms the $27$-dimensional
${\cal G}_0$-module of highest weight $\lmd_6$. According to (2.10)
and (2.11), ${\cal G}_-$ forms a ${\cal G}_0$-module with respect to
the adjoint representation, and the linear map from ${\cal A}_1$ to
${\cal G}_0$ determined by $x_i\mapsto \eta_i$ for $i\in\ol{1,27}$
gives a ${\cal G}_0$-module isomorphism. Thus $\vf$ can also be
viewed as a ${\cal G}_0$-module homomorphism from $\widehat M$ to
$U({\cal G}_0)(M)$. Moreover,
$$\vf(\widehat M_{\la k\ra})=(U({\cal G}_-)(
M))_{{\langle}k\rangle}\qquad\for\;\;k\in\mbb{N}.\eqno(5.8)$$

Observe that the fundamental weights of ${\cal G}^{E_6}$ are:
$$\lmd_1=\frac{1}{3}(4\al_1+3\al_2+5\al_3+6\al_4+4\al_5+2\al_6),\eqno(5.9)$$
$$\lmd_2=\al_1+2\al_2+2\al_3+3\al_4+2\al_5+\al_6,\eqno(5.10)$$
$$\lmd_3=\frac{1}{3}(5\al_1+6\al_2+10\al_3+12\al_4+8\al_5+4\al_6),\eqno(5.11)$$
$$\lmd_4=2\al_1+3\al_2+4\al_3+6\al_4+4\al_5+2\al_6,\eqno(5.12)$$
$$\lmd_5=\frac{1}{3}(4\al_1+6\al_2+8\al_3+12\al_4+10\al_5+5\al_6),\eqno(5.13)$$
$$\lmd_6=\frac{1}{3}(2\al_1+3\al_2+4\al_3+6\al_4+5\al_5+4\al_6).\eqno(5.14)$$
Using the above expressions, we can view $\lmd_i\in{\cal H}_{E_6}$
by (2.16). Then Casimir element of ${\cal G}^{E_6}$ is
$$\omega=\sum_{i=1}^6\lmd_i\al_i-\sum_{\be\in\Phi_{E_6}}E_{\be}E_{-\be}\in
U({\cal G}^{E_6})\eqno(5.15)$$ due to (2.12) and (2.13).
 The algebra
$U({\cal G}^{E_6})$  can be imbedded into the tensor algebra
$U({\cal G}^{E_6})\otimes U({\cal G}^{E_6})$ by the associative
algebra homomorphism $\dlt: U({\cal G}^{E_6}) \rightarrow U({\cal
G}^{E_6})\otimes_\mbb{C} U({\cal G}^{E_6}))$ determined  by
$$\dlt(u)=u\otimes 1 +1 \otimes u \qquad \mbox{ for} \ u\in
{\cal G}^{E_6}.\eqno(5.16)$$ Set
$$\td\omega=\frac{1}{2}(\dlt(\omega)-\omega\otimes 1-1\otimes
\omega)\in U({\cal G}^{E_6})\otimes_\mbb{C} U({\cal
G}^{E_6}).\eqno(5.17)$$ Since $\sum_{i=1}^6\lmd_i\otimes \al_i$ is
symmetric with respect to $\{\al_1,...,\al_6\}$ by (5.9)-(5.14),
$$\td\omega=\sum_{i=1}^6\lmd_i\otimes\al_i-\sum_{\be\in\Phi_{E_6}}E_{\be}\otimes
E_{-\be}.\eqno(5.18)$$ Furthermore, $\td{\omega}$ acts on $\widehat
M$ as a ${\cal G}^{E_6}$-module homomorphism via
$$(u_1\otimes u_2)(fw)=u_1(f)u_2(w)\qquad\for\;\;u_1,u_2\in
{\cal G}^{E_6},\;f\in{\cal A},\;w\in M.\eqno(5.19)$$

{\bf Lemma 5.1}. {\it We have $\vf|_{\widehat M_{\la
1\ra}}=(\td\omega-c)|_{\widehat M_{\la 1\ra}}$.}

{\it Proof.} By (2.87) with Table 1 and (5.9)-(5.14)
$$\lmd_r|_{\cal
A}=\sum_{i=1}^{27}c_{i,r}x_i\ptl_{x_i}\qquad\for\;\;r\in\ol{1,6}\eqno(5.20)$$
with $c_{i,r}\in\mbb{C}$, for instance,
$$c_{1,1}=\frac{2}{3},\;c_{1,2}=1,\;c_{1,3}=\frac{4}{3},\;c_{1,4}=2,\;c_{1,5}=
\frac{5}{3},\;c_{1,6}=\frac{4}{3}.\eqno(5.21)$$ On the other hand,
$$\al_7=\frac{1}{3}[\widehat\al-(2\al_1+3\al_2+4\al_3+6\al_4+5\al_5+4\al_6)]\eqno(5.22)$$
by (2.88). According to (4.10), (4.42), (5.19) and (5.21),
$$(\sum_{i=1}^6\lmd_i\otimes\al_i)(x_1w)=cx_1w-x_1\al_7(w)\qquad\for\;w\in M.\eqno(5.23)$$
Thus (2.50)-(2.86) yield
\begin{eqnarray*}\td\omega(x_1w)&=&[cx_1+x_2E_{\al_6}+x_3E_{(0,0,0,0,1,1)}
+x_4E_{(0,0,0,1,1,1)}+x_5E_{(0,0,1,1,1,1)} +x_6E_{(0,1,0,1,1,1)}\\ &
&-x_1\al_7(w)+x_7E_{(0,1,1,1,1,1)}+x_8E_{(1,0,1,1,1,1)}+x_9E_{(0,1,1,2,1,1)}+x_{10}E_{(1,1,1,1,1,1)}
\\ & &+x_{11}E_{(0,1,1,2,2,1)}+x_{12}E_{(1,1,1,2,1,1)}+x_{13}E_{(1,1,1,2,2,1)}+x_{15}E_{(1,1,2,2,1,1)}\\
&
&+x_{16}E_{(1,1,2,2,2,1)}+x_{18}E_{(1,1,2,3,2,1)}+x_{20}E_{(1,2,2,3,2,1)}](w).\hspace{3.2cm}(5.24)
\end{eqnarray*}
Comparing (4.13) and (5.24), we get $\vf(x_1w)=(\td\omega-c)(x_1w)$.
According to (3.7), $-E_{-\al_6}(x_1)=x_2$. So for $v\in M$,
\begin{eqnarray*}\qquad& &\vf(x_2v)-\vf(x_1E_{-\al_6}(v))=
-(E_{-\al_6})(\vf(x_1v))\\
&=&-(E_{-\al_6})[(\td\omega-c)(x_1v)]=(\td\omega-c)[-E_{-\al_6}(x_1v)]
\\
&=&(\td\omega-c)[x_2v-x_1E_{-\al_6}(v)]\\
&=&(\td\omega-c)(x_2v)-(\td\omega-c)(x_1E_{-\al_6})(v))
\\ &=&(\td\omega-c)(x_2v)-\vf(x_1E_{-\al_6}(v)),\hspace{7.7cm}(5.25)\end{eqnarray*}
equivalently, $\vf(x_2v)=(\td\omega-c)(x_2v)$. By (3.2)-(3.7),
similar arguments as (5.25) and induction on $i$, we can prove
$$\vf(x_iv)=(\td\omega-c)(x_iv)\qquad\for\;\;i\in\ol{1,27},\;v\in M,\eqno(5.26)$$
that is, the lemma holds.$\qquad\Box$\psp

Recall that $\mbb{N}$ is the set of nonnegative integers and the set
of dominate integral weights is
$$\Lmd^+=\sum_{r=1}^6\mbb{N}\lmd_r.\eqno(5.27)$$
Denote
$$\rho=\frac{1}{2}\sum_{\be\in\Phi_{E_6}^+}\be.\eqno(5.28)$$
Then
$$\rho=\sum_{r=1}^6\lmd_r\eqno(5.29)$$
(e.g., cf. [H]). For any $\mu\in \Lmd^+$, we denote by $V(\mu)$ the
finite-dimensional irreducible ${\cal G}^{E_6}$-module with the
highest weight $\mu$ and have
$$\omega|_{V(\mu)}=(\mu+2\rho,\mu)\mbox{Id}_{V(\mu)}\eqno(5.30)$$
by (5.15). According to (2.87) and Table 1, the weight set of the
${\cal G}^{E_6}$-module ${\cal A}_1$ is
$$\Pi({\cal A}_1)=\{\sum_{r=1}^6a_{i,r}\lmd_r\mid
i\in\ol{1,27}\}.\eqno(5.31)$$
 Fixing
$\lmd\in\Lmd^+$, we define
$$\Upsilon(\lmd)=\{\lmd+\mu\mid\mu\in\Pi({\cal
A}_1),\;\lmd+\mu\in\Lmd^+\}.\eqno(5.32)$$ \pse

{\bf Lemma 5.2}. {\it We have}:
$${\cal A}_1\otimes V(\lmd)\cong \bigoplus_{\lmd'\in
\Upsilon(\lmd)}V(\lmd').\eqno(5.33)$$

{\it Proof}. Note that all the weight subspaces of ${\cal A}_1$ are
one-dimensional. Thus all the irreducible components of ${\cal
A}_1\otimes V(\lmd)$ are of multiplicity one. Since
$$\rho+\lmd+\mu\in\Lmd^+\qquad\for\;\;\mu\in \Pi({\cal
A}_1),\eqno(5.34)$$ the tensor theory of finite-dimensional
irreducible modules over a finite-dimensional simple Lie algebra
(e.g., cf. [H]) says that $V(\lmd')$ is a component of ${\cal
A}_1\otimes V(\lmd)$ if and only if $\lmd'\in
\Upsilon(\lmd).\qquad\Box$\psp

Recall
$$\mbox{the highest weight of}\;{\cal
A}_1=\lmd_6\eqno(5.35)$$ by Table 1. Thus the eigenvalues of
$\td{\omega}|_{\widehat{V(\lmd)}_{\la1\ra}}$ are
$$\{[(\lmd'+2\rho,\lmd')-(\lmd+2\rho,\lmd)-(\lmd_6+2\rho,\lmd_6)]/2\mid\lmd'\in
\Upsilon(\lmd)\}\eqno(5.36)$$ by (5.17) and (5.19). Define
$$\ell_\omega(\lmd)=\min\{[(\lmd'+2\rho,\lmd')-(\lmd+2\rho,\lmd)-(\lmd_6+2\rho,\lmd_6)]/2\mid\lmd'\in
\Upsilon(\lmd)\},\eqno(5.37)$$ which will be used to determine the
irreducibility of $\widehat{V(\lmd)}$. If
$\lmd'=\lmd+\lmd_6-\al\in\Upsilon(\lmd)$ with $\al\in\Phi_{D_6}^+$,
then
$$(\lmd'+2\rho,\lmd')-(\lmd+2\rho,\lmd)-(\lmd_6+2\rho,\lmd_6)
=2[(\lmd,\lmd_6)+1-(\rho+\lmd+\lmd_6,\al)].\eqno(5.38)$$

Recall the differential operators $P_1,...,P_{27}$ given in
(3.58)-(3.102). We also view the elements of ${\cal A}$ as the
multiplication operators on ${\cal A}$. Recall $\zeta_1$ in (3.1).
It turns out that we need the following lemma in order to determine
the irreducibility of $\widehat{V(\lmd)}$.\psp

{\bf Lemma 5.3}. {\it As operators on} ${\cal A}$:
\begin{eqnarray*}& &P_{14}x_1+P_1x_{14}+P_{11}x_2+P_2x_{11}+P_9x_3+P_3x_9+P_7x_4+P_4x_7
-P_6x_5-P_5x_6\\
&=&\zeta_1(D-8)+\chi\ptl_{x_{27}}.\hspace{10.3cm}(5.39)\end{eqnarray*}

{\it Proof}. According to  (3.58), (3.77)-(3.82), (3.84), (3.86) and
(3.89), we find that
\begin{eqnarray*}&
&P_{14}x_1+P_1x_{14}+P_{11}x_2+P_2x_{11}+P_9x_3+P_3x_9+P_7x_4+P_4x_7
-P_6x_5-P_5x_6\\
&=&-8\zeta_1+x_1P_{14}+x_{14}P_1+x_2P_{11}+x_{11}P_2+x_3P_9\\
& &+x_9P_3+x_4P_7+x_7P_4
-x_5P_6-x_6P_5\hspace{7.3cm}(5.40)\end{eqnarray*} and
\begin{eqnarray*}&
&x_1P_{14}+x_{14}P_1+x_2P_{11}+x_{11}P_2+x_3P_9+x_9P_3+x_4P_7+x_7P_4
-x_5P_6-x_6P_5
\\\!\!\!\!\!\!&=&\!\!\!x_1(x_{14}D-\zeta_1\ptl_{x_1}-\zeta_8\ptl_{x_8}-\zeta_{10}\ptl_{x_{10}}
+\zeta_{12}\ptl_{x_{12}}-\zeta_{13}\ptl_{x_{13}}-\zeta_{15}\ptl_{x_{15}}
+\zeta_{16}\ptl_{x_{16}}-\zeta_{18}\ptl_{x_{18}}\\ &
&-\zeta_{20}\ptl_{x_{20}}
+\zeta_{27}\ptl_{x_{27}})+x_{14}(x_1D-\zeta_1\ptl_{x_{14}}-\zeta_2\ptl_{x_{17}}-\zeta_3\ptl_{x_{19}}
+\zeta_4\ptl_{x_{21}}-\zeta_5\ptl_{x_{22}}
-\zeta_6\ptl_{x_{23}}\\&&+\zeta_7\ptl_{x_{24}}-\zeta_9\ptl_{x_{25}}-\zeta_{11}\ptl_{x_{26}}
-\zeta_{14}\ptl_{x_{27}})+x_2(x_{11}D-\zeta_1\ptl_{x_2}+\zeta_5\ptl_{x_8}+\zeta_7\ptl_{x_{10}}
-\zeta_9\ptl_{x_{12}}\\
& &+\zeta_{11}\ptl_{x_{15}}+\zeta_{13}\ptl_{x_{17}}
-\zeta_{16}\ptl_{x_{19}}+\zeta_{18}\ptl_{x_{21}}+\zeta_{20}\ptl_{x_{23}}
+\zeta_{26}\ptl_{x_{27}})+x_{11}(x_2D-\zeta_1\ptl_{x_{11}}\\
& &-\zeta_2\ptl_{x_{13}}-\zeta_3\ptl_{x_{16}}
+\zeta_4\ptl_{x_{18}}-\zeta_6\ptl_{x_{20}}-\zeta_8\ptl_{x_{22}}
+\zeta_{10}\ptl_{x_{24}}-\zeta_{12}\ptl_{x_{25}}-\zeta_{15}\ptl_{x_{26}}
+\zeta_{17}\ptl_{x_{27}})
\\ &&+x_3(x_9D-\zeta_1\ptl_{x_3}-\zeta_4\ptl_{x_8}-\zeta_6\ptl_{x_{10}}
+\zeta_9\ptl_{x_{13}}-\zeta_{12}\ptl_{x_{17}}-\zeta_{11}\ptl_{x_{16}}
+\zeta_{15}\ptl_{x_{19}}-\zeta_{18}\ptl_{x_{22}}\\ &
&-\zeta_{20}\ptl_{x_{24}}
+\zeta_{25}\ptl_{x_{27}})+x_9(x_3D-\zeta_1\ptl_{x_9}-\zeta_2\ptl_{x_{12}}-\zeta_3\ptl_{x_{15}}
+\zeta_5\ptl_{x_{18}}-\zeta_7\ptl_{x_{20}}-\zeta_8\ptl_{x_{21}}
\\ &&+\zeta_{10}\ptl_{x_{23}}-\zeta_{13}\ptl_{x_{25}}-\zeta_{16}\ptl_{x_{26}}
+\zeta_{19}\ptl_{x_{27}})+x_4(x_7D-\zeta_1\ptl_{x_4}+\zeta_3\ptl_{x_8}+\zeta_6\ptl_{x_{12}}
-\zeta_7\ptl_{x_{13}}\\ &
&+\zeta_{10}\ptl_{x_{17}}+\zeta_{11}\ptl_{x_{18}}
-\zeta_{15}\ptl_{x_{21}}+\zeta_{16}\ptl_{x_{22}}+\zeta_{20}\ptl_{x_{25}}
+\zeta_{24}\ptl_{x_{27}})+x_7(x_4D-\zeta_1\ptl_{x_7}\\
& &-\zeta_2\ptl_{x_{10}}-\zeta_4\ptl_{x_{15}}
+\zeta_5\ptl_{x_{16}}-\zeta_8\ptl_{x_{19}}-\zeta_9\ptl_{x_{20}}
+\zeta_{12}\ptl_{x_{23}}-\zeta_{13}\ptl_{x_{24}}-\zeta_{18}\ptl_{x_{26}}
+\zeta_{21}\ptl_{x_{27}})\\ &
&-x_5(x_6D+\zeta_1\ptl_{x_5}+\zeta_2\ptl_{x_8}-\zeta_6\ptl_{x_{15}}
+\zeta_7\ptl_{x_{16}}-\zeta_9\ptl_{x_{18}}-\zeta_{10}\ptl_{x_{19}}
+\zeta_{12}\ptl_{x_{21}}\\
& &-\zeta_{13}\ptl_{x_{22}}+\zeta_{20}\ptl_{x_{26}}
-\zeta_{23}\ptl_{x_{27}})-x_6(x_5D+\zeta_1\ptl_{x_6}-\zeta_3\ptl_{x_{10}}+\zeta_4\ptl_{x_{12}}
\\ &&-\zeta_5\ptl_{x_{13}}+\zeta_8\ptl_{x_{17}}+\zeta_{11}\ptl_{x_{20}}
-\zeta_{15}\ptl_{x_{23}}+\zeta_{16}\ptl_{x_{24}}-\zeta_{18}\ptl_{x_{25}}
-\zeta_{22}\ptl_{x_{27}})
\\\!\!\!\!\!\!&=&\!\!\!\zeta_1(2D-\sum_{i=1}^7x_i\ptl_{x_i}-x_9\ptl_{x_9}-x_{11}\ptl_{x_{11}}-x_{14}\ptl_{x_{14}})
-(x_1\zeta_8-x_2\zeta_5+x_3\zeta_4-x_4\zeta_3\\
&&+x_5\zeta_2)\ptl_{x_8}-(x_1\zeta_{10}-x_2\zeta_7+x_3\zeta_6+x_7\zeta_2-x_6\zeta_3)\ptl_{x_{10}}+(x_1\zeta_{12}-x_2\zeta_9-x_9\zeta_2
\\&&+x_4\zeta_6-x_6\zeta_4)\ptl_{x_{12}}-(x_1\zeta_{13}+x_{11}\zeta_2-x_3\zeta_9+x_4\zeta_7-x_6\zeta_5)\ptl_{x_{13}}
-(x_1\zeta_{15}-x_2\zeta_{11}\\&&+x_9\zeta_3+x_7\zeta_4-x_5\zeta_6)\ptl_{x_{15}}
+(x_1\zeta_{16}-x_{11}\zeta_3-x_3\zeta_{11}+x_7\zeta_5-x_5\zeta_7)\ptl_{x_{16}}
-(x_{14}\zeta_2\\&&-x_2\zeta_{13}+x_3\zeta_{12}-x_4\zeta_{10}+x_6\zeta_8)\ptl_{x_{17}}
-(x_1\zeta_{18}-x_{11}\zeta_4-x_9\zeta_5-x_4\zeta_{11}-x_5\zeta_9)\ptl_{x_{18}}\\&&
-(x_{14}\zeta_3+x_2\zeta_{16}-x_3\zeta_{15}+x_7\zeta_8-x_5\zeta_{10})\ptl_{x_{19}}
-(x_1\zeta_{20}+x_{11}\zeta_6+x_9\zeta_7+x_7\zeta_9\\&&+x_6\zeta_{11})\ptl_{x_{20}}
+(x_{14}\zeta_4+x_2\zeta_{18}-x_9\zeta_8-x_4\zeta_{15}-x_5\zeta_{12})\ptl_{x_{21}}
-(x_{14}\zeta_5+x_{11}\zeta_8\\&
&+x_3\zeta_{18}-x_4\zeta_{16}-x_5\zeta_{13})\ptl_{x_{22}}
-(x_{14}\zeta_6-x_2\zeta_{20}-x_9\zeta_{10}-x_7\zeta_{12}-x_6\zeta_{15})\ptl_{x_{23}}
\\& &+(x_{14}\zeta_7+x_{11}\zeta_{10}-x_3\zeta_{20}-x_7\zeta_{13}-x_6\zeta_{16})\ptl_{x_{24}}
-(x_{14}\zeta_9+x_{11}\zeta_{12}+x_9\zeta_{13}\\&
&-x_4\zeta_{20}-x_6\zeta_{18})\ptl_{x_{25}}
-(x_{14}\zeta_{11}+x_{11}\zeta_{15}+x_9\zeta_{16}+x_7\zeta_{18}+x_5\zeta_{20})\ptl_{x_{26}}
+(x_1\zeta_{27}\\&&-x_{14}\zeta_{14}+x_2\zeta_{26}+x_{11}\zeta_{17}+x_3\zeta_{25}
+x_9\zeta_{19}+x_4\zeta_{24}+x_7\zeta_{21}+x_5\zeta_{23}+x_6\zeta_{22})\ptl_{x_{27}}
\\\!\!\!\!\!\!&=&\!\!\!\zeta_1D+\chi\ptl_{x_{27}}\hspace{11.9cm}(5.41)
\end{eqnarray*}
because
$$x_1\zeta_8-x_2\zeta_5+x_3\zeta_4-x_4\zeta_3+x_5\zeta_2=\zeta_1x_8\eqno(5.42)$$
by (3.1), (3.8)-(3.11) and (3.14),
$$x_1\zeta_{10}-x_2\zeta_7+x_3\zeta_6+x_7\zeta_2-x_6\zeta_3=\zeta_1x_{10}\eqno(5.43)$$
by (3.1), (3.8), (3.9), (3.12), (3.13) and (3.16),
$$x_1\zeta_{12}-x_2\zeta_9-x_9\zeta_2+x_4\zeta_6-x_6\zeta_4=-\zeta_1x_{12}\eqno(5.44)$$
by (3.1), (3.8), (3.10), (3.12), (3.15) and (3.18),
$$x_1\zeta_{13}+x_{11}\zeta_2-x_3\zeta_9+x_4\zeta_7-x_6\zeta_5=\zeta_1x_{13}\eqno(5.45)$$
by (3.1), (3.8), (3.11), (3.13), (3.15) and (3.19),
$$x_1\zeta_{15}-x_2\zeta_{11}+x_9\zeta_3+x_7\zeta_4-x_5\zeta_6=\zeta_1x_{15}\eqno(5.46)$$
by (3.1), (3.9), (3.10), (3.12), (3.17) and (3.21),
$$x_1\zeta_{16}-x_{11}\zeta_3-x_3\zeta_{11}+x_7\zeta_5-x_5\zeta_7=-\zeta_1x_{16}\eqno(5.47)$$
by (3.1), (3.9), (3.11), (3.13), (3.17) and (3.22),
$$x_{14}\zeta_2-x_2\zeta_{13}+x_3\zeta_{12}-x_4\zeta_{10}+x_6\zeta_8=\zeta_1x_{17}\eqno(5.48)$$
by (3.1), (3.8), (3.14), (3.16), (3.18) and (3.19),
$$x_1\zeta_{18}-x_{11}\zeta_4-x_9\zeta_5-x_4\zeta_{11}-x_5\zeta_9=\zeta_1x_{18}\eqno(5.49)$$
by (3.1), (3.10), (3.11), (3.15), (3.17) and (3.24),
$$x_{14}\zeta_3+x_2\zeta_{16}-x_3\zeta_{15}+x_7\zeta_8-x_5\zeta_{10}=\zeta_1x_{19}\eqno(5.50)$$
by (3.1), (3.9), (3.14), (3.16), (3.21) and (3.22),
$$x_1\zeta_{20}+x_{11}\zeta_6+x_9\zeta_7+x_7\zeta_9+x_6\zeta_{11}=\zeta_1x_{20}\eqno(5.51)$$
by (3.1), (3.12), (3.13), (3.15), (3.17) and (3.26),
$$x_{14}\zeta_4+x_2\zeta_{18}-x_9\zeta_8-x_4\zeta_{15}-x_5\zeta_{12}=-\zeta_1x_{21}\eqno(5.52)$$
by (3.1), (3.10), (3.14), (3.18), (3.21) and (3.24),
$$x_{14}\zeta_5+x_{11}\zeta_8+x_3\zeta_{18}-x_4\zeta_{16}-x_5\zeta_{13}=\zeta_1x_{22}\eqno(5.53)$$
by (3.1), (3.11), (3.14), (3.19), (3.22) and (3.24),
$$x_{14}\zeta_6-x_2\zeta_{20}-x_9\zeta_{10}-x_7\zeta_{12}-x_6\zeta_{15}=\zeta_1x_{23}\eqno(5.54)$$
by (3.1), (3.12), (3.16), (3.18), (3.21) and (3.26),
$$x_{14}\zeta_7+x_{11}\zeta_{10}-x_3\zeta_{20}-x_7\zeta_{13}-x_6\zeta_{16}=-\zeta_1x_{24}\eqno(5.55)$$
by (3.1), (3.13), (3.16), (3.19), (3.22) and (3.26),
$$x_{14}\zeta_9+x_{11}\zeta_{12}+x_9\zeta_{13}-x_4\zeta_{20}-x_6\zeta_{18}=\zeta_1x_{25}\eqno(5.56)$$
by (3.1), (3.15), (3.18), (3.19), (3.24) and (3.26),
$$x_{14}\zeta_{11}+x_{11}\zeta_{15}+x_9\zeta_{16}+x_7\zeta_{18}+x_5\zeta_{20}=\zeta_1x_{26}\eqno(5.57)$$
by (3.1), (3.17), (3.21), (3.22), (3.24) and (3.26),
\begin{eqnarray*}\qquad
&&x_1\zeta_{27}-x_{14}\zeta_{14}+x_2\zeta_{26}+x_{11}\zeta_{17}+x_3\zeta_{25}
+x_9\zeta_{19}+x_4\zeta_{24}\\
& &+x_7\zeta_{21}+x_5\zeta_{23}+x_6\zeta_{22}
=-\zeta_1x_{27}+\chi\hspace{6.8cm}(5.58)\end{eqnarray*} by (3.1),
(3.20), (3.23), (3.25), (3.27)-(3.33) and (3.137).$\qquad\Box$\psp

{\bf Lemma 5.4}. {\it As operators on ${\cal A}$,}
$$\sum_{14\neq
i\in\ol{1,26}}P_i\zeta_{28-i}-P_{14}\zeta_{14}-P_{27}\zeta_1=\chi(24-5D).\eqno(5.59)$$

{\it Proof}. We calculate it by (3.1), (3.8)-(3.33), (3.58) and
(3.77)-(3.102). In particular,
$$2(\zeta_1\zeta_{14}-\zeta_2\zeta_{11}-\zeta_3\zeta_9+\zeta_4\zeta_7-\zeta_5\zeta_6)=-2\chi x_1,\eqno(5.60)$$
which is the coefficient of $\ptl_{x_1}$ in addition to the term
containing $D$. According to (2.138), the operator on the left hand
side of (5.59) is a ${\cal G}^{E_6}$-invariant differential
operator. By symmetry,
$$\mbox{the coefficient of}\;\ptl_{x_i}=-2\chi x_i\;\;\for\;\;i\in\ol{1,27}.\qquad\Box\eqno(5.61)$$
\pse

We define the multiplication
$$f(gv)=(fg)v\qquad\for\;\;f,g\in{\cal A},\;v\in M.\eqno(5.62)$$
 Furthermore, we have:\psp

{\bf Lemma 5.5}. {\it As operators on $\widehat M$,}
$$[(\sum_{14\neq
i\in\ol{1,26}}\iota(\eta_i)\zeta_{28-i}-\iota(\eta_{14})\zeta_{14}-\iota(\eta_{27})\zeta_1]|_{\widehat
M}=\chi(24-5D+3c).\eqno(5.63)$$

{\it Proof}. By (4.13)-(4.39), the coefficient of $E_{\al_6}$ is
$$x_2\zeta_{27}-x_{14}\zeta_{17}-x_{17}\zeta_{15}-x_{19}\zeta_{12}-x_{21}\zeta_{10}-x_{23}\zeta_8.\eqno(5.64)$$
Moreover, we use (3.14), (3.16), (3.18), (3.21), (3.23) and (3.33)
to find that (5.64) is equal to zero. Since  the left hand side of
(5.63) is a ${\cal G}^{E_6}$-invariant differential operator, it is
invariant under the action of the $E_6$ Weyl group. The transitivity
the Weyl group on $\Phi_{E_6}$ yields that
$$\mbox{the coefficient of}\;E_\be=0\;\;\mbox{for
any}\;\be\in\Phi_{E_6}.\eqno(5.65)$$ By (4.10), (4.42), (5.22) and
Lemma 5.4,
$$\sum_{14\neq
i\in\ol{1,26}}\iota(\eta_i)\zeta_{28-i}-\iota(\eta_{14})\zeta_{14}-\iota(\eta_{27})\zeta_1
=\chi(24-5D+3c)+\sum_{r=1}^6f_r\al_r\eqno(5.66)$$ as the elements of
$\widehat {\cal K}$ acting on $\widehat M$ (cf. (4.7)). The ${\cal
G}^{E_6}$-invariancy implies
$$[E_\be|_{\widehat M},\sum_{r=1}^6f_r\al_r]=0\;\;\mbox{for
any}\;\be\in\Phi_{E_6}.\eqno(5.68)$$ Thus $\sum_{r=1}^6f_r\al_r=0$,
that is, (3.63) holds. $\qquad\Box$\psp

Next we calculate
\begin{eqnarray*}T_1&=&\iota(\eta_{14})x_1+\iota(\eta_1)x_{14}+\iota(\eta_{11})x_2+\iota(\eta_2)x_{11}+\iota(\eta_9)x_3
\\ &&+\iota(\eta_3)x_9+\iota(\eta_7)x_4+\iota(\eta_4)x_7
-\iota(\eta_6)x_5-\iota(\eta_5)x_6\hspace{5cm}\end{eqnarray*}\begin{eqnarray*}
\quad&=&\zeta_1(D-8)+\chi\ptl_{x_{27}}+\zeta_2E_{\al_1}-\zeta_3E_{(1,0,1)}+\zeta_4E_{(1,0,1,1)}-\zeta_5E_{(1,0,1,1,1)}+\zeta_6E_{(1,1,1,1)}
\\&
&-\zeta_7E_{(1,1,1,1,1)}+\zeta_8E_{(1,0,1,1,1,1)}+\zeta_9E_{(1,1,1,2,1)}+\zeta_{10}E_{(1,1,1,1,1,1)}-\zeta_{11}E_{(1,1,2,2,1)}
\\ & &-\zeta_{12}E_{(1,1,1,2,1,1)}+\zeta_{13}E_{(1,1,1,2,2,1)}+\zeta_{15}E_{(1,1,2,2,1,1)}-\zeta_{16}E_{(1,1,2,2,2,1)}+\zeta_{18}E_{(1,1,2,3,2,1)}
\\& &+\zeta_{20}E_{(1,2,2,3,2,1)}-2\zeta_1
c+\frac{\zeta_1}{3}(4\al_1+3\al_2+5\al_3+6\al_4+4\al_5+2\al_6)
 \hspace{1.8cm}(5.69)\end{eqnarray*}
by (3.1), (3.8)-(3.19), (3.21), (3.22), (3.24), (3.26),
(4.13)-(4.19), (4.21), (4.23), (4.26), (5.20) and Lemma 5.3.

We define a ${\cal G}^{E_6}$-module
 structure on the space $\mbox{End}\:\widehat M$ of linear
 transformations on $\widehat M$ by
 $$\iota(u)(T)=[\iota(u),T]=\iota(u)T-T\iota(u)\qquad\for\;\;u\in
{\cal G}^{E_6},\;T\in \mbox{End}\:\widehat M\eqno(5.70)$$ (cf.
(5.1)). It can be verified that $T_1$ is a ${\cal G}^{E_6}$-singular
vector with weight $\lmd_1$ in $\mbox{End}\:\widehat M$. So it
generates the 27-dimensional  module of highest weight $\lmd_1$.
  We set
\begin{eqnarray*}T_2&=&-[\iota(E_{-\al_1}),T_1]=
\zeta_2(D-8)-\chi\ptl_{x_{26}}-\zeta_1E_{-\al_1}+\zeta_3E_{\al_3}-\zeta_4E_{(0,0,1,1)}\\&
&+\zeta_5E_{(0,0,1,1,1)}-\zeta_6E_{(0,1,1,1)}
+\zeta_7E_{(0,1,1,1,1)}-\zeta_8E_{(0,0,1,1,1,1)}-\zeta_9E_{(0,1,1,2,1)}\\
& &-\zeta_{10}E_{(0,1,1,1,1,1)}-\zeta_{14}E_{(1,1,2,2,1)}
+\zeta_{12}E_{(0,1,1,2,1,1)}-\zeta_{13}E_{(0,1,1,2,2,1)}\\&
&-\zeta_{17}E_{(1,1,2,2,1,1)}+\zeta_{19}E_{(1,1,2,2,2,1)}-\zeta_{21}E_{(1,1,2,3,2,1)}
-\zeta_{23}E_{(1,2,2,3,2,1)}\\& &-2\zeta_2
c+\frac{\zeta_2}{3}(\al_1+3\al_2+5\al_3+6\al_4+4\al_5+2\al_6)
 \hspace{5cm}(5.71)\end{eqnarray*}
by (2.10), (2.11), (3.2) and (3.41),
\begin{eqnarray*}T_3&=&-[\iota(E_{-\al_3}),T_2]=
\zeta_3(D-8)-\chi\ptl_{x_{25}}+\zeta_1E'_{(1,0,1)}-\zeta_2E_{-\al_3}+\zeta_4E_{\al_4}\\&
&-\zeta_5E_{(0,0,0,1,1)}+\zeta_6E_{(0,1,0,1)}
-\zeta_7E_{(0,1,0,1,1)}+\zeta_8E_{(0,0,0,1,1,1)}+\zeta_{11}E_{(0,1,1,2,1)}\\
& &+\zeta_{10}E_{(0,1,0,1,1,1)}+\zeta_{14}E_{(1,1,1,2,1)}
-\zeta_{15}E_{(0,1,1,2,1,1)}+\zeta_{16}E_{(0,1,1,2,2,1)}\\&
&+\zeta_{17}E_{(1,1,1,2,1,1)}-\zeta_{19}E_{(1,1,1,2,2,1)}+\zeta_{22}E_{(1,1,2,3,2,1)}
+\zeta_{24}E_{(1,2,2,3,2,1)}\\& &-2\zeta_3
c+\frac{\zeta_3}{3}(\al_1+3\al_2+2\al_3+6\al_4+4\al_5+2\al_6)
 \hspace{5cm}(5.72)\end{eqnarray*}
by (2.10), (2.11), (3.4) and (3.43),
\begin{eqnarray*}T_4&=&-[\iota(E_{-\al_4}),T_3]=
\zeta_4(D-8)-\chi\ptl_{x_{24}}-\zeta_1E'_{(1,0,1,1)}+\zeta_2E'_{(0,0,1,1)}-\zeta_3E_{-\al_4}\\&
&+\zeta_5E_{\al_5}+\zeta_6E_{\al_2}
-\zeta_9E_{(0,1,0,1,1)}-\zeta_8E_{(0,0,0,0,1,1)}+\zeta_{11}E_{(0,1,1,1,1)}\\
& &+\zeta_{12}E_{(0,1,0,1,1,1)}+\zeta_{14}E_{(1,1,1,1,1)}
-\zeta_{15}E_{(0,1,1,1,1,1)}+\zeta_{18}E_{(0,1,1,2,2,1)}\\&
&+\zeta_{17}E_{(1,1,1,1,1,1)}-\zeta_{21}E_{(1,1,1,2,2,1)}+\zeta_{22}E_{(1,1,2,2,2,1)}
-\zeta_{25}E_{(1,2,2,3,2,1)}\\& &-2\zeta_4
c+\frac{\zeta_4}{3}(\al_1+3\al_2+2\al_3+3\al_4+4\al_5+2\al_6)
 \hspace{5cm}(5.73)\end{eqnarray*}
by (2.10), (2.11), (3.5) and (3.44),
\begin{eqnarray*}T_5&=&-[\iota(E_{-\al_5}),T_4]=
\zeta_5(D-8)-\chi\ptl_{x_{23}}+\zeta_1E'_{(1,0,1,1,1)}-\zeta_2E'_{(0,0,1,1,1)}\\&
&+\zeta_3E'_{(0,0,0,1,1)}-\zeta_4E_{-\al_5}+\zeta_7E_{\al_2}
-\zeta_9E_{(0,1,0,1)}+\zeta_8E_{\al_6}+\zeta_{11}E_{(0,1,1,1)}\\
& &+\zeta_{13}E_{(0,1,0,1,1,1)}+\zeta_{14}E_{(1,1,1,1)}
-\zeta_{16}E_{(0,1,1,1,1,1)}+\zeta_{18}E_{(0,1,1,2,1,1)}\\&
&+\zeta_{19}E_{(1,1,1,1,1,1)}-\zeta_{21}E_{(1,1,1,2,1,1)}+\zeta_{22}E_{(1,1,2,2,1,1)}
+\zeta_{26}E_{(1,2,2,3,2,1)}\\& &-2\zeta_5
c+\frac{\zeta_5}{3}(\al_1+3\al_2+2\al_3+3\al_4+\al_5+2\al_6)
 \hspace{5.1cm}(5.74)\end{eqnarray*}
by (2.10), (2.11), (3.6) and (3.45),
\begin{eqnarray*}T_6&=&-[\iota(E_{-\al_2}),T_4]=
\zeta_6(D-8)-\chi\ptl_{x_{22}}-\zeta_1E'_{(1,1,1,1)}+\zeta_2E'_{(0,1,1,1)}-\zeta_3E'_{(0,1,0,1)}\\&
&+\zeta_7E_{\al_5}-\zeta_4E_{-\al_2}
+\zeta_9E_{(0,0,0,1,1)}-\zeta_{10}E_{(0,0,0,0,1,1)}-\zeta_{11}E_{(0,0,1,1,1)}\\
& &-\zeta_{12}E_{(0,0,0,1,1,1)}-\zeta_{14}E_{(1,0,1,1,1)}
+\zeta_{15}E_{(0,0,1,1,1,1)}-\zeta_{20}E_{(0,1,1,2,2,1)}\\&
&-\zeta_{17}E_{(1,0,1,1,1,1)}+\zeta_{23}E_{(1,1,1,2,2,1)}-\zeta_{24}E_{(1,1,2,2,2,1)}
-\zeta_{25}E_{(1,1,2,3,2,1)}\\& &-2\zeta_6
c+\frac{\zeta_6}{3}(\al_1+2\al_3+3\al_4+4\al_5+2\al_6)
 \hspace{6.1cm}(5.75)\end{eqnarray*}
by (2.10), (2.11), (3.3) and (3.42),
\begin{eqnarray*}T_7&=&-[\iota(E_{-\al_5}),T_6]=
\zeta_7(D-8)-\chi\ptl_{x_{21}}+\zeta_1E'_{(1,1,1,1,1)}-\zeta_2E'_{(0,1,1,1,1)}\\&
&+\zeta_3E'_{(0,1,0,1,1)}-\zeta_6E_{-\al_5}-\zeta_5E_{-\al_2}
+\zeta_9E_{\al_4}+\zeta_{10}E_{\al_6}-\zeta_{11}E_{(0,0,1,1)}\\
& &-\zeta_{13}E_{(0,0,0,1,1,1)}-\zeta_{14}E_{(1,0,1,1)}
+\zeta_{16}E_{(0,0,1,1,1,1)}-\zeta_{20}E_{(0,1,1,2,1,1)}\\&
&-\zeta_{19}E_{(1,0,1,1,1,1)}+\zeta_{23}E_{(1,1,1,2,1,1)}-\zeta_{24}E_{(1,1,2,2,1,1)}
+\zeta_{26}E_{(1,1,2,3,2,1)}\\& &-2\zeta_7
c+\frac{\zeta_7}{3}(\al_1+2\al_3+3\al_4+\al_5+2\al_6)
 \hspace{6.3cm}(5.76)\end{eqnarray*}
by (2.10), (2.11), (3.6) and (3.45),
\begin{eqnarray*}T_8&=&-[\iota(E_{-\al_6}),T_5]=
\zeta_8(D-8)-\chi\ptl_{x_{20}}-\zeta_1E'_{(1,0,1,1,1,1)}+\zeta_2E'_{(0,0,1,1,1,1)}\\&
&-\zeta_3E'_{(0,0,0,1,1,1)}+\zeta_4E'_{(0,0,0,0,1,1)}-\zeta_5E_{-\al_6}+\zeta_{10}E_{\al_2}
-\zeta_{12}E_{(0,1,0,1)}\\
&&+\zeta_{13}E_{(0,1,0,1,1)}+\zeta_{15}E_{(0,1,1,1)}-\zeta_{17}E_{(1,1,1,1)}
-\zeta_{16}E_{(0,1,1,1,1)}+\zeta_{18}E_{(0,1,1,2,1)}\\&
&+\zeta_{19}E_{(1,1,1,1,1)}-\zeta_{21}E_{(1,1,1,2,1)}+\zeta_{22}E_{(1,1,2,2,1)}
-\zeta_{27}E_{(1,2,2,3,2,1)}\\& &-2\zeta_8
c+\frac{\zeta_8}{3}(\al_1+3\al_2+2\al_3+3\al_4+\al_5-\al_6)
 \hspace{5.3cm}(5.77)\end{eqnarray*}
by (2.10), (2.11), (3.7) and (3.46),
\begin{eqnarray*}T_9&=&-[\iota(E_{-\al_4}),T_7]=
\zeta_9(D-8)-\chi\ptl_{x_{19}}-\zeta_1E'_{(1,1,1,2,1)}+\zeta_2E'_{(0,1,1,2,1)}\\&
&+\zeta_4E'_{(0,1,0,1,1)}-\zeta_6E'_{(0,0,0,1,1)}+\zeta_5E'_{(0,1,0,1)}
-\zeta_7E_{-\al_4}-\zeta_{11}E_{\al_3}+\zeta_{12}E_{\al_6}\\
& &+\zeta_{13}E_{(0,0,0,0,1,1)}-\zeta_{14}E_{(1,0,1)}
+\zeta_{18}E_{(0,0,1,1,1,1)}-\zeta_{20}E_{(0,1,1,1,1,1)}\\&
&-\zeta_{21}E_{(1,0,1,1,1,1)}+\zeta_{23}E_{(1,1,1,1,1,1)}+\zeta_{25}E_{(1,1,2,2,1,1)}
+\zeta_{26}E_{(1,1,2,2,2,1)}\\& &-2\zeta_9
c+\frac{\zeta_9}{3}(\al_1+2\al_3+\al_5+2\al_6)
 \hspace{7.4cm}(5.78)\end{eqnarray*}
by (2.10), (2.11), (3.5) and (3.44),
\begin{eqnarray*}T_{10}&=&-[\iota(E_{-\al_2}),T_8]=
\zeta_{10}(D-8)-\chi\ptl_{x_{18}}-\zeta_1E'_{(1,1,1,1,1,1)}+\zeta_2E'_{(0,1,1,1,1,1)}\\&
&-\zeta_3E'_{(0,1,0,1,1,1)}+\zeta_6E'_{(0,0,0,0,1,1)}-\zeta_7E_{-\al_6}-\zeta_8E_{-\al_2}
+\zeta_{12}E_{\al_4}\\
&&-\zeta_{13}E_{(0,0,0,1,1)}-\zeta_{15}E_{(0,0,1,1)}+\zeta_{17}E_{(1,0,1,1)}
+\zeta_{16}E_{(0,0,1,1,1)}-\zeta_{20}E_{(0,1,1,2,1)}\\&
&-\zeta_{19}E_{(1,0,1,1,1)}+\zeta_{23}E_{(1,1,1,2,1)}-\zeta_{24}E_{(1,1,2,2,1)}
-\zeta_{27}E_{(1,1,2,3,2,1)}\\& &-2\zeta_{10}
c+\frac{\zeta_{10}}{3}(\al_1+2\al_3+3\al_4+\al_5-\al_6)
 \hspace{6cm}(5.79)\end{eqnarray*}
by (2.10), (2.11), (3.3) and (3.42),
\begin{eqnarray*}T_{11}&=&[\iota(E_{-\al_3}),T_9]=
\zeta_{11}(D-8)-\chi\ptl_{x_{17}}+\zeta_1E'_{(1,1,2,2,1)}-\zeta_3E'_{(0,1,1,2,1)}\\&
&-\zeta_4E'_{(0,1,1,1,1)}+\zeta_6E'_{(0,0,1,1,1)}-\zeta_5E'_{(0,1,1,1)}
+\zeta_7E'_{(0,0,1,1)}+\zeta_9E_{-\al_3}\hspace{5cm}\end{eqnarray*}\begin{eqnarray*}
\qquad
&&+\zeta_{14}E_{\al_1}+\zeta_{15}E_{\al_6}+\zeta_{16}E_{(0,0,0,0,1,1)}
+\zeta_{18}E_{(0,0,0,1,1,1)}-\zeta_{20}E_{(0,1,0,1,1,1)}\\&
&-\zeta_{22}E_{(1,0,1,1,1,1)}+\zeta_{24}E_{(1,1,1,1,1,1)}+\zeta_{25}E_{(1,1,1,2,1,1)}
+\zeta_{26}E_{(1,1,1,2,2,1)}\\& &-2\zeta_{11}
c+\frac{\zeta_{11}}{3}(\al_1-\al_3+\al_5+2\al_6)
 \hspace{7.3cm}(5.80)\end{eqnarray*}
by (2.10), (2.11), (3.4) and (3.43),
\begin{eqnarray*}T_{12}&=&-[\iota(E_{-\al_4}),T_{10}]=
\zeta_{12}(D-8)-\chi\ptl_{x_{16}}+\zeta_1E'_{(1,1,1,2,1,1)}-\zeta_2E'_{(0,1,1,2,1,1)}\\&
&-\zeta_4E'_{(0,1,0,1,1,1)}+\zeta_6E'_{(0,0,0,1,1,1)}+\zeta_8E'_{(0,1,0,1)}-\zeta_9E_{-\al_6}
-\zeta_{10}E_{-\al_4}\\
&&+\zeta_{13}E_{\al_5}-\zeta_{15}E_{\al_3}+\zeta_{17}E_{(1,0,1)}
+\zeta_{18}E_{(0,0,1,1,1)}-\zeta_{20}E_{(0,1,1,1,1)}\\&
&-\zeta_{21}E_{(1,0,1,1,1)}+\zeta_{23}E_{(1,1,1,1,1)}+\zeta_{25}E_{(1,1,2,2,1)}
-\zeta_{27}E_{(1,1,2,2,2,1)}\\& &-2\zeta_{12}
c+\frac{\zeta_{12}}{3}(\al_1+2\al_3+\al_5-\al_6)
 \hspace{7cm}(5.81)\end{eqnarray*}
by (2.10), (2.11), (3.5) and (3.44),
\begin{eqnarray*}T_{13}&=&-[\iota(E_{-\al_5}),T_{12}]=
\zeta_{13}(D-8)-\chi\ptl_{x_{15}}-\zeta_1E'_{(1,1,1,2,2,1)}-\zeta_2E'_{(0,1,1,2,2,1)}\\&
&-\zeta_5E'_{(0,1,0,1,1,1)}+\zeta_7E'_{(0,0,0,1,1,1)}-\zeta_8E'_{(0,1,0,1,1)}-\zeta_9E'_{(0,0,0,0,1,1)}
+\zeta_{10}E'_{(0,0,0,1,1)}\\
&&-\zeta_{12}E_{-\al_5}-\zeta_{16}E_{\al_3}+\zeta_{19}E_{(1,0,1)}
+\zeta_{18}E_{(0,0,1,1)}-\zeta_{20}E_{(0,1,1,1)}\\&
&-\zeta_{21}E_{(1,0,1,1)}+\zeta_{23}E_{(1,1,1,1)}-\zeta_{26}E_{(1,1,2,2,1)}
-\zeta_{27}E_{(1,1,2,2,1,1)}\\& &-2\zeta_{13}
c+\frac{\zeta_{13}}{3}(\al_1+2\al_3-2\al_5-\al_6)
 \hspace{6.8cm}(5.82)\end{eqnarray*}
by (2.10), (2.11), (3.6) and (3.45),
\begin{eqnarray*}T_{14}&=&-[\iota(E_{-\al_1}),T_{11}]=
\zeta_{14}(D-8)+\chi\ptl_{x_{14}}+\zeta_2E'_{(1,1,2,2,1)}-\zeta_3E'_{(1,1,1,2,1)}\\&
&-\zeta_4E'_{(1,1,1,1,1)}+\zeta_6E'_{(1,0,1,1,1)}-\zeta_5E'_{(1,1,1,1)}
+\zeta_7E'_{(1,0,1,1)}+\zeta_9E'_{(1,0,1)}\\
&&-\zeta_{11}E_{-\al_1}-\zeta_{17}E_{\al_6}-\zeta_{19}E_{(0,0,0,0,1,1)}
-\zeta_{21}E_{(0,0,0,1,1,1)}+\zeta_{23}E_{(0,1,0,1,1,1)}\\&
&+\zeta_{22}E_{(0,0,1,1,1,1)}-\zeta_{24}E_{(0,1,1,1,1,1)}-\zeta_{25}E_{(0,1,1,2,1,1)}
-\zeta_{26}E_{(0,1,1,2,2,1)}\\& &-2\zeta_{14}
c+\frac{\zeta_{14}}{3}(-2\al_1-\al_3+\al_5+2\al_6)
 \hspace{6.7cm}(5.83)\end{eqnarray*}
by (2.10), (2.11), (3.2) and (3.41),
\begin{eqnarray*}T_{15}&=&[\iota(E_{-\al_3}),T_{12}]=
\zeta_{15}(D-8)-\chi\ptl_{x_{13}}-\zeta_1E'_{(1,1,2,2,1,1)}+\zeta_3E'_{(0,1,1,2,1,1)}\\&
&+\zeta_4E'_{(0,1,1,1,1,1)}-\zeta_6E'_{(0,0,1,1,1,1)}-\zeta_8E'_{(0,1,1,1)}
+\zeta_{10}E'_{(0,0,1,1)}-\zeta_{11}E_{-\al_6}\\
&&+\zeta_{12}E_{-\al_3}+\zeta_{16}E_{\al_5}-\zeta_{17}E_{\al_1}
+\zeta_{18}E_{(0,0,0,1,1)}-\zeta_{20}E_{(0,1,0,1,1)}\\&
&-\zeta_{22}E_{(1,0,1,1,1)}+\zeta_{24}E_{(1,1,1,1,1)}+\zeta_{25}E_{(1,1,1,2,1)}
-\zeta_{27}E_{(1,1,1,2,2,1)}\\& &-2\zeta_{15}
c+\frac{\zeta_{15}}{3}(\al_1-\al_3+\al_5-\al_6)
 \hspace{7.4cm}(5.84)\end{eqnarray*}
by (2.10), (2.11), (3.4) and (3.43),
\begin{eqnarray*}T_{16}&=&[\iota(E_{-\al_3}),T_{13}]=
\zeta_{16}(D-8)-\chi\ptl_{x_{12}}+\zeta_1E'_{(1,1,2,2,2,1)}+\zeta_3E'_{(0,1,1,2,2,1)}\\&
&+\zeta_5E'_{(0,1,1,1,1,1)}-\zeta_7E'_{(0,0,1,1,1,1)}+\zeta_8E'_{(0,1,1,1,1)}
-\zeta_{10}E'_{(0,0,1,1,1)}-\zeta_{11}E'_{(0,0,0,0,1,1)}\\
&&+\zeta_{13}E_{-\al_3}-\zeta_{15}E_{-\al_5}-\zeta_{19}E_{\al_1}
+\zeta_{18}E_{\al_4}-\zeta_{20}E_{(0,1,0,1)}\\&
&-\zeta_{22}E_{(1,0,1,1)}+\zeta_{24}E_{(1,1,1,1)}-\zeta_{26}E_{(1,1,1,2,1)}
-\zeta_{27}E_{(1,1,1,2,1,1)}\\& &-2\zeta_{16}
c+\frac{\zeta_{16}}{3}(\al_1-\al_3-2\al_5-\al_6)
 \hspace{7cm}(5.85)\end{eqnarray*}
by (2.10), (2.11), (3.4) and (3.43),
\begin{eqnarray*}T_{17}&=&[\iota(E_{-\al_1}),T_{15}]=
\zeta_{17}(D-8)-\chi\ptl_{x_{11}}+\zeta_2E'_{(1,1,2,2,1,1)}-\zeta_3E'_{(1,1,1,2,1,1)}\\&
&-\zeta_4E'_{(1,1,1,1,1,1)}+\zeta_6E'_{(1,0,1,1,1,1)}+\zeta_8E'_{(1,1,1,1)}
-\zeta_{10}E'_{(1,0,1,1)}-\zeta_{12}E'_{(1,0,1)}\\
&&+\zeta_{14}E_{-\al_6}+\zeta_{15}E_{-\al_1}+\zeta_{19}E_{\al_5}
+\zeta_{21}E_{(0,0,0,1,1)}-\zeta_{23}E_{(0,1,0,1,1)}\\&
&-\zeta_{22}E_{(0,0,1,1,1)}+\zeta_{24}E_{(0,1,1,1,1)}+\zeta_{25}E_{(0,1,1,2,1)}
-\zeta_{27}E_{(0,1,1,2,2,1)}\\& &-2\zeta_{17}
c-\frac{\zeta_{17}}{3}(2\al_1+\al_3-\al_5+\al_6)
 \hspace{7.4cm}(5.86)\end{eqnarray*}
by (2.10), (2.11), (3.2) and (3.41),
\begin{eqnarray*}T_{18}&=&-[\iota(E_{-\al_4}),T_{16}]=
\zeta_{18}(D-8)-\chi\ptl_{x_{10}}-\zeta_1E'_{(1,1,2,3,2,1)}+\zeta_4E'_{(0,1,1,2,2,1)}\\&
&-\zeta_5E'_{(0,1,1,2,1,1)}-\zeta_8E'_{(0,1,1,2,1)}-\zeta_9E'_{(0,0,1,1,1,1)}
-\zeta_{11}E'_{(0,0,0,1,1,1)}-\zeta_{12}E'_{(0,0,1,1,1)}\\
&&-\zeta_{13}E'_{(0,0,1,1)}-\zeta_{15}E'_{(0,0,0,1,1)}-\zeta_{16}E_{-\al_4}
-\zeta_{20}E_{\al_2}-\zeta_{21}E_{\al_1}\\&
&-\zeta_{22}E_{(1,0,1)}-\zeta_{25}E_{(1,1,1,1)}-\zeta_{26}E_{(1,1,1,1,1)}
-\zeta_{27}E_{(1,1,1,1,1,1)}\\& &-2\zeta_{16}
c+\frac{\zeta_{18}}{3}(\al_1-\al_3-3\al_4-2\al_5-\al_6)
 \hspace{6cm}(5.87)\end{eqnarray*}
by (2.10), (2.11), (3.5) and (3.44),
\begin{eqnarray*}T_{19}&=&[\iota(E_{-\al_1}),T_{16}]=
\zeta_{19}(D-8)-\chi\ptl_{x_9}-\zeta_2E'_{(1,1,2,2,2,1)}-\zeta_3E'_{(1,1,1,2,2,1)}\\&
&-\zeta_5E'_{(1,1,1,1,1,1)}+\zeta_7E'_{(1,0,1,1,1,1)}-\zeta_8E'_{(1,1,1,1,1)}
+\zeta_{10}E'_{(1,0,1,1,1)}-\zeta_{13}E'_{(1,0,1)}\\
&&+\zeta_{14}E'_{(0,0,0,0,1,1)}+\zeta_{16}E_{-\al_1}-\zeta_{17}E_{-\al_5}
+\zeta_{21}E_{\al_4}-\zeta_{23}E_{(0,1,0,1)}\\&
&-\zeta_{22}E_{(0,0,1,1)}+\zeta_{24}E_{(0,1,1,1)}-\zeta_{26}E_{(0,1,1,2,1)}
-\zeta_{27}E_{(0,1,1,2,1,1)}\\& &-2\zeta_{19}
c-\frac{\zeta_{19}}{3}(2\al_1+\al_3+2\al_5+\al_6)
 \hspace{7cm}(5.88)\end{eqnarray*}
by (2.10), (2.11), (3.2) and (3.41),
\begin{eqnarray*}T_{20}&=&[\iota(E_{-\al_2}),T_{18}]=
\zeta_{20}(D-8)-\chi\ptl_{x_8}-\zeta_1E'_{(1,2,2,3,2,1)}-\zeta_6E'_{(0,1,1,2,2,1)}\\&
&+\zeta_7E'_{(0,1,1,2,1,1)}+\zeta_{10}E'_{(0,1,1,2,1)}+\zeta_9E'_{(0,1,1,1,1,1)}
+\zeta_{11}E'_{(0,1,0,1,1,1)}+\zeta_{12}E'_{(0,1,1,1,1)}\\
&&+\zeta_{13}E'_{(0,1,1,1)}+\zeta_{15}E'_{(0,1,0,1,1)}+\zeta_{16}E'_{(0,1,0,1)}
+\zeta_{18}E_{-\al_2}-\zeta_{23}E_{\al_1}\\&
&-\zeta_{24}E_{(1,0,1)}-\zeta_{25}E_{(1,0,1,1)}-\zeta_{26}E_{(1,0,1,1,1)}
-\zeta_{27}E_{(1,0,1,1,1,1)}\\& &-2\zeta_{20}
c+\frac{\zeta_{20}}{3}(\al_1-3\al_2-\al_3-3\al_4-2\al_5-\al_6)
 \hspace{4.9cm}(5.89)\end{eqnarray*}
by (2.10), (2.11), (3.3) and (3.42),
\begin{eqnarray*}T_{21}&=&-[\iota(E_{-\al_4}),T_{19}]=
\zeta_{21}(D-8)-\chi\ptl_{x_7}+\zeta_2E'_{(1,1,2,3,2,1)}-\zeta_4E'_{(1,1,1,2,2,1)}\\&
&+\zeta_5E'_{(1,1,1,2,1,1)}+\zeta_8E'_{(1,1,1,2,1)}+\zeta_9E'_{(1,0,1,1,1,1)}
+\zeta_{12}E'_{(1,0,1,1,1)}+\zeta_{13}E'_{(1,0,1,1)}\\
&&+\zeta_{14}E'_{(0,0,0,1,1,1)}-\zeta_{17}E'_{(0,0,0,1,1)}+\zeta_{18}E_{-\al_1}
-\zeta_{19}E_{-\al_4}-\zeta_{22}E_{\al_3}\\&
&-\zeta_{23}E_{\al_2}-\zeta_{25}E_{(0,1,1,1)}-\zeta_{26}E_{(0,1,1,1,1)}
-\zeta_{27}E_{(0,1,1,1,1,1)}\\& &-2\zeta_{21}
c-\frac{\zeta_{21}}{3}(2\al_1+\al_3+3\al_4+2\al_5+\al_6)
 \hspace{5.9cm}(5.90)\end{eqnarray*}
by (2.10), (2.11), (3.5) and (3.44),
\begin{eqnarray*}T_{22}&=&[\iota(E_{-\al_3}),T_{21}]=
\zeta_{22}(D-8)-\chi\ptl_{x_6}-\zeta_3E'_{(1,1,2,3,2,1)}+\zeta_4E'_{(1,1,2,2,2,1)}\\&
&-\zeta_5E'_{(1,1,2,2,1,1)}-\zeta_8E'_{(1,1,2,2,1)}+\zeta_{11}E'_{(1,0,1,1,1,1)}
+\zeta_{15}E'_{(1,0,1,1,1)}+\zeta_{16}E'_{(1,0,1,1)}\hspace{5cm}\end{eqnarray*}\begin{eqnarray*}
\qquad
&&-\zeta_{14}E'_{(0,0,1,1,1,1)}+\zeta_{17}E'_{(0,0,1,1,1)}+\zeta_{18}E'_{(1,0,1)}
+\zeta_{19}E'_{(0,0,1,1)}+\zeta_{21}E_{-\al_3}\\&
&-\zeta_{24}E_{\al_2}-\zeta_{25}E_{(0,1,0,1)}-\zeta_{26}E_{(0,1,0,1,1)}
-\zeta_{27}E_{(0,1,0,1,1,1)}\\& &-2\zeta_{22}
c-\frac{\zeta_{22}}{3}(2\al_1+4\al_3+3\al_4+2\al_5+\al_6)
 \hspace{5.8cm}(5.91)\end{eqnarray*}
by (2.10), (2.11), (3.4) and (3.43),
\begin{eqnarray*}T_{23}&=&[\iota(E_{-\al_2}),T_{21}]=
\zeta_{23}(D-8)-\chi\ptl_{x_5}+\zeta_2E'_{(1,2,2,3,2,1)}+\zeta_6E'_{(1,1,1,2,2,1)}\\&
&-\zeta_7E'_{(1,1,1,2,1,1)}-\zeta_{10}E'_{(1,1,1,2,1)}-\zeta_9E'_{(1,1,1,1,1,1)}
-\zeta_{12}E'_{(1,1,1,1,1)}-\zeta_{13}E'_{(1,1,1,1)}\\
&&-\zeta_{14}E'_{(0,1,0,1,1,1)}+\zeta_{17}E'_{(0,1,0,1,1)}+\zeta_{19}E'_{(0,1,0,1)}+\zeta_{20}E_{-\al_1}
+\zeta_{21}E_{-\al_2}\\&
&-\zeta_{24}E_{\al_3}-\zeta_{25}E_{(0,0,1,1)}-\zeta_{26}E_{(0,0,1,1,1)}
-\zeta_{27}E_{(0,0,1,1,1,1)}\\& &-2\zeta_{23}
c-\frac{\zeta_{23}}{3}(2\al_1+3\al_2+\al_3+3\al_4+2\al_5+\al_6)
 \hspace{4.9cm}(5.92)\end{eqnarray*}
by (2.10), (2.11), (3.3) and (3.42),
\begin{eqnarray*}T_{24}&=&[\iota(E_{-\al_3}),T_{23}]=
\zeta_{24}(D-8)-\chi\ptl_{x_4}-\zeta_3E'_{(1,2,2,3,2,1)}-\zeta_6E'_{(1,1,2,2,2,1)}\\&
&+\zeta_7E'_{(1,1,2,2,1,1)}+\zeta_{10}E'_{(1,1,2,2,1)}-\zeta_{11}E'_{(1,1,1,1,1,1)}
+\zeta_{14}E'_{(0,1,1,1,1,1)}-\zeta_{15}E'_{(1,1,1,1,1)}\\
&&-\zeta_{16}E'_{(1,1,1,1)}-\zeta_{17}E'_{(0,1,1,1,1)}
-\zeta_{19}E'_{(0,1,1,1)}+\zeta_{20}E'_{(1,0,1)}+\zeta_{22}E_{-\al_2}\\&
&+\zeta_{23}E_{-\al_3}-\zeta_{25}E_{\al_4}-\zeta_{26}E_{(0,0,0,1,1)}
-\zeta_{27}E_{(0,0,0,1,1,1)}\\& &-2\zeta_{24}
c-\frac{\zeta_{24}}{3}(2\al_1+3\al_2+4\al_3+3\al_4+2\al_5+\al_6)
 \hspace{4.6cm}(5.93)\end{eqnarray*}
by (2.10), (2.11), (3.4) and (3.43),
\begin{eqnarray*}T_{25}&=&[\iota(E_{-\al_4}),T_{24}]=
\zeta_{25}(D-8)-\chi\ptl_{x_3}+\zeta_4E'_{(1,2,2,3,2,1)}-\zeta_6E'_{(1,1,2,3,2,1)}\\&
&-\zeta_9E'_{(1,1,2,2,1,1)}-\zeta_{11}E'_{(1,1,1,2,1,1)}-\zeta_{12}E'_{(1,1,2,2,1)}
+\zeta_{14}E'_{(0,1,1,2,1,1)}-\zeta_{15}E'_{(1,1,1,2,1)}\\
&&+\zeta_{18}E'_{(1,1,1,1)}-\zeta_{17}E'_{(0,1,1,2,1)}
+\zeta_{20}E'_{(1,0,1,1)}+\zeta_{21}E'_{(0,1,1,1)}+\zeta_{22}E'_{(0,1,0,1)}\\&
&+\zeta_{23}E'_{(0,0,1,1)}+\zeta_{24}E_{-\al_4}-\zeta_{26}E_{\al_5}
-\zeta_{27}E_{(0,0,0,0,1,1)}\\& &-2\zeta_{25}
c-\frac{\zeta_{25}}{3}(2\al_1+3\al_2+4\al_3+6\al_4+2\al_5+\al_6)
 \hspace{4.6cm}(5.94)\end{eqnarray*}
by (2.10), (2.11), (3.5) and (3.44),
\begin{eqnarray*}T_{26}&=&[\iota(E_{-\al_5}),T_{25}]=
\zeta_{26}(D-8)-\chi\ptl_{x_2}-\zeta_5E'_{(1,2,2,3,2,1)}+\zeta_7E'_{(1,1,2,3,2,1)}\\&
&-\zeta_9E'_{(1,1,2,2,2,1)}-\zeta_{11}E'_{(1,1,1,2,2,1)}+\zeta_{13}E'_{(1,1,2,2,1)}
+\zeta_{14}E'_{(0,1,1,2,2,1)}\\&
&+\zeta_{16}E'_{(1,1,1,2,1)}+\zeta_{18}E'_{(1,1,1,1,1)}+\zeta_{19}E'_{(0,1,1,2,1)}
+\zeta_{20}E'_{(1,0,1,1,1)}\\&
&+\zeta_{21}E'_{(0,1,1,1,1)}+\zeta_{22}E'_{(0,1,0,1,1)}+\zeta_{23}E'_{(0,0,1,1,1)}+\zeta_{24}E'_{(0,0,0,1,1)}+\zeta_{25}E_{-\al_5}
\\& &-\zeta_{27}E_{\al_6}-2\zeta_{26}
c-\frac{\zeta_{26}}{3}(2\al_1+3\al_2+4\al_3+6\al_4+5\al_5+\al_6)
 \hspace{2.8cm}(5.95)\end{eqnarray*}
by (2.10), (2.11), (3.6) and (3.45),
\begin{eqnarray*}T_{27}&=&[\iota(E_{-\al_6}),T_{26}]=
\zeta_{27}(D-8)-\chi\ptl_{x_1}+\zeta_8E'_{(1,2,2,3,2,1)}-\zeta_{10}E'_{(1,1,2,3,2,1)}\\&
&+\zeta_{12}E'_{(1,1,2,2,2,1)}+\zeta_{13}E'_{(1,1,2,2,1,1)}+\zeta_{15}E'_{(1,1,1,2,2,1)}
+\zeta_{17}E'_{(0,1,1,2,2,1)}\\&
&+\zeta_{16}E'_{(1,1,1,2,1,1)}+\zeta_{18}E'_{(1,1,1,1,1,1)}+\zeta_{19}E'_{(0,1,1,2,1,1)}
+\zeta_{20}E'_{(1,0,1,1,1,1)}\\&
&+\zeta_{21}E'_{(0,1,1,1,1,1)}+\zeta_{22}E'_{(0,1,0,1,1,1)}+\zeta_{23}E'_{(0,0,1,1,1,1)}+\zeta_{24}E'_{(0,0,0,1,1,1)}+\zeta_{26}E_{-\al_6}
\\& &+\zeta_{25}E'_{(0,0,0,0,1,1)}-2\zeta_{27}
c-\frac{\zeta_{27}}{3}(2\al_1+3\al_2+4\al_3+6\al_4+5\al_5+4\al_6)
 \hspace{1.2cm}(5.96)\end{eqnarray*}
by (2.10), (2.11), (3.7) and (3.46). Then ${\cal
T}=\sum_{r=1}^{27}\mbb{C}T_r$ forms the 27-dimensional ${\cal
G}^{E_6}$ of highest weight $\lmd_1$. Indeed, the map $\eta_r\mapsto
T_r$ for $r\in\ol{1,27}$ determines a ${\cal G}^{E_6}$-module
isomorphism from  $U$ to ${\cal T}$ (cf. (3.34)).

Denote
$$T'_r=T_r-\zeta_r(D-2c-8)+\chi\ptl_{x_{28-r}}\qquad\for\;\;r\in\ol{1,27}.\eqno(5.97)$$
Easily see that  ${\cal T}'=\sum_{r=1}^{27}\mbb{C}T_r'$ forms the
27-dimensional ${\cal G}^{E_6}$-module of highest weight $\lmd_1$.
So we have the ${\cal G}^{E_6}$-module isomorphism from
$U=\sum_{r=1}^{27}\mbb{C}\zeta_r$ to ${\cal T}'$ determined by
$\zeta_r\mapsto T_r'$ for $r\in\ol{1,27}$. The weight set of $U$ is
$$\Pi(U)=\{\sum_{s=1}^6b_{r,s}\lmd_s\mid r\in\ol{1,27}\}\eqno(5.98)$$
(cf. (5.9)-(5.14) and Table 2 in Section 3). Let $\lmd\in\Lmd^+$.
Denote
$$\Upsilon'(\lmd)=\{\lmd+\mu\mid
\mu\in\Pi(U),\;\;\lmd+\mu\in\Lmd^+\}\eqno(5.99)$$ (cf. (5.27)). Take
$M=V(\lmd)$, the irreducible ${\cal G}^{E_6}$-module of highest
weight $\lmd$. It is known that
$$UV(\lmd)=U\otimes_{\mbb{C}}V(\lmd)\cong
\bigoplus_{\lmd'\in \Upsilon'(\lmd)}V(\lmd').\eqno(5.100)$$

Given $\lmd'\in\Upsilon'(\lmd)$, we pick a singular vector
$$u=\sum_{r=1}^{27}\zeta_ru_r\eqno(5.101)$$
of weight $\lmd'$ in $UV(\lmd)$, where $u_r\in V(\lmd)$. Moreover,
any singular vector of  weight $\lmd'$ in $UV(\lmd)$ is a scalar
multiple of $u$. Note that the vector
$$w=\sum_{r=1}^{27}T'_r(u_r)\eqno(5.102)$$
is also a  singular vector of  weight $\lmd'$ if it is not zero.
Thus
$$w=\flat_{\lmd'}u,\qquad \flat_{\lmd'}\in\mbb{C}.\eqno(5.103)$$
Set
$$\flat(\lmd)=\min\{\flat_{\lmd'}\mid\lmd'\in\Upsilon'(\lmd)\}.\eqno(5.104)$$
\psp

{\bf Theorem 5.6}. {\it The ${\cal G}^{E_7}$-module
$\widehat{V(\lmd)}$ is irreducible if}
$$c\in\mbb{C}\setminus\{-8+5\mbb{N}/3,
(1/2)(\flat(\lmd)+\mbb{N})-4,\ell_\omega(\lmd)+\mbb{N}\}.\eqno(5.105)$$

{\it Proof}. Recall that the ${\cal G}^{E_7}$-submodule $U({\cal
G}_-)(V(\lmd))$ is irreducible by Proposition 4.2. It is enough to
prove $\widehat{V(\lmd)}=U({\cal G}_-)(V(\lmd))$. It is obvious that
$$\widehat{V(\lmd)}_{\la 0\ra}=V(\lmd)=(U({\cal G}_-)(V(\lmd)))_{\la
0\ra}\eqno(5.106)$$ (cf. (5.3) and (5.4) with $M=V(\lmd)$).
Moreover, Lemma 5.1 with $M=V(\lmd)$,  (5.37) and (5.105) imply that
$\vf|_{\widehat{V(\lmd)}_{\la 1\ra}}$ is invertible, equivalently,
$$\widehat{V(\lmd)}_{\la 1\ra}=(U({\cal G}_-)(V(\lmd)))_{\la
1\ra}.\eqno(5.107)$$ Suppose that
$$\widehat{V(\lmd)}_{\la i\ra}=(U({\cal G}_-)(V(\lmd)))_{\la
i\ra}\eqno(5.108)$$ for $i\in\ol{0,k}$ with $1\leq k\in\mbb{N}$.

For any $v\in V(\lmd)$ and $\al\in\mbb{N}^{16}$ such that
$|\al|=k-1$, we have
$$T_r(x^\al v)=x^\al[(|\al|-2c-8)\zeta_r+T'_r](v)+\chi \ptl_{x_{28-r}}(x^\al)v\in (U({\cal G}_-)(V(\lmd)))_{\la
k+1\ra}\eqno(5.109)$$ for $r\in\ol{1,27}$ by (5.97), (5.108) with
$i=k-1,k$. If $k=1$, then $\al=0$. So $\chi
\ptl_{x_{28-r}}(x^\al)v=0$. When $k>1$,
$$\ptl_{x_{28-r}}(x^\al)v\in
\widehat{V(\lmd)}_{\la k-2\ra}=(U({\cal G}_-)(V(\lmd)))_{\la
k-2\ra}\eqno(5.110)$$ and so \begin{eqnarray*}\qquad &
&[(\sum_{14\neq
p\in\ol{1,26}}\iota(\eta_p)\zeta_{28-p}-\iota(\eta_{14})\zeta_{14}-\iota(\eta_{27})\zeta_1](\ptl_{x_{28-r}}(x^\al)v)
\\&=&\chi(24-5(|\al|-1)+3c)(\ptl_{x_{28-r}}(x^\al)v)\in (U({\cal G}_-)(V(\lmd)))_{\la
k+1\ra}\hspace{2.4cm}(5.111)\end{eqnarray*} by Lemma 5.5. Thus
(5.105) gives
$$\chi\ptl_{x_{28-r}}(x^\al)v\in (U({\cal G}_-)(V(\lmd)))_{\la
k+1\ra}\qquad\for\;\;r\in\ol{1,27}.\eqno(5.112)$$ Hence in any case,
$$T_r(x^\al v)=x^\al[(|\al|-2c-8)\zeta_r+T'_r](v)\in (U({\cal G}_-)(V(\lmd)))_{\la
k+1\ra}\qquad\for\;\;r\in\ol{1,27}.\eqno(5.113)$$

On the other hand,
$$V'=\mbox{Span}\{[(|\al|-2c-8)\zeta_r+T'_r](v)\mid r\in\ol{1,27},\;v\in
V(\lmd)\}\eqno(5.114)$$ forms a ${\cal G}^{E_6}$-submodule of
$UV(\lmd)$ with respect to the action in (5.1). Let $u$ be a ${\cal
G}^{E_6}$-singular vector in (5.101). Then
$$V'\ni
\sum_{r=1}^{27}[(|\al|-2c-8)\zeta_r+T'_r](u_r)=(|\al|-2c-8)u+w=(|\al|-2c-8+\flat_{\lmd'})u\eqno(5.115)$$
by (5.102) and (5.103). Moreover, (5.104) and (5.105) yield $u\in
V'$. Since $UV(\lmd)$ is a ${\cal G}^{E_6}$-module generated by all
the singular vectors, we have $V'=UV(\lmd)$. So
$$x^\al UV(\lmd)\subset (U({\cal G}_-)(V(\lmd)))_{\la
k+1\ra}.\eqno(5.116)$$ The arbitrariness of $\al$ implies
$$\zeta_r\widehat{V(\lmd)}_{\la k-1\ra}\subset (U({\cal G}_-)(V(\lmd)))_{\la
k+1\ra}\qquad\for\;\;r\in\ol{1,27}.\eqno(5.117)$$

Given any $f\in{\cal A}_k$ and $v\in V(\lmd)$, we have
$$\zeta_r\ptl_{x_s}(f)v\in \zeta_r\widehat{V(\lmd)}_{\la k-1\ra}\subset (U({\cal G}_-)(V(\lmd)))_{\la
k+1\ra}\qquad\for\;\;r,s\in\ol{1,27}.\eqno(5.118)$$ Moreover,
\begin{eqnarray*}\qquad\eta_s(fv)&=&\iota(\eta_s)(fv)=P_s(fv)+f(\td\omega-c)(x_sv)\\
&\equiv&
f(k+\td\omega-c)(x_sv)\;\;(\mbox{mod}\;\sum_{r=1}^{27}\zeta_r\widehat{V(\lmd)}_{\la
k-1\ra})\hspace{4.3cm}(5.119)\end{eqnarray*} for $s\in\ol{1,27}$ by
(3.58), (3.77)-(3.102), (4.13)-(4.39) and Lemma 5.1. According to
(5.37), (5.105), (5.115) and (5.117), we get
$$x_sfv\in (U({\cal G}_-)(V(\lmd)))_{\la
k+1\ra}\qquad\for\;\;s\in\ol{1,27}.\eqno(5.120)$$ Thus (5.108) holds
for $i=k+1$. By induction on $k$, (5.108) holds for any
$i\in\mbb{N}$, that is, $\widehat{V(\lmd)}=U({\cal
G}_-)(V(\lmd)).\qquad\Box$\psp

When $\lmd=0$, $V(0)$ is the one-dimensional trivial module and
$\ell_\omega(0)=\flat(0)=0$. So we have:\psp

{\bf Corollary 6.7}. {\it The ${\cal G}^{E_7}$-module
$\widehat{V(0)}$ is irreducible if}
$c\in\mbb{C}\setminus\{5\mbb{N}/3-8, \mbb{N}/2-4\}.$ \psp

Next we consider $\lmd=k\lmd_1$. In this case,
$$\Upsilon(\lmd)=\{k\lmd_1+\lmd_6,(k-1)\lmd_1+\lmd_2,(k-1)\lmd_1\}
\eqno(5.121)$$ by (5.32) and Tables 1, 2. Thus we have
$$\ell_\omega(\lmd)=-16-\frac{4k}{3}\eqno(5.122)$$ by (5.37).
Moreover,
$$\Upsilon'(\lmd)=\{(k+1)\lmd_1,(k-1)\lmd_1+\lmd_3,(k-1)\lmd_1+\lmd_6\}\eqno(5.123)$$
by Table 2 and (5.99).

We define a representation of ${\cal G}^{E_6}$ on ${\cal
C}=\mbb{C}[z_1,...,z_{27}]$ determined via (3.35)-(3.47) with $U$
replaced by ${\cal C}$ and $\zeta_i$ replaced by $z_i$. In fact,
$$E_{\al_1}|_{\cal C}=z_1\ptl_{z_2}+z_{11}\ptl_{z_{14}}-z_{15}\ptl_{z_{17}}
-z_{16}\ptl_{z_{19}}-z_{18}\ptl_{z_{21}}-z_{20}\ptl_{z_{23}},\eqno(5.124)$$
$$E_{\al_2}|_{\cal C}=z_4\ptl_{z_6}+z_5\ptl_{z_7}+z_8\ptl_{z_{10}}
-z_{18}\ptl_{z_{20}}-z_{21}\ptl_{z_{23}}-z_{22}\ptl_{z_{24}},\eqno(5.125)$$
$$E_{\al_3}|_{\cal C}=z_2\ptl_{z_3}-z_9\ptl_{z_{11}}
-z_{12}\ptl_{z_{15}}-z_{13}\ptl_{z_{16}}-z_{21}\ptl_{z_{22}}-z_{23}\ptl_{z_{24}},\eqno(5.126)$$
$$E_{\al_4}|_{\cal C}=z_3\ptl_{z_4}+z_7\ptl_{z_9}+z_{10}\ptl_{z_{12}}
+z_{16}\ptl_{z_{18}}+z_{19}\ptl_{z_{21}}-z_{24}\ptl_{z_{25}},\eqno(5.127)$$
$$E_{\al_5}|_{\cal C}=z_4\ptl_{z_5}+z_6\ptl_{z_7}+z_{12}\ptl_{z_{13}}+z_{15}\ptl_{z_{16}}
+z_{17}\ptl_{z_{19}}-z_{25}\ptl_{z_{26}},\eqno(5.128)$$
$$E_{\al_6}|_{\cal C}=z_5\ptl_{z_8}+z_7\ptl_{z_{10}}+z_9\ptl_{z_{12}}+z_{11}
\ptl_{z_{15}}
-z_{14}\ptl_{z_{17}}-z_{26}\ptl_{z_{27}},\eqno(5.129)$$
$$E_{(1,0,1)}|_{\cal C}=-z_1\ptl_{z_3}-z_9\ptl_{z_{14}}+z_{12}\ptl_{z_{17}}
+z_{13}\ptl_{z_{19}}-z_{18}\ptl_{z_{22}}-z_{20}\ptl_{z_{24}},\eqno(5.130)$$
$$E_{(0,1,0,1)}|_{\cal C}=z_3\ptl_{z_6}-z_5\ptl_{z_9}-z_8\ptl_{z_{12}}-z_{16}\ptl_{z_{20}}
-z_{19}\ptl_{z_{23}}-z_{22}\ptl_{z_{25}},\eqno(5.131)$$
$$E_{(0,0,1,1)}|_{\cal C}=-z_2\ptl_{z_4}-z_7\ptl_{z_{11}}-z_{10}\ptl_{z_{15}}+z_{13}\ptl_{z_{18}}
-z_{19}\ptl_{z_{22}}-z_{23}\ptl_{z_{25}},\eqno(5.132)$$
$$E_{(0,0,0,1,1)}|_{\cal C}=-z_3\ptl_{z_5}+z_6\ptl_{z_9}-z_{10}\ptl_{z_{13}}+z_{15}\ptl_{z_{18}}
+z_{17}\ptl_{z_{21}}-z_{24}\ptl_{z_{26}},\eqno(5.133)$$
$$E_{(0,0,0,0,1,1)}|_{\cal C}=-z_4\ptl_{z_8}-z_6\ptl_{z_{10}}+z_9\ptl_{z_{13}}+z_{11}\ptl_{z_{16}}-z_{14}\ptl_{z_{19}}
+z_{25}\ptl_{z_{27}},\eqno(5.134)$$
$$E_{(1,0,1,1)}|_{\cal C}=z_1\ptl_{z_4}-z_7\ptl_{z_{14}}+z_{10}\ptl_{z_{17}}-z_{13}\ptl_{z_{21}}-z_{16}\ptl_{z_{22}}-z_{20}\ptl_{z_{25}},
\eqno(5.135)$$
$$E_{(0,1,1,1)}|_{\cal C}=-z_2\ptl_{z_6}+z_5\ptl_{z_{11}}+z_8\ptl_{z_{15}}-z_{13}\ptl_{z_{20}}
+z_{19}\ptl_{z_{z_{24}}}-z_{21}\ptl_{z_{25}},\eqno(5.136)$$
$$E_{(0,1,0,1,1)}|_{\cal C}=-z_3\ptl_{z_7}-
z_4\ptl_{z_9}+z_8\ptl_{z_{13}}-z_{15}\ptl_{z_{20}}-z_{17}\ptl_{z_{23}}-z_{22}\ptl_{z_{26}},\eqno(5.137)$$
$$E_{(0,0,1,1,1)}|_{\cal C}=z_2\ptl_{z_5}-z_6\ptl_{z_{11}}+z_{10}\ptl_{z_{16}}+z_{12}\ptl_{z_{18}}
-z_{17}\ptl_{z_{22}}-z_{23}\ptl_{z_{26}},\eqno(5.138)$$
$$E_{(0,0,0,1,1,1)}|_{\cal A}=z_3\ptl_{z_8}-z_6\ptl_{z_{12}}-z_7\ptl_{z_{13}}+z_{11}\ptl_{z_{18}}
-z_{14}\ptl_{z_{21}} +z_{24}\ptl_{z_{27}},\eqno(5.139)$$
$$E_{(1,1,1,1)}|_{\cal C}=z_1\ptl_{z_6}+z_5\ptl_{z_{14}}
-z_8\ptl_{z_{17}}+z_{13}\ptl_{z_{23}}+z_{16}\ptl_{z_{24}}-z_{18}\ptl_{z_{25}},\eqno(5.140)$$
$$E_{(1,0,1,1,1)}|_{\cal C}=-z_1\ptl_{z_5}-z_6\ptl_{z_{14}}-z_{10}\ptl_{z_{19}}-z_{12}\ptl_{z_{21}}-z_{15}\ptl_{z_{22}}
-z_{20}\ptl_{z_{26}},\eqno(5.141)$$
$$E_{(0,1,1,1,1)}|_{\cal C}=z_2\ptl_{z_7}+z_4\ptl_{z_{11}}-z_8\ptl_{z_{16}}-z_{12}\ptl_{z_{20}}
+z_{17}\ptl_{z_{24}}-z_{21}\ptl_{z_{26}},\eqno(5.142)$$
$$E_{(0,1,0,1,1,1)}|_{\cal C}=z_3\ptl_{z_{10}}+z_4\ptl_{z_{12}}+z_5\ptl_{z_{15}}-z_{11}\ptl_{z_{20}}+z_{14}\ptl_{z_{23}}
+z_{22}\ptl_{z_{27}},\eqno(5.143)$$
$$E_{(0,0,1,1,1,1)}|_{\cal C}=-z_2\ptl_{z_8}+z_6\ptl_{z_{15}}+z_7\ptl_{z_{16}}+z_9\ptl_{z_{18}}+z_{14}\ptl_{z_{22}}+z_{23}\ptl_{z_{27}},
\eqno(5.144)$$
$$E_{(1,1,1,1,1}|_{\cal C}=-z_1\ptl_{z_7}+z_4\ptl_{z_{14}}+z_8\ptl_{z_{19}}+z_{12}\ptl_{z_{23}}+z_{15}\ptl_{z_{24}}
-z_{18}\ptl_{z_{26}},\eqno(5.145)$$
$$E_{(1,0,1,1,1,1)}|_{\cal C}=z_1\ptl_{z_8}-z_6\ptl_{z_{17}}-z_7\ptl_{z_{19}}-z_9\ptl_{z_{21}}-z_{11}\ptl_{z_{22}}
+z_{20}\ptl_{z_{27}},\eqno(5.146)$$
$$E_{(0,1,1,2,1)}|_{\cal C}=-z_2\ptl_{z_9}+z_3\ptl_{z_{11}}+z_8\ptl_{z_{18}}-z_{10}\ptl_{z_{20}}
+z_{17}\ptl_{z_{25}}+z_{19}\ptl_{z_{26}},\eqno(5.147)$$
$$E_{(0,1,1,1,1,1)}|_{\cal C}=-z_2\ptl_{z_{10}}-z_4\ptl_{z_{15}}-z_5\ptl_{z_{16}}-z_9\ptl_{z_{20}}-z_{14}\ptl_{z_{24}}+z_{21}\ptl_{z_{27}},
\eqno(5.148)$$
$$E_{(1,1,1,2,1)}|_{\cal C}=z_1\ptl_{z_9}+z_3\ptl_{z_{14}}-z_8\ptl_{z_{21}}+z_{10}\ptl_{z_{23}}
+z_{15}\ptl_{z_{25}}+z_{16}\ptl_{z_{26}},\eqno(5.149)$$
$$E_{(1,1,1,1,1,1)}|_{\cal C}=z_1\ptl_{z_{10}}
+z_4\ptl_{z_{17}}+z_5\ptl_{z_{19}}+z_9\ptl_{z_{23}}+z_{11}\ptl_{z_{24}}+z_{18}\ptl_{z_{27}},
\eqno(5.150)$$
$$E_{(0,1,1,2,1,1)}|_{\cal C}=z_2\ptl_{z_{12}}-z_3\ptl_{z_{15}}
+z_5\ptl_{z_{18}}-z_7\ptl_{z_{20}}-z_{14}\ptl_{z_{25}}+z_{19}\ptl_{z_{27}},
\eqno(5.151)$$
$$E_{(1,1,2,2,1)}|_{\cal C}=-z_1\ptl_{z_{11}}
-z_2\ptl_{z_{14}}+z_8\ptl_{z_{22}}-z_{10}\ptl_{z_{24}}
+z_{12}\ptl_{z_{25}}+z_{13}\ptl_{z_{26}},\eqno(5.152)$$
$$E_{(1,1,1,2,1,1)}|_{\cal C}=-z_1\ptl_{z_{12}}
+z_3\ptl_{z_{17}}-z_5\ptl_{z_{21}}+z_7\ptl_{z_{23}}+z_{11}\ptl_{z_{25}}+z_{16}\ptl_{z_{27}},
\eqno(5.153)$$
$$E_{(0,1,1,2,2,1)}|_{\cal C}=-z_2\ptl_{z_{13}}+z_3\ptl_{z_{16}}
+z_4\ptl_{z_{18}}-z_6\ptl_{z_{20}}-z_{14}\ptl_{z_{26}}+z_{17}\ptl_{z_{27}},
\eqno(5.154)$$
$$E_{(1,1,2,2,1,1)}|_{\cal C}=z_1\ptl_{z_{15}}
-z_2\ptl_{z_{17}}+z_5\ptl_{z_{22}}-z_7\ptl_{z_{24}}
+z_9\ptl_{z_{25}}+z_{13}\ptl_{z_{27}},\eqno(5.155)$$
$$E_{(1,1,1,2,2,1)}|_{\cal C}=z_1\ptl_{z_{13}}
-z_3\ptl_{z_{19}}-z_4\ptl_{z_{21}}+z_6\ptl_{z_{23}}+z_{11}\ptl_{z_{26}}+z_{15}\ptl_{z_{27}},
\eqno(5.156)$$
$$E_{(1,1,2,2,2,1)}|_{\cal C}=-z_1\ptl_{z_{16}}
+z_2\ptl_{z_{19}}+z_4\ptl_{z_{22}}-z_6\ptl_{z_{24}}
+z_9\ptl_{z_{26}}+z_{12}\ptl_{z_{27}},\eqno(5.157)$$
$$E_{(1,1,2,3,2,1)}|_{\cal C}=z_1\ptl_{z_{18}}
-z_2\ptl_{z_{21}}+z_3\ptl_{z_{22}}-z_6\ptl_{z_{25}}+z_7\ptl_{z_{26}}+z_{10}\ptl_{z_{27}},
\eqno(5.158)$$
$$E_{(1,2,2,3,2,1)}|_{\cal C}=z_1\ptl_{z_{20}}
-z_2\ptl_{z_{23}}+z_3\ptl_{z_{24}}-z_4\ptl_{z_{25}}+z_5\ptl_{z_{26}}+z_8\ptl_{z_{27}},
\eqno(5.159)$$
$$\al_r|_{\cal C}=\sum_{i=1}^{27}b_{i,r}z_i\ptl_{z_i}\qquad\for\;\;r\in\ol{1,6}\eqno(5.160)$$
with $b_{i,r}$ given Table 2,
$$E_{-\al}|_{\cal
C}=-\tau'(E_{\al})\qquad\for\;\;\al\i\Phi_{E_6},\eqno(5.161)$$ where
$\tau'$ is a linear transformation on the differential operators in
$z_1,...,z_{27}$ such that
$$\tau'(z_i\ptl_{z_j})=z_j\ptl_{z_i}\qquad\for\;\;i,j\in\ol{1,27}.\eqno(5.162)$$
 Then the ${\cal G}^{E_6}$-submodule ${\cal N}_k$
generated by $z_1^k$ is isomorphic to $V(k\lmd_1)$.

We calculate:
$$z_1^{k-1}z_2=-\frac{1}{k}E_{-\al_1}(z_1^k)\in{\cal
N}_k,\;\;z_1^{k-1}z_3=-E_{-\al_3}(z_1^{k-1}z_2)\in{\cal
N}_k,\eqno(5.163)$$
$$z_1^{k-1}z_4=-E_{-\al_4}(z_1^{k-1}z_3)\in{\cal N}_k,\;\;z_1^{k-1}z_5=-E_{-\al_5}(z_1^{k-1}z_4)\in{\cal
N}_k,\eqno(5.164)$$
$$z_1^{k-1}z_6=-E_{-\al_2}(z_1^{k-1}z_4)\in{\cal N}_k,\;\;z_1^{k-1}z_7=-E_{-\al_2}(z_1^{k-1}z_5)\in{\cal
N}_k,\eqno(5.165)$$
$$z_1^{k-1}z_9=-E_{-\al_4}(z_1^{k-1}z_7)\in{\cal
N}_k,\;\;z_1^{k-1}z_{11}=E_{-\al_3}(z_1^{k-1}z_9)\in{\cal
N}_k,\eqno(5.166)$$
$$-E_{-\al_1}(z_1^{k-1}z_{11})=z_1^{k-1}z_{14}+(k-1)z_1^{k-2}z_2z_{11}\in{\cal
N}_k,\eqno(5.167)$$
$$-E_{-\al_3}(z_1^{k-1}z_{14}+(k-1)z_1^{k-2}z_2z_{11})=(k-1)z_3z_{11}\in{\cal
N}_k,\eqno(5.168)$$
$$E_{\al_3}((k-1)z_3z_{11})=(k-1)(z_2z_{11}-z_3z_9)\in{\cal
N}_k,\eqno(5.169)$$
$$E_{\al_4}E_{-\al_4}[(k-1)(z_2z_{11}-z_3z_9)]=(k-1)(z_3z_9+z_4z_7)\in{\cal
N}_k,\eqno(5.170)$$
$$-E_{\al_5}E_{-\al_5}[(k-1)(z_3z_9+z_4z_7)]=(k-1)(z_4z_7+z_5z_6)\in{\cal
N}_k.\eqno(5.171)$$

Now we take $M={\cal N}_k$ in our earlier settings. First
$\zeta_1z_1^k$ is a singular vector of weight $(k+1)\lmd_1$ in
$U{\cal N}_k$. By (5.69), (5.160) and Table 2,
$$T_1'(z_1^k)=\frac{\zeta_1}{3}(4\al_1+3\al_2+5\al_3+6\al_4+4\al_5+2\al_6)(z_1^k)=\frac{4k}{3}\zeta_1z_1^k.\eqno(5.172)$$
So $\flat_{(k+1)\lmd_1}=4k/3$. Next
$\zeta_1z_1^{k-1}z_2-\zeta_2z_1^k$ is a singular vector of weight
$(k-1)\lmd_1+\lmd_3$ in $U{\cal N}_k$. According to (5.69), (5.71),
(5.160) and Table 2:
\begin{eqnarray*}& &T'_1(z_1^{k-1}z_2)-T'_2(z_1^k)\\
&=&\zeta_2z_1^k+\frac{4k-3}{3}\zeta_1z_1^{k-1}z_2-k\zeta_1z_1^{k-1}z_2-\frac{k}{3}\zeta_2z_1^k
=\frac{k-3}{3}(\zeta_1z_1^{k-1}z_2-\zeta_2z_1^k),
\hspace{1.2cm}(5.173)\end{eqnarray*} which gives
$\flat_{(k-1)\lmd_1+\lmd_3}=k/3-1$.

Expressions (3.35)-(3.40), (5.124)-(5.129) and (5.163)-(5.171) show
that
\begin{eqnarray*}\qquad
u&=&\zeta_1[4z_1^{k-1}z_{14}+(k-1)z_1^{k-2}(z_2z_{11}+z_3z_9-z_4z_7+z_5z_6)]-(k+3)z_1^{k-1}\\&&\times[\zeta_2z_{11}+\zeta_3z_9
-\zeta_4z_7+\zeta_5z_6+\zeta_6z_5-\zeta_7z_4+\zeta_9z_3+\zeta_{11}z_2-\zeta_{14}z_1]\hspace{1.2cm}(5.174)\end{eqnarray*}
is a singular vector of weight $(k-1)\lmd_1+\lmd_6$ in $U{\cal
N}_k$. We find
\begin{eqnarray*}&
&T_1'[4z_1^{k-1}z_{14}+(k-1)z_1^{k-2}(z_2z_{11}+z_3z_9-z_4z_7+z_5z_6)]\\
&
&-(k+3)[T_2'(z_1^{k-1}z_{11})+T_3'(z_1^{k-1}z_9)-T_4'(z_1^{k-1}z_7)
+T_5'(z_1^{k-1}z_6)\\
&&+T_6'(z_1^{k-1}z_5)-T_7'(z_1^{k-1}z_4)+T_9'(z_1^{k-1}z_3)+T_{11}'(z_1^{k-1}z_2)-T_{14}'(z_1^k)]
\\
&=&(k+3)z_1^{k-1}(\zeta_2z_{11}+\zeta_3z_9-\zeta_4z_7+\zeta_5z_6+\zeta_6z_5-\zeta_7z_4+\zeta_9z_3+\zeta_{11}z_2)
\\ & &+(4k/3-2)\zeta_1[4z_1^{k-1}z_{14}+(k-1)z_1^{k-2}(z_2z_{11}+z_3z_9-z_4z_7+z_5z_6)]
\\ & &-(k+3)z_1^{k-2}[2z_1(\zeta_1z_{14}-\zeta_3z_9+\zeta_4z_7-\zeta_5z_6-\zeta_6z_5+\zeta_7z_4-\zeta_9z_3+\zeta_{14}z_1)\\
&
&+2z_1(\zeta_1z_{14}-\zeta_2z_{11}+\zeta_4z_7-\zeta_5z_6-\zeta_6z_5+\zeta_7z_4-\zeta_{11}z_2+\zeta_{14}z_1)\\
& &+2z_1(\zeta_1z_{14}-\zeta_2z_{11}
-\zeta_3z_9-\zeta_5z_6-\zeta_6z_5-\zeta_9z_3-\zeta_{11}z_2+\zeta_{14}z_1)\\
& &+2z_1(\zeta_1z_{14}-\zeta_2z_{11}
-\zeta_3z_9+\zeta_4z_7+\zeta_7z_4-\zeta_9z_3-\zeta_{11}z_2+\zeta_{14}z_1)
\\ & &+2(k-1)\zeta_1(z_2z_{11}+z_3z_9-z_4z_7+z_5z_6)-(2k/3+1)z_1(\zeta_2z_{11}+\zeta_3z_9
\\ &
&-\zeta_4z_7+\zeta_5z_6+\zeta_6z_5-\zeta_7z_4+\zeta_9z_3+\zeta_{11}z_2)+(2k/3)\zeta_{14}z_1^2]
\hspace{5cm}\end{eqnarray*}\begin{eqnarray*}  &
=&(k+3)z_1^{k-1}(\zeta_2z_{11}+\zeta_3z_9-\zeta_4z_7+\zeta_5z_6+\zeta_6z_5-\zeta_7z_4+\zeta_9z_3+\zeta_{11}z_2)
\\ & &+(4k/3-2)\zeta_1[4z_1^{k-1}z_{14}+(k-1)z_1^{k-2}(z_2z_{11}+z_3z_9-z_4z_7+z_5z_6)]
-(k+3)z_1^{k-2}\\ &
&\times[2(4\zeta_1z_1z_{14}+(k-1)\zeta_1(z_2z_{11}+z_3z_9-z_4z_7+z_5z_6))
+(2k/3+8)\zeta_{14}z_1^2 \\ & &-(2k/3+7)z_1(\zeta_2z_{11}+\zeta_3z_9
-\zeta_4z_7+\zeta_5z_6+\zeta_6z_5-\zeta_7z_4+\zeta_9z_3+\zeta_{11}z_2)]
\\ &=&-(2k/3+8)u
 \hspace{11.3cm}(5.175)\end{eqnarray*}
by (5.69), (5.71)-(5.96), (5.124)-(5.161) and Table 2. Thus
$\flat_{(k-1)\lmd_1+\lmd_6}=-(2k/3+8)$. Therefore,
$$\flat(k\lmd_1)=-(2k/3+8).\eqno(5.176)$$
\pse

{\bf Corollary 5.8}. {\it The ${\cal G}^{E_7}$-module
$\widehat{V(k\lmd_1)}$ is irreducible if}
$$c\in\mbb{C}\setminus\{-8+5\mbb{N}/3,
\mbb{N}/2-8-k/3,\mbb{N}-16-4k/3\}.\eqno(5.177)$$

\vspace{1cm}

\noindent{\Large \bf References}

\hspace{0.5cm}

\begin{description}

\item[{[A]}] J. Adams, {\it Lectures on Exceptional Lie Groups}, The
University of Chicago Press Ltd., London, 1996.

\item[{[BKR]}] L. Brink, S. Kim and P. Ramond, $E_{7(7)}$ on the light cone, {\it
J. High Energy Phys.} (2008), no. 6, 034, 18pp.

\item[{[Br1]}] R. Brown, A minimal representation for the Lie algebra of type
$E_7$, {\it Illiois J. Math.} {\bf 12} (1968), 190-200.

\item[{[Br2]}] R. Brown,  Groups of type $E_7$, {\it J. Reine Angew. Math.} {\bf
236} (1969), 79-102.

\item[{[BK]}] R. Brylinski and B. Kostant, Minimal representations of
$E_6,\;E_7$, and $E_8$ and the generalized Capelli identity, {\it
Proc. Nar. Acad. Sci. U.S.A.} {\bf 91}, no. 7., 2469-2472.

\item[{[Bl]}] L. Borsten, $E_{7(7)}$ invariant measures of entanglement, {\it
Fortschr. Phys.} {\bf 56} (2008), no. 7-9, 842-848.

\item[{[CY]}] Y. Choi and S. Yoon, Homology of the double and the triple loop
spaces of $E_7,E_7$ and $E_8$, {\it Manuscripta Math.} {\bf 103}
(2000), 101-116.

\item[{[Cv]}] V. Chernousov, The kernel of the Rost invariant, Serre's
conjecture II and the Hasse principle for quasi-split groups
$D_4,\;E_6,\;E_7$, {\it Math. Ann.} {\bf 326} (2003), no. 2,
297-330.

\item[{[Cb]}] B. Cooperstein, The fifty-six-dimensional module for $E_7$.I. A
four form for $E_7$, {\it J. Algebra} {\bf 173} (1995), 361-389.

\item[{[Co]}] O. Cvitanovi\'{c}, Negative dimensions and $E_7$ symmetry, {\it
Nuclear Phys. B} {\bf 188} (1981), no. 2, 373-396.

\item[{[DFT]}] R. D'Auria, S. Ferrara and M. Trigiante, $E_{7(7)}$ symmetry and
dual gauge algebras of M-theory on a twisted even-torus, {\it
Nuclear Phys. B} {\bf 732} (2006),  389-400.

\item[{[Dt]}] T. De Medts, A characterization of quadratic forms of type
$E_6,\;E_7$ and $E_8$, {\it J. Algebra} {\bf 252} (2002), no. 2,
394-410.

\item[{[Dl]}] L. Dickson, A class of groups in an arbitrary realm
connected with the configuration of the 27 lines on a cubic surface,
{\it J. Math.} {\bf 33} (1901), 145-123.

\item[{[Dd1]}] D. Dokovi\'{c}, Explicit Cayley triples in real forms of $E_7$,
{\it Pacific J. Math.} {\bf 191} (1999), 1-23.

\item[{[Dd2]}] D. Dokovi\'{c}, The closure daigram for nilpotent orbits of the
split real form of $E_7$, {\it Reprent. Theory} {\bf 5} (2001),
284-316.

\item[{[DF]}] M. Duff and S. Ferrara, $E_7$ and the tripartite entanglement of
seven qubits, {\it Phys. Rev. D} {\bf 76} (2007), no. 2, 025018,
7pp.

\item[{[Fau]}]J. Faulkner, A geometry for $E_7$, {\it Trans. Amer. Math. Soc.}
{\bf 167} (1972), 49-58.

\item[{[FFP]}] N. Fern\'{a}dez, W. Garcia Fuerres and A. Perelomov, Quantum
trogonometric Calogero-Sutherland model, irreducible characters and
Clebsch-Gordan sries for the exceptional algebra $E_7$, {\it J.
Math. Phys.} {\bf 46} (2005), no. 10, 103505, 27pp.

\item[{[Fer]}] J. Ferrar, On the classification of Freudenthal triple systems and
Lie algebras of type $E_7$, {\it J. Algebra} {\bf 62} (1980),
276-282.

\item[{[Gr]}] R. Garibaldi, Structurable algebras and groups of type $E_6$ and
$E_7$, {\it J. Algebra} {\bf 236} (2001), no. 2, 651-691.

\item[{[Gd]}] D. Ginzburg, On standard $L$-functions for $E_6$ and $E_7$, {\it
J. Reine Angew. Math.} {\bf 465} (1995), 101-131.

\item[{[GR]}] R. Griess and A. Ryba, Embeddings of $U_3(8),\;Sz(8)$ and the
Rudvalis group in algebraic groups of type $E_7$, {\it Invent.
Math.} {\bf 116} (1994), 215-141.

\item[{[HKKT]}] S. Han, J. Kim, I. Koh and Y. Tanii, Supersymmetrization of the
six-dimensional anomaly-free $E_6\times E_7\times U(1)$ theory with
Lorentz Chern-Simons term, {\it Phys. Lett. B} {\bf 177} (1986), no.
2, 167-170.

\item[{[H]}] J. E. Humphreys, {\it Introduction to Lie Algebras and Representation Theory},
 Springer-Verlag New York Inc., 1972.

\item[{[K]}] V. Kac, {\it Infinite Dimensional Lie Algebras,} Third Edition,
Cambridge University Press, 1990.

\item[{[KK]}] A. Kato and Y. Kitazawa, $E_7$-type modular invariant Wess-Zumino
theory and Gepner's string compactification, {\it Nuclear Phys. B}
{\bf 319} (1989), no. 2, 474-490.

\item[{[KR]}] P. Kleidman and A. Ryba, Kostant's conjecture holds for $E_7:
L_2(37)<E_7(C)$, {\it J. Algebra} {\bf 161} (1993), 535-540.

\item[{[KMR]}] P. Kleidman, U. Meierfrankenfeld and A. Ryba, $HS<E_7(5)$, {\it J.
London Math. Soc. (2)} {\bf 60} (1999), 95-107.

\item[{[KLN]}] A. Kono, J. Lin and O. Nishimura, Characterization of the mod 3
cohomology of $E_7$, {\it Proc. Amer. Math. Soc.} {\bf 131} (2003),
no. 10, 3289-3295.

\item[{[M]}] Z. Ma, The spectrum-dependent solutions of the Yang-Baxter
equation for quantum $E_6$ and $E_7$, {\it J. Phys. A} {\bf 23}
(1990), no. 23, 5513-5522.

\item[{[P]}] E. Plotkin, On the stability of the $K_1$-functor for Chevalley
groups of type $E_7$, {\it J. Algebra} {\bf 210} (1998), 67-85.

\item[{[R]}] P. Ramond, Is there an exceptional group in your future? $E_7$ and
the travails of the symmetry breaking, {\it Nuclear Phys. B} {\bf
126} (1977), no. 3, 509-524.

\item[{[Sj]}]J. Sekiguchi, Configurations of seven lines on the real projective
plane and the root system of type $E_7$, {\it J. Math. Soc. Japan}
{\bf 51} (1999), no. 4, 987-1013.

\item[{[Se]}]E. Shult, Embeddings and hyperplanes of the Lie incidence geometry
of type $E_{7,1}$, {\it J. Geom.} {\bf 59} (1997), 152-172.

\item[{[U]}] K. Ukai, $b$-functions of the prehomogeneous vector space arising
from a cuspidal character sheaf of $E_7$, {\it J. Algebra} {\bf 237}
(2001), 358-381.

\item[{[W]}] R. Weiss, Moufang quadangles of type $E_6$ and $E_7$,
{\it J. Reine Angew. Math.} {\bf 590} (2006), 189-226.

\item[{[X1]}] X. Xu, {\it Kac-Moody Algebras and Their Representations}, China
Science Press, 2007.

\item[{[X2]}] X. Xu, A cubic $E_6$-generalization of the classical theorem on harmonic
polynomials, {\it J. Lie Theory} {\bf 21} (2011), 145-164.

\item[{[X3]}] X. Xu, Representations of Lie algebras and coding
theory, {\it J. Lie Theory} {\bf 22} (2012), no. 3, 647-682.

\item[{[X4]}] X. Xu, A new functor from $D_5$-{\bf Mod} to $E_6$-{\bf
Mod}, {\it arXiv:1112.3792v1[math.RT]}.

\item[{[XZ]}] X. Xu and Y. Zhao, Generalized conformal representations
of orthogonal Lie algebras, {\it arXiv:1105.1254v1[math.RT].}

\item[{[ZX]}] Y. Zhao and X. Xu, Generalized  projective  representations
 for sl(n+1), {\it J. Algebra} {\bf 328} (2011), 132-154.

\end{description}

\end{document}